%% file: maincomb.tex
%%%%%%%%%%%%%%%%%%%%%%%%%%%%%%%%%%%%%%%%%%%%%%%%%%%%%%%%%%%%%%%%%%%%%%%%%%
%%      FILE PRINCIPALE PER ARTICOLO SUL PETTINE                                                           %%
%%%%%%%%%%%%%%%%%%%%%%%%%%%%%%%%%%%%%%%%%%%%%%%%%%%%%%%%%%%%%%%%%%%%%%%%%
%%                                                                     %%
%%           DEFINIZIONI E MACRO                                       %%
%%                                                                     %%
								       
%\magnification= \magstep1
\input coman.tex
\input amssym.def
\input pictex.tex

%\bozze
\riferimentifuturi
\numerazionedoppia
%\indice

\def\d{{\,\rm d}} 
\def\phi{{\varphi}}
\def\epsilon{{\varepsilon}}
\def\re{{\rm{\bf Re}\ }}
\def\im{{\rm{\bf Im}\ }}
\def\r1{{R_1(\xi,t_n)}}
\def\uo{{u_o(\xi)}}
\def\zo{{z_o(\xi)}}
\def\zt{{z(\xi,t)}}
\def\prob{{\Bbb P}}

\def\deg{{\rm deg}}

\def\asym#1{{\,{\buildrel {#1}\over\sim}\,}}

\def\eps{\varepsilon}
\def\ra{\rightarrow}

\def\ovl{\overline}

\def\sbs{\subset}

\def\reali{{\Bbb R}}
\def\complessi{{\Bbb C}}
\def\zzz{{\Bbb Z}}
\def\naturali{{\Bbb N}}

\def\ident{{\mathchoice {\rm 1\mskip-4mu l} {\rm 1\mskip-4mu l}
{\rm 1\mskip-4.5mu l} {\rm 1\mskip-5mu l}}}
\def\complessi{{\mathchoice {\setbox0=\hbox{$\displaystyle\rm C$}\hbox{\hbox
to0pt{\kern0.4\wd0\vrule height0.9\ht0\hss}\box0}}
{\setbox0=\hbox{$\textstyle\rm C$}\hbox{\hbox
to0pt{\kern0.4\wd0\vrule height0.9\ht0\hss}\box0}}
{\setbox0=\hbox{$\scriptstyle\rm C$}\hbox{\hbox
to0pt{\kern0.4\wd0\vrule height0.9\ht0\hss}\box0}}
{\setbox0=\hbox{$\scriptscriptstyle\rm C$}\hbox{\hbox
to0pt{\kern0.4\wd0\vrule height0.9\ht0\hss}\box0}}}}
\def\bbbe{{\mathchoice {\setbox0=\hbox{\smalletextfont e}\hbox{\raise
0.1\ht0\hbox to0pt{\kern0.4\wd0\vrule width0.3pt height0.7\ht0\hss}\box0}}
{\setbox0=\hbox{\smalletextfont e}\hbox{\raise
0.1\ht0\hbox to0pt{\kern0.4\wd0\vrule width0.3pt height0.7\ht0\hss}\box0}}
{\setbox0=\hbox{\smallescriptfont e}\hbox{\raise
0.1\ht0\hbox to0pt{\kern0.5\wd0\vrule width0.2pt height0.7\ht0\hss}\box0}}
{\setbox0=\hbox{\smallescriptscriptfont e}\hbox{\raise
0.1\ht0\hbox to0pt{\kern0.4\wd0\vrule width0.2pt height0.7\ht0\hss}\box0}}}}
\def\razionali{{\mathchoice {\setbox0=\hbox{$\displaystyle\rm Q$}\hbox{\raise
0.15\ht0\hbox to0pt{\kern0.4\wd0\vrule height0.8\ht0\hss}\box0}}
{\setbox0=\hbox{$\textstyle\rm Q$}\hbox{\raise
0.15\ht0\hbox to0pt{\kern0.4\wd0\vrule height0.8\ht0\hss}\box0}}
{\setbox0=\hbox{$\scriptstyle\rm Q$}\hbox{\raise
0.15\ht0\hbox to0pt{\kern0.4\wd0\vrule height0.7\ht0\hss}\box0}}
{\setbox0=\hbox{$\scriptscriptstyle\rm Q$}\hbox{\raise
0.15\ht0\hbox to0pt{\kern0.4\wd0\vrule height0.7\ht0\hss}\box0}}}}
\def\bbbt{{\mathchoice {\setbox0=\hbox{$\displaystyle\rm
T$}\hbox{\hbox to0pt{\kern0.3\wd0\vrule height0.9\ht0\hss}\box0}}
{\setbox0=\hbox{$\textstyle\rm T$}\hbox{\hbox
to0pt{\kern0.3\wd0\vrule height0.9\ht0\hss}\box0}}
{\setbox0=\hbox{$\scriptstyle\rm T$}\hbox{\hbox
to0pt{\kern0.3\wd0\vrule height0.9\ht0\hss}\box0}}
{\setbox0=\hbox{$\scriptscriptstyle\rm T$}\hbox{\hbox
to0pt{\kern0.3\wd0\vrule height0.9\ht0\hss}\box0}}}}
\def\bbbs{{\mathchoice
{\setbox0=\hbox{$\displaystyle     \rm S$}\hbox{\raise0.5\ht0\hbox
to0pt{\kern0.35\wd0\vrule height0.45\ht0\hss}\hbox
to0pt{\kern0.55\wd0\vrule height0.5\ht0\hss}\box0}}
{\setbox0=\hbox{$\textstyle        \rm S$}\hbox{\raise0.5\ht0\hbox
to0pt{\kern0.35\wd0\vrule height0.45\ht0\hss}\hbox
to0pt{\kern0.55\wd0\vrule height0.5\ht0\hss}\box0}}
{\setbox0=\hbox{$\scriptstyle      \rm S$}\hbox{\raise0.5\ht0\hbox
to0pt{\kern0.35\wd0\vrule height0.45\ht0\hss}\raise0.05\ht0\hbox
to0pt{\kern0.5\wd0\vrule height0.45\ht0\hss}\box0}}
{\setbox0=\hbox{$\scriptscriptstyle\rm S$}\hbox{\raise0.5\ht0\hbox
to0pt{\kern0.4\wd0\vrule height0.45\ht0\hss}\raise0.05\ht0\hbox
to0pt{\kern0.55\wd0\vrule height0.45\ht0\hss}\box0}}}}
\def\boxf{{$\ulcorner \! \! \lrcorner \! \! \! \! \llcorner \! \! \urcorner$}}

\def\QED{{\ifmmode\sq\else{\unskip\nobreak\hfil
\penalty50\hskip1em\null\nobreak\hfil{\boxf}
\parfillskip=0pt\finalhyphendemerits=0\endgraf}\fi} \vskip 12 pt} 
%%%%%%%%%%%%%%%%%%%%%%%%%%%%%%%%%%%%%%%%%%%%%%%%%%%%%%%%%%%%%%%%%%%%%%%%%%%%%%%%%%%%%%%%%%%%%%
%%%%%%%%%%%%%%%%%%%%%%%%%%%%%%%%%%%%%%%%%%%%%%%%%%%%%%%%%%%%%%%%%%%%%%%%%%%%%%%%%%%%%%%%
\input biblio.tex

\font\fontgrande=cmss17
\centerline{\fontgrande Uniform asymptotic estimates of transition
	probabilities on combs}
\vskip 10 pt
{\noindent
\centerline{D.~Bertacchi - F.~Zucca}\par\noindent
\centerline{Universit\`a degli Studi di Milano}\par\noindent
\centerline{Dipartimento di Matematica F.~Enriques}\par\noindent
\centerline{Via Saldini 50, 20133 Milano, Italy.}\par}
\vskip 12 pt
{\leftskip=80pt \rightskip=80pt \noindent {\bf Abstract.}\
We investigate the asymptotical behaviour of the transition probabilities
of the simple random walk on the $2$-comb. In particular we obtain space-time
uniform asymptotical estimates which show the lack of symmetry of this
walk better than local limit estimates. Our results also point out the
impossibility of getting Jones-type non-Gaussian estimates.
\par}



\vskip 18pt

\input intro.tex

\input tech.tex
\input fungen.tex

\input y.tex

\input x.tex

\input fr.tex

\vskip18pt
{\bf Acknowledgement}

\noindent We are deeply grateful to Wolfgang Woess for his unceasing
encouragement which helped us to reach our final goal.

\vskip36pt
{\bf Bibliography}
\bigskip
\insertbibliografia

\end

%% file: biblio.tex
\autobibliografia

\biblitem{Cassi-Regina}     
        D.~Cassi, S.~Regina, {\it Random walks on $d$-dimensional
        comb lattices},  \hfill\break 
        Modern Phys.~Lett.~B  {\bf 6} (1992), 1397-1403.

\biblitem{Avez}     
        A.~Avez, {\it Limite des quotients pour des marches al\'eatoires sur
                des groupes}, \hfill\break 
                C.~R.~Acad.~Sci.~Paris S\'er.~A {\bf 276}
                (1973), 317-320.

%%%%%%%%% qui iniziano quelli che ho messo io x la tesi
\biblitem{Bender}
        E.~A.~Bender, {\it Asymptotic method in enumeration}, 
        SIAM review, {\bf 4} (1974), 485-515.

\biblitem{Lal1}
        S.~P.~Lalley, {\it Saddlepoint approximations and space-time
		Martin boundary for nearest neighbour random walk on a
		homogeneous tree}, 
        J.~Theoret.~Probab. {\bf 4} (1991), 701-723.

\biblitem{Lal2}
        S.~P.~Lalley, {\it Finite range random walks on free groups and
		homogeneous trees}, 
        Ann.~Probab. {\bf 21} (1993), 2087-2130.

\biblitem{Survey}     
        W.~Woess, {\it Random walks on infinite graphs and groups - A 
        survey on selected topics}, Bull.~London Math.~Soc. {\bf 26} (1994),
        1-60.

\biblitem{Woess-rosa}
        W.~Woess, {\it Catene di Markov e Teoria del Potenziale nel
        Discreto}, Quaderno U.M.I {\bf 41}, Ed.~Pitagora, Bologna (1996).

\biblitem{Woess-giallo}
	W.~Woess, {\it Random walks on infinite graphs and groups}, 
        Cambridge Tracts in Mathematics {\bf 138}, Cambridge
		Univ. Press, 2000.

\biblitem{Ber-Zuc}
        D.~Bertacchi, F.~Zucca, {\it Equidistribution of random walks
        on spheres}, J.~Stat.~Phys. {\bf 94} (1999), 91-111.

\biblitem{Weiss-Havlin}
        G.~H.~Weiss, S.~Havlin, 
        %{?}, 
        Physica {\bf 134A} (1986), 474.
        %-?.

\biblitem{Alexopoulos}
        G.~Alexopoulos, {\it A lower estimate for central probabilities
        on polycyclic groups},  \hfill\break 
        Can.~J.~Math. {\bf 44} (5) (1992),
        897-910.

\biblitem{CKW}
        D.~I.~Cartwright, V.~A.~Kaimanovich, W.~Woess, {\it Random walks on the
        affine group of local fields and of homogeneous trees}, 
        Ann.~Inst.~Fourier (Grenoble) {\bf 44} (1994), 1243-1288.

\biblitem{Bi} 
        P.~Billingsley, {\it Convergence of probability measures},
        Wiley, New York, 1968.

\biblitem{F}
        W.~Feller, {\it An introduction to probability theorey and its
        applications}, Vol II, Wiley, New York, 1966.

\biblitem{Cartier}
        P.~Cartier, {\it Fonctions harmoniques sur un arbre},
        Symposia Math. {\bf 9} (1972), 203-270.

\biblitem{Ca-Ku-St}
        E.~Carlen, S.~Kusuoka, D.~Stroock, {\it Upper bounds for symmetric
        transition functions}, Ann.~Inst. H.~Poincar\'e Probab.~Statist.
        {\bf 23} (1987), 245-287.

\biblitem{Co1}
        T.~Coulhon, {\it Ultracontractivity and Nash type inequalities},
        J.~Funct.~Analysis {\bf 141} (1996), 510-539.

\biblitem{Dod}
        J.~Dodziuk, {\it Difference equations, isoperimetric inequality,
        and transience of certain random walks}, Trans.~Amer.~Math.~Soc.
        {\bf 284} (1984),787-794.

\biblitem{Ge10}
        P.~Gerl, {\it Random walks on graphs with a strong isoperimetric
        inequality}, J.~Theoret.~Probab. {\bf 1} (1988), 171-187.

\biblitem{Ki}
        J.~F.~C.~Kingman, {\it The ergodic theory of subadditive
        processes}, J.~Royal Stat.~Soc., Ser.~B {\bf 30} (1968),
        499-510.

\biblitem{Ba-Ba2}
        M.~T.~Barlow, R.~F.~Bass, {\it Random walks on graphical
        Sierpi\'nski carpets}, Symposia Math., in press.

\biblitem{Jo}
        O.~D.~Jones, {\it Transition probabilities for the simple
        random walk on the Sierpi\'nski graph}, Stochastic Proc.~Appl.
        {\bf 61} (1996), 4-69.

\biblitem{Te1}
        A.~Telcs, {\it Spectra of graphs and fractal dimensions. $I$},
        Probab.~Th.~Rel.~Fields {\bf 85} (1990), 489-497.

\biblitem{Te2}
        A.~Telcs, {\it Spectra of graphs and fractal dimensions. $I\!I$},
        Probab.~Th.~Rel.~Fields {\bf 8} (1995), 77-96.

\biblitem{CSV}
        T.~Coulhon, L.~Saloff-Coste, N.~Th.~Varopoulos, {\it Analysis
        and geometry on groups}, Cambridge Tracts in Mathematics,
        {\bf 100}, 1992.

\biblitem{Gerl2}
        P.~Gerl, {\it Natural spanning trees of $\zzz^d$ are recurrent},
        Discrete Math. {\bf 61} (1986), 2-3, 333-336.

\biblitem{Spitzer}
        F.~Spitzer, {\it Principles of
        Random Walks}, Springer, New York, 1976.

%%%%%%%%%%%%%%%%%%%%%%%%%%%%%%%%%%%%%%%%%%%%%%%%%%%%%%%%%%%%%%%%%%%%%%%%%%%%%%%%%%%%%%%%

%% file: intro.tex
%%%%%%%%%%%%%%%%%%%%%%%%%%%%%%%%%%%%%%%%%%%%%%%%%%%%%%%%%%%%%%%%%%%%%%%%%%%%
%%%    INTRODUZIONE
%%%%%%%%%%%%%%%%%%%%%%%%%%%%%%%%%%%%%%%%%%%%%%%%%%%%%%%%%%%%%%%%%%%%%%%%%%%%
\autosez{intro} Introduction

Given a random walk $(Z_n)_{n\ge0}$ on a graph $X$, there are many related
questions regarding its behaviour when the discrete time parameter $n$
goes to infinity.
Classical questions of this kind are, for instance: will the random walk
visit a given vertex of the graph only a finite number of times (with
probability one)? will it leave any bounded set after a finite time
(with probability one)?
Moreover, if we denote by $p^{(n)}(x,y)$ the $n$-step probabilities of the
random walk from the vertex $x$ to the vertex $y$ of $X$, we can study some
features of the sequence $(p^{(n)}(x,y))_n$, for instance answer to the 
question: is it asymptotic to some numerical sequence?

Answers to the first two questions are theorems and criteria for recurrence
and transience; answers to the latter question are provided by {\it local
limit theorems}, that is theorems which give a numerical estimate of
$p^{(n)}(x,y)$ for fixed $x,y$ as $n$ tends to infinity.

The present work studies the asymptotic behaviour of the
transition probabilities of the simple random walk on the 2-comb, which is
a graph obtained attaching at each point of $\zzz$ another copy of $\zzz$
by its origin. 
In other words, this graph is obtained from the usual square grid $\zzz^2$
by deleting all edges that are parallel to the $x$-axis except those on
the $x$-axis itself.

A local limit theorem for the $2$-comb (and in general for $d$-dimensional 
combs) and $x=y$ is well known: we observe how to extend it to the asymptotic
estimate of $(p^{(n)}(x,y))_n$ for any $x$ and $y$ (equation~\eqref{loc-lim2}).
The simple random walk on the $2$-comb lacks of symmetry (more precisely it
is not isotropic -- see Bertacchi and Zucca~\cite{Ber-Zuc}) and this
feature is not shown by local limit estimates.
To stress the different behaviour the random walk has in the two principal directions
(vertical and horizontal: see Figure~1) one needs space-time estimates.

A {\it space-time asymptotic estimate} is
a result which provides an estimate of $p^{(n)}(x,y)$ as $n$ tends to
infinity, uniform with respect to the quotient $d(x,y)/n$ lying in a suitable
range. Of course local limit theorems can be derived from space-time
estimates.
We provide space-time asymptotic estimates for the $p^{(n)}(x,y)$ when
$y=o:=(0,0)$ and $x=(k,0)$ or $x=(0,k)$, that is results of the form
$$
  p^{(n)}(x,o)\asym{n} Cf\left({d(x,o)\over n},n\right)n^{-d}
$$
where $C$, $d$ are constants and $f$ is a real valued function which all
depend on the range of $d(x,o)/n$. From these results all known limit
theorems for the 2-comb can be derived as a particular case.

The technique we exploit is essentially a Laplace-type estimate of integrals,
since the transition probabilities can be written as integrals thanks to the
Cauchy formula for the coefficients of the power series of an holomorphic
function:
$$
  p^{(n)}(x,o)={1\over 2\pi i}\int_\gamma{G(x,o|z)\over z^{2n+1}}\d z
$$
where $G(x,o|z):=\sum_{n=0}^\infty p^{(n)}(x,o)z^n$ is the Green
function associated to the walk and has an explicit expression, and $\gamma$
is a positively oriented, simple closed curve in $\complessi$ surrounding
$0$.
The basic idea is that the main part of the integral is given by integration
on the part of the curve which is closer to the singularity $z=1$ of 
$G(x,0|z)$. To develop this idea into mathematical terms we first separate in the
integrand the part with algebraic behaviour and the one with exponential
behaviour, then we choose a suitable curve of integration and its
parametrization. Afterwards, we write the Taylor expansion of the argument of
the exponential part of the integrand as a function of the parameter of the
curve and we finally show that it is possible to choose a piece of the curve
on which integration gives the asymptotic behaviour of the transition
probabilities.
%%%%%%%%%%%%%%%%%%%%%%%%%%%%%%%%%%%%%%%%%%%%%%%%%%%%%%%%%%%%%%
It is remarkable that the above mentioned Taylor expansions are very
different in the two cases $x=(k,0)$ and $x=(0,k)$. This reflects in two 
different asymptotic behaviours of the transition probabilities when
$d(x,o)/n$ tends to $0$: in the first case 
$$
p^{(n)}(x,o)\asym{n}c_1\exp(c_2\,n(d(x,o)/n)^{4/3})\,n^{-3/4},
$$
while in the second case
$$
p^{(n)}(x,o)\asym{n}c_3\exp(c_4\,n(d(x,o)/n)^{2})\,n^{-3/4}.
$$
In particular this shows that for the 2-comb it is impossible to find a
uniform estimate involving the spectral dimension and the walk dimension
such as the one found by Jones for the $2$-dimensional Sierpinski graph 
(see Paragraph~\sref{fr}).

\medskip
%%%%%%%%%%%%%%%%%%%%%%%%%%%%%%%%%%%%%%%%%%%%%%%%%%%%%%%%%%%%%%%%%%%%%%%%%%%%

We give a brief outline of the paper: in Paragraph~\sref{tech} we listed
the definition of uniform estimate with respect to a parameter and two
Lebesgue-type theorems which are needed to obtain such estimates.
Paragraph~\sref{fungen} recalls local limit theorems and generating functions for
combs.
Then in Paragraphs~\sref{y} and \sref{xi<a} uniform estimates are proved for
the vertical direction, while in Paragraphs~\sref{x}, \sref{x-xi<a},
\sref{xi<n3/4} and \sref{xi>n3/4} we prove uniform estimates for the
horizontal direction.
In the last paragraph we discuss the obtained results and possible extensions.

%% file: tech.tex
%%%%%%%%%%%%%%%%%%%%%%%%%%%%%%%%%%%%%%%%%%%%%%%%%%%%%%%%%%%%%%%%%%%%%%%%%%%%%%%%%%%%%%%%%%%%%%
%% tech
%%%%%%%%%%%%%%%%%%%%%%%%%%%%%%%%%%%%%%%%%%%%%%%%%%%%%%%%%%%%%%%%%%%%%%%%%%%%%%%%%%%%%%%%%%%%%%
\autosez{tech}{Asymptotic estimates: definitions and technical results}

In this paper we are concerned with asymptotic estimates of transition
probabilities. We then give some useful definitions and theorems.

\definition{unif1}{
Given two sequences $\big(a_n(x)\big)_{n}$ and $\big(b_n(x)\big)_{n}$
of complex 
functions defined on a space $X$ and a sequence $(A_n)_{n}$ of 
subsets of $X$, we say that $a_n$ is asymptotic to $b_n$ (and we write 
$a_n \asym{n} b_n$) as $n$ tends to infinity, uniformly 
with respect to $x \in A_n$ 
if and only if there exists a sequence of
complex functions
$\big(o_n(x)\big)_{n}$ and $n_0 \in \naturali$ such that:
\item{(i)}
$a_n(x)=b_n(x)(1+o_n(x))$
for every $n \ge n_0$, $x \in A_n$; 
\item{(ii)} for every $\eps >0$ there exists $n_\eps\in \naturali$,
$n_\eps \ge n_0$, such that for every $n>n_\eps$, $x \in A_n$ we have
$|o_n(x)| < \eps$.
}

This definition extends the usual one (let $X$ be a one point space
and $A_n\equiv X$ for all $n\in\naturali$).

\definition{unif2}{
Let us consider a sequence $\big(a_n(x)\big)_n$ of   
functions defined on a space $X$ with values on a metric space $(Y,d)$ and
let $(A_n)_n$ be a sequence of subsets of $X$. 
Let $b:X \rightarrow (Y,d)$, we say that 
$(a_n)_n$ converges to $b$ when $n$ tends to infinity, uniformly 
with respect to $x \in A_n$ if and only if 
for every $\eps >0$ there exists $n_\eps \in \naturali$
such that for every $n\ge n_\eps$, $x \in A_n$ we have
$d(a_n(x), b(x)) < \eps$.
}

This is an extension of the well-known definition of uniform convergence
(usually stated for $A_n=A\sbs X$).
The following result provides an extension of Lebesgue's bounded convergence
theorem to the case of uniform convergence (in the preceding sense).

\theorem{dominunif}{Let $(X, \Sigma, \mu)$ be a $\sigma$-finite measure
space and let $(f_n)_n: (X, \Sigma, \mu) \times Y \ra\complessi$ 
be a sequence of complex functions which are measurable with respect to 
$x \in X$ for every fixed $y \in Y$
and $(A_n)_n \in {\cal P}({Y})$ 
such that $f_n(x,y) \ra f(x,y)$ 
when $n$ tends to infinity uniformly 
with respect to $x \in A_n$. Then $f$ is measurable with respect to $x \in X$
for every fixed $y \in Y$ and if there exists $g \in L^1(X, \Sigma, \mu)$
such that $|f_n(x,y)| \le g(x)$, $\mu$ a.e., for every $n \in \naturali$
and for every $y \in A_n$, then
$f_n(\cdot, y) \rightarrow f(\cdot, y)$ in $L^1(X, \Sigma, \mu)$ uniformly
with respect to $y \in A_n$.
}
\proof
Let us fix $\epsilon >0$ and a non decreasing covering of $X$,
$(M_n)_{n}$ with finite measure measurable sets.
From Lebesgue's monotone convergence theorem, if 
$h \in L^1(X, \Sigma, \mu)$, then there exists $n_\epsilon \in \naturali$ 
such that for every $n \geq n_\epsilon$, we have $\int_X |h| \d \mu <
\epsilon$.
Moreover, by the absolute continuity of the Lebesgue integral, there exists
$\delta_\epsilon$ such that if $E \in \Sigma$, $\mu(E) < \delta_\epsilon$
then $\int_E |h| \d \mu < \epsilon$.

Now let 
$$
M_m(n):= \left\{ x \in M_m: \forall p \geq n, |f_p(x,y) - f(x,y)| < {\epsilon 
	\over\mu(M_m)}, \forall y \in A_p\right\},
$$
where $m,n \in \naturali$. It is easy to show that $M_m(n) \supseteq
M_m(n-1)$ and $\mu(M_m \setminus \cup_{n \in \naturali} M_m(n))=0$, then
if we take $m, n_0 \in \naturali$ 
such that $\int_{M_m^c} g \d \mu < \epsilon$
and $\mu(M_m \setminus M_m(n))< \delta_\epsilon$, $\forall n \geq n_0$ then
$$
\eqalign{
& \int_X |f_n(x,y)-f(x,y)| \d \mu =\cr
&= \int_{M_m^c} 
|f_n(x,y)-f(x,y)| \d \mu
+\int_{M_m} |f_n(x,y)-f(x,y)| \d \mu \leq \cr
& \leq 2 \int_{M_m^c} g(x) \d \mu + 2 \int_{M_m \setminus M_m(n)}
g(x) \d \mu + \int_{M_m(n)} |f_n(x,y)-f(x,y)| \d \mu < 5 \epsilon \cr
}
$$
for all $n \geq n_0$, and for all $y \in A_n$.
\QED

Sometimes it will be impossible to exhibit a limit function for our
estimates, but we will be able to find a sequence of functions asymptotic to 
the given one and much simpler. In that direction works the following theorem.

\theorem{dominscon1}
{Let $(X, \Sigma, \mu)$ be a measure space and let 
$(f_n)_n, (h_n)_n, (o_n)_n : (X, \Sigma, \mu) \times Y \rightarrow
\complessi$ be three sequences of functions which are measurable with
respect to $x \in X$ for every fixed $y \in Y$. Let  
$(A_n)_n$ be a sequence of sets in ${\cal P} (Y)$ such that
$h_n(x,y) = f_n(x,y)(1+ o_n(x,y))$, $\mu$ a.e., for all $n \in \naturali$, 
for all $y \in A_n$, where
$o_n(x,y) \rightarrow 0$, when $n$ tends to infinity, $\mu$ a.e., 
uniformly with respect to $y \in A_n$.
If there exist $g,g_1 \in L^1(X, \Sigma, \mu)$ such that
$$
|h_n (x,y)| \leq g(x), \quad |f_n(x,y)| \leq g_1(x)
$$
$\mu$ a.e., for all $n \in \naturali$, for all $y \in A_n$ and we can
choose $c>0$ such that $\left | \int_X f_n(x,y) \d \mu \right |\ge c$,
for all $n \in \naturali$ and for all $y \in A_n$ then
$$
\int_X f_n(x,y) \d \mu \asym{n}\! \int_X h_n(x,y) \d\mu
$$
uniformly with respect to $y \in A_n$. 
}
\proof
Let us fix $\epsilon>0$ and define $X_n(\epsilon)=\{x \in X: 
|o_m(x,y)|< \epsilon, \forall m \geq n, \forall y \in A_m\}$ then
$X_n(\epsilon) \supseteq X_{n-1}(\epsilon)$ and $\mu(X \setminus 
\cup_{n \in \naturali} X_n(\epsilon))=0$. For all $n \in \naturali$ we
define 
$$\eqalign{
  & a_n(\epsilon):= \int_{X_n^c(\epsilon)} g \d \mu,\cr
  & b_n(\epsilon):= \int_{X_n^c(\epsilon)} g_1 \d \mu,
}$$
and it is easy to show
that $a_n(\epsilon) \rightarrow 0$, $b_n(\epsilon) \rightarrow 0$
if $n$ tends to infinity.
Moreover it is useful to fix 
$$\matrix{
I(y,n):= \int_X f_n(x,y) \d \mu,\hfill & I_1(y,n):= \int_X h_n(x,y) \d \mu,\hfill\cr
I_2(y,n):= \int_{X_n(\epsilon)} f_n(x,y) \d \mu,\hfill & \hfill\cr
}$$
 (the thesis is that
$I(y,n)\asym{n}I_1(y,n)$ uniformly with respect to $y\in A_n$);
then
$$ |I(y,n) - I_1(y,n)| \leq |I(y,n)-I_2(y,n)|+|I_1(y,n)-I_2(y,n)|.$$
Now $$
|I(y,n) - I_2(y,n)| \le \int_{X_n^c(\epsilon)} |f_n(x,y)| \d \mu
\leq b_n(\epsilon), \quad \forall n, \forall y \in A_n;$$
besides
$$\eqalign{
|I_1(y,n) - I_2(y,n)| &\le \int_{X_n(\epsilon)} |f_n(x,y)-h_n(x,y)| 
    \d \mu + \int_{X_n^c(\epsilon)} |h_n(x,y)| \d \mu \le\cr
	& \le\epsilon \int_{X_n(\epsilon)} |f_n(x,y)| \d \mu +
	a_n(\epsilon) \le \cr
	&\le \epsilon  
	\|g_1\|_1 + a_n(\epsilon), \quad \forall n, \forall y \in A_n.\cr
}$$
Keeping in mind that $|I(n,y)| \ge c$ we have
$$
\eqalign{
\left|1-{I_1(y,n)\over I(y,n)}\right| &
	\le {|I(y,n)-I_2(y,n)|\over |I(y,n)|}+
	{|I_1(y,n)-I_2(y,n)|\over |I(y,n)|} \le \cr
	&\le {1 \over c} (b_n(\epsilon)+
	\epsilon \|g_1\|_1)+{1 \over c} a_n(\epsilon), \quad
	\forall n, \forall y \in A_n.\cr
}
$$
Now we take $n_\epsilon$ such that for all $n \geq n_\epsilon$ we have
$a_n(\epsilon)+b_n(\epsilon) \leq 2 \epsilon \|g\|_1$
and we finally obtain
$$
\left|1-{I_1(y,n)\over I(y,n)}\right| 
\leq 3 \epsilon c^{-1} \|g\|_1, \quad \forall n \geq 
n_\epsilon, \forall y \in A_n.$$
\QED

The last (non-standard) proposition we need deals with
the triangular inequality for power series: we are interested
in the cases where a strict inequality holds.

\proposition{seriepot}{Let $f(z)=\sum_{n=0}^\infty a_nz^n$ be a power series
with non negative coefficients (and positive radius of convergence): then 
$|f(z)|\le f(|z|)$. Moreover, if there
exists $n\in\naturali$ such that $a_n$ and $a_{n+1}$ are both strictly 
positive, then $|f(z)|=f(|z|)$ if and only if $z=|z|$.}

\proof
It is well known that for every measurable space $(X,{\cal M},\mu)$ and for
every measurable function $g$, $|\int_Xg(x)\d\mu(x)|\le\int_X|g(x)|\d\mu(x)$
and it can be easily shown that equality holds if and only if there exists 
$\theta\in [0,2\pi)$ such that $|g(x)|=g(x)e^{i\theta}$ $\mu$-almost everywhere.
If we apply the first fact to our case, that is to $X=\naturali$, ${\cal M}=
{\cal P}(\naturali)$, $\mu(\{n\})=a_n$ and $g(n)=z^n$, then
$\int_Xg(x)\d\mu(x)=\sum_{n=1}^\infty a_nz^n=f(z)$  and the first part of the
claim is proved.

Equality obviously holds if $z=|z|$. Now suppose equality holds: this can be only 
if there exists $\theta\in [0,2\pi)$ such that $|z|^n=z^ne^{i\theta}$ for all
$n$ such that $a_n\neq 0$. By assumption, there exists $k\in\naturali$ such
that either $a_{2k}\neq 0$ and $a_{2k+1}\neq0$ or $a_{2k-1}\neq 0$ and 
$a_{2k}\neq0$. If we have the first case, then $|z|^{2k}=z^{2k}e^{i\theta}$
and $|z|^{2k+1}=z^{2k}e^{i\theta}$. Hence there exist $m,p\in\naturali$ such
that $2k\arg(z)+\theta=2m\pi$ and $(2k+1)\arg(z)+\theta=2p\pi$, that is
$\arg(z)=2(m-p)\pi$ and $z=|z|$. The proof in the case $a_{2k-1}\neq 0$,
$a_{2k}\neq0$ is completely analogous. \QED

%% file: fungen.tex
%%%%%%%%%%%%%%%%%%%%%%%%%%%%%%%%%%%%%%%%%%%%%%%%%%%%%%%%%%%%%%%%%%%%%%%%%
%%%%%%%%%%%%%%%%%%%%%%%%%%%%%%%FUNGEN.TEX%%%%%%%%%%%%%%%%%%%%%%%%%%%%%%%%%
%%%%%%%%%%%%%%%%%%%%%%%%%%%%%%%%%%%%%%%%%%%%%%%%%%%%%%%%%%%%%%%%%%%%%%%%%%%%%%%%%%%%%%%%%%%%%%
\autosez{fungen}{Local limit theorems and generating functions}

The {\it $2$-comb lattice} $C_2$ is a {\it spanning tree} of $\zzz^2$,
that is a subgraph of $\zzz^2$ which is a tree and contains all vertices.
More precisely, we obtain the $2$-comb 
by deleting in $\zzz^2$ all edges that are parallel to the $x$-axis besides 
those on the $x$-axis itself.
This construction provides a natural choice of coordinates on the comb
(that is $(x,y)\in C_2$ indicates the same point of $\zzz^2$, now thought
as belonging to $C_2$).
Alternatively, we construct it by attaching at each point of $\zzz$ a
two way-infinite path (that is a copy of $\zzz$).

\beginpicture
%%%%%%%%IL 2-COMB
\setcoordinatesystem units <.4cm,.4cm> point at  0 0
\put{$\,$} at 0 7
\setlinear
\plot 10  0  26  0   /
\plot 12 -6  12 6 /
\plot 15 -6  15  6 /
\plot 18 -6  18  6 /
\plot 21 -6  21  6  /
\plot 24  -6  24  6  /

\setdots <2.5pt>
\put{$\,$} at -2 0
\put{\it Figure~1: the $2$-comb} [l] at 0 0
\put{$\scriptstyle(0,0)$} [rt] at 17.5 -.5 
\plot 17.5 -.5  18 0 /
\put{$\scriptstyle(-1,0)$} [rt] at 14.5 -.5 
\plot 14.5 -.5  15 0 /
\put{$\scriptstyle(-2,0)$} [rt] at 11.5 -.5 
\plot 11.5 -.5  12 0 /
\put{$\scriptstyle(1,0)$} [rt] at 20.5 -.5 
\plot 20.5 -.5  21 0 /
\put{$\scriptstyle(2,0)$} [rt] at 23.5 -.5 
\plot 23.5 -.5  24 0 /
\put{$\scriptstyle(0,1)$} [rt] at 17.5 2.5 
\plot 17.5 2.5  18 3 /

\put{$\scriptscriptstyle\bullet$} at 12 0
\put{$\scriptscriptstyle\bullet$} at 15 0
\put{$\scriptscriptstyle\bullet$} at 18 0
\put{$\scriptscriptstyle\bullet$} at 21 0
\put{$\scriptscriptstyle\bullet$} at 24 0
\put{$\scriptscriptstyle\bullet$} at 12 3
\put{$\scriptscriptstyle\bullet$} at 12 -3
\put{$\scriptscriptstyle\bullet$} at 15 3
\put{$\scriptscriptstyle\bullet$} at 15 -3
\put{$\scriptscriptstyle\bullet$} at 18 3
\put{$\scriptscriptstyle\bullet$} at 18 -3
\put{$\scriptscriptstyle\bullet$} at 21 3
\put{$\scriptscriptstyle\bullet$} at 21 -3
\put{$\scriptscriptstyle\bullet$} at 24 3
\put{$\scriptscriptstyle\bullet$} at 24 -3

\put{$\,$} at 0 -7
\endpicture

More generally, $d$-dimensional comb lattices $C_d$ are the spanning trees
of $\zzz^d$ obtained inductively by attaching at each point of $C_{d-1}$
a copy of $\zzz$.

The estimate of the asymptotic behaviour of the transition probabilities
of the simple random walk on comb lattices passes through the
knowledge of the generating functions, which we now recall.

Given two vertices of the comb, $x$ and $y$, and the Markov chain 
$(Z_n)_n$ representing the simple random walk
on the comb starting at $x$, the $n$-step transition 
probabilities are defined by
$$
p^{(n)}(x,y)=\prob [Z_n=y|Z_0=x].
$$
The Green function is then the power series
$$
	G(x,y|z)=\sum_{n=0}^\infty p^{(n)}(x,y) z^n, \quad z\in\complessi,
$$
while 
$$
	F(x,y|z)=\sum_{n=0}^\infty \prob [{\bf s}^y=n|Z_0=x] z^n,
 \quad z\in\complessi,
$$
where ${\bf s}^y:=\inf [n\ge 0|Z_n=y]$.
Then it is well known (see for instance \cite{Survey} and 
\cite{Woess-rosa})
that in the common domain of convergence of these power series
$G(x,y|z)=F(x,y|z)G(y,y|z)$.

The Green function of the $d$-dimensional comb, $G_d$, can be obtained recursively by the following 
formula (see \cite{Cassi-Regina})
$$
G_d(z) = {d 
\over \sqrt{ \left ( 1+{d-1 \over G_{d-1}(z)} \right )^2 -z^2}},
\autoeqno{greencomb1}
$$
recalling that $G_1(z)=1/\sqrt{1-z^2}$.

From equation~\eqref{greencomb1} it is easy to see that $G_d(z) = G_d(-z)$,
and that $\{+1, -1\}$ are algebraic singularities (see
\cite{Bender} pag.~498) since these properties hold for $G_1$.

By using mathematical induction on eq.~\eqref{greencomb1}, 
it is not difficult to see that 
$$
G_d(z) \,\asym{z \rightarrow 1^\pm}\, d2^{1/2^d-1}(1\pm z)^{-1/2^d}.
\autoeqno{greencomb3}
$$
Local limit theorems, that is numerical estimates of $p^{(n)}(o,o)$ where
by $o$ we denote the origin of the grid $\zzz^d$ where $C_d$ is embedded,
were studied by Weiss and Havlin (see \cite{Weiss-Havlin}) in the case
$d=2$ or $d=3$, and extended to the general case by Cassi and Regina
(see \cite{Cassi-Regina}).
In fact applying Theorem 4 of \cite{Bender} to equation~\eqref{greencomb3}
we have that
$$
p^{(n)}(o,o) \asym{n}
\cases{\displaystyle 0  \hfill & if $n$ is odd \cr
\cr
\displaystyle { 2^{1 / 2^d+1}  \over 
\Gamma \left ( 1 \over 2^d \right )} n^{1/2^d -1} &
if $n$ is even.
}
$$
From this particular estimate one can derive an
estimate of the  general transition probabilities (see Bertacchi and
Zucca~\cite{Ber-Zuc} Paragraph 6):
$$
p^{(n)}(x,y) \asym{n}
\cases{\displaystyle 0  \hfill & if $n+d(x,y)$ is odd \cr
\cr
\displaystyle { 2^{1 / 2^d-1} \deg(y) \over 
\Gamma \left ( 1 \over 2^d \right )} n^{1/2^d -1} &
if $n+d(x,y)$ is even.
}\autoeqno{loc-lim2}
$$

Our goal is to obtain asymptotic estimates of
$p^{(n)}((0,k),(0,0))$ and of $p^{(n)}((k,0),(0,0))$, $k\ge0$, {\it 
uniform} with 
respect to the parameter $\xi={k\over n}$ (in order to avoid discussions 
about the parity of $n$ and $k$ we will only deal with the case where they 
both are even, but this is no severe restriction ---
see Paragraphs~\sref{y} and \sref{fr}).

We will make use of Cauchy's integral formula
$$
  p^{(n)}(x,y)={1\over 2\pi i}\int_{\gamma}{G(x,y|z)
	\over z^{n+1}}\d z
$$
($\gamma$ is a positively oriented, simple closed curve in $\complessi$
with $0$ in its interior), and then estimate this integral.
Hence we need
the explicit expressions of $G((k,0),(0,0)|z)$ and
$G((0,k),(0,0)|z)$, which are the following:
$$
 \eqalign{
	G\big((k,0),(0,0)|z\big) & = {\sqrt 2\over\sqrt{1-z^2+\sqrt{1-z^2}}}
	\cdot   \left({1+\sqrt{1-z^2}-\sqrt 2
		\sqrt{1-z^2+\sqrt{1-z^2}}\over z}\right)^{|k|}\cr
	G\big((0,k),(0,0)|z\big) & = {\sqrt 2\over\sqrt{1-z^2+\sqrt{1-z^2}}}
	\cdot   \left({1-\sqrt{1-z^2}\over z}   \right)^{|k|}.
}$$
Computation uses well known techniques for the generating functions on 
graphs involving explicit expressions of $G(o,o|z)$, $F((k,0),(0,0)|z)$ and
$F((0,k),(0,0)|z)$ (see for instance Woess~\cite{Woess-rosa} Theorem 1.4).

Now for simplicity we rename
$$\matrix{
	\widetilde G(z):={\sqrt 2\over\sqrt{1-z^2+\sqrt{1-z^2}}}; \hfill
	 & G(z):=\widetilde G(\sqrt z)={\sqrt 2\over\sqrt{1-z+\sqrt{1-z}}};\hfill\cr
	\widetilde F_1(z):={1-\sqrt{1-z^2}\over z};\hfill
	  & F_1(z):=\widetilde F_1(\sqrt z)={1-\sqrt{1-z}\over\sqrt z};\hfill\cr
	\widetilde F_2(z):={1+\sqrt{1-z^2}-\sqrt 2\sqrt{1-z^2+
	  \sqrt{1-z^2}}\over z};\hfill
	  & F_2(z):=\widetilde F_2(\sqrt z)={1+\sqrt{1-z}-\sqrt 2\sqrt{1-z+
	     \sqrt{1-z}}\over \sqrt z}.\hfill\cr
}$$
We note that these generating functions all contain radicals, hence we must 
pay attention to their polidromy. 
{\bf Remark}
The functions
$G(z)$, $F_1^2(z)$ and $F_2^2(z)$ are
holomorphic in the open ball with radius $1$ centered in $0$, that is, $0$
is not a singularity for any of these functions.
Moreover, $z=1$ is a branch point for all of them and their only singularity
in the complex plane.

{\bf Choice of the determination of the square root:} we choose, once for all, 
the determination of the square root with
argument between $-\pi/2$ and $\pi/2$, that is the function
$h(w):=\sqrt{|w|}\exp(i\arg(w)/2)$ where $\arg(w)$ is chosen in the interval
$[-\pi,\pi)$.
That means that $\sqrt{1-z}$ is an holomorphic function
defined in the open set $A:=\complessi\setminus\{z\in\reali\,:\,z>1\}$, and
we extend it to $\{z\in\reali\,:\,z>1\}$ by continuity from the upper
half plane (then we will not have continuity from the lower half plane).
Note that this choice allows us to define $\sqrt{1-z+\sqrt{1-z}}$ as
an holomorphic function in $A$ as well.

%% file: y.tex
%%%%%%%%%%%%%%%%%%%%%%%%%%%%%%%%%%%%%%%%%%%%%%%%%%%%%%%%%%%%%%%%%%%%%%%%%%%%%
%% INIZIANO I CONTI LUNGO y
%%%%%%%%%%%%%%%%%%%%%%%%%%%%%%%%%%%%%%%%%%%%%%%%%%%%%%%%%%%%%%%%%%%%%%%%%%%%%
\autosez{y}{Estimates along the $y$-axis: the case $\xi\in [a,1-c]$}

We first deal with the case of the asymptotic estimate 
for the transition probabilities
$p^{(n)}((0,k),(0,0))$, which we will rename for sake of simplicity                                      
$p^{(n)}(k,0)$.
We first note that if $n$ and $k$ do not have the same parity, then
$p^{(n)}(k,0)=0$, while
$$
  p^{(2n+1)}(2k+1,0) ={1\over 2}\{ p^{(2n)}(2k+2,0)+p^{(2n)}(2k,0)\}.
$$
Hence it will be enough to estimate $p^{(2n)}(2k,0)$.

By Cauchy's integral formula,
$$
  p^{(2n)}(2k,0)={1\over 2\pi i}\int_{\gamma^\prime}{\widetilde G(z)\cdot
	\widetilde F_1(z)^{2k}\over z^{2n+1}}\d z
$$
where $\gamma^\prime$ is a positively
oriented, simple closed curve in $\complessi$ which has $0$ in its interior
and $1$ in its exterior.

The substitution $u=z^2$ appears to be useful and when $\gamma^\prime\,:\,z=z(t)$
describes one circuit about $0$ then $\gamma\,:\,u=u(t)$ describes 
two times the corresponding circuit, hence
$$
%\eqalign{
  p^{(2n)}(2k,0) 
%& ={2\over 2\pi i}\int_{\gamma}{G(u)F_1(u)^{2k}\over u^n\sqrt u}\cdot
%	{\d u\over 2\sqrt u}\cr &
	 ={1\over 2\pi i}\int_\gamma{G(z)F(z)^{2k}\over z^{n+1}}\d z,
%}
$$
where $F(z)=F_1(z)$. 

We first stress the exponential part of the integrand, rewriting
$$
\eqalign{
 {F(z)^{2k}\over z^n}& =\exp\{2k\log(F(z))-n\log(z)\}=\cr
%	& =\exp\{2k\log(1-\sqrt{1-z})-(k+n)\log(z)\}=\cr
	& =\exp\{n[2\xi\log(1-\sqrt{1-z})-(\xi+1)\log(z)]\},
}$$
where $\xi:={k\over n}$.
We now introduce the auxiliary function 
$$
\Psi_\xi(z):=2\xi\log(1-\sqrt{1-z})-(\xi+1)\log(z)
\autoeqno{psiy}$$
 and rewrite the
expression for the transition probabilities:
$$
  p^{(2n)}(2k,0) ={1\over 2\pi i}
  \int_\gamma{G(z)\over z}\exp\{n\Psi_\xi(z)\}\d z.
\autoeqno{y1}
$$
\lemma{phixi}{The function $\Psi_\xi(z)$ has a unique minimum in $[0,1]$, 
namely $z_o(\xi)=1-\xi^2$.
Let $\phi(\xi)$ be this minimum: then 
$\phi(\xi)=\log\big((1-\xi)^{\xi-1}(1+\xi)^{-\xi-1}\big)$.
The second and the third order derivatives of $\Psi_\xi(z)$
take the following values in the generic point $z$ and in $z_o(\xi)$:
$$
\matrix{
  \Psi^{\prime\prime}_\xi(z) ={2(1-z)^{3/2}+\xi(3z-2)\over 2z^2(1-z)^{3/2}};\hfill
   & \Psi^{\prime\prime}_\xi(z_o(\xi)) ={1\over 2\xi^2(1-\xi^2)};\hfill\cr
  \Psi^{\prime\prime\prime}_\xi(z) =
	{\xi(15z^2-20z+8)-8(1-z)^{5/2}\over 4z^3(1-z)^{5/2}};\hfill
   & \Psi^{\prime\prime\prime}_\xi(z_o(\xi))={3-7\xi^2\over 4\xi^4(1-\xi^2)^2}.
	\hfill\cr
}$$}

%\proof
%We note that $\lim_{z\ra 0^+}\Psi_\xi(z)=+\infty$, $\Psi_\xi(1)=0$ and that
%for $z\in (0,1)$ the first derivative of $\Psi_\xi(z)$ is
%$$
%  \Psi^\prime_\xi(z) 
%	 ={1\over z}\left({\xi\over\sqrt{1-z}}-1\right).
%$$
%Then the derivative is negative if $\xi<\sqrt{1-z}$ that is if $0<z<1-\xi^2$, 
%positive if $1-\xi^2<z<1$ and $0$ if $z=1-\xi^2$, which is then the value for
%$z_o(\xi)$. The value of $\phi(\xi)=\Psi_\xi(z_o(\xi))$ follows from 
%elementary computations. \QED

%In the following, we will also need some other derivatives of $\Psi_\xi$,
%in particular evaluated in $z_o(\xi)$.
%\lemma{derPsi1}{The second and the third order derivatives of $\Psi_\xi(z)$
%take the following values in the generic point $z$ and in $z_o(\xi)$:
%$$
%\matrix{
%  \Psi^{\prime\prime}_\xi(z) ={2(1-z)^{3/2}+\xi(3z-2)\over 2z^2(1-z)^{3/2}};\hfill
%   & \Psi^{\prime\prime}_\xi(z_o(\xi)) ={1\over 2\xi^2(1-\xi^2)};\hfill\cr
%  \Psi^{\prime\prime\prime}_\xi(z) =
%	{\xi(15z^2-20z+8)-8(1-z)^{5/2}\over 4z^3(1-z)^{5/2}};\hfill
%   & \Psi^{\prime\prime\prime}_\xi(z_o(\xi))={3-7\xi^2\over 4\xi^4(1-\xi^2)^2}.
%	\hfill\cr
%}$$}
Here is the first estimate of our transition probabilities.

\theorem{y-xinonzero}{Let $a,c$ be positive numbers such that $a<1-c$ and let
$\xi\in[a,1-c]$. Then
$$
  p^{(2n)}(2k,0)\asym{n}
	{\xi\over\sqrt\pi\sqrt{1-\xi^2}}\,G(z_o(\xi))e^{n\phi(\xi)}\,n^{-1/2}=
  \sqrt{2\xi\over (1-\xi^2)(1+\xi)}            
  {e^{n\phi(\xi)}\over\sqrt\pi}   \,n^{-1/2}
$$
uniformly with respect to $\xi\in [a,1-c]$.}

\proof
 We first choose the curve of integration  for equation~\eqref{y1}
and split the integral into two
parts (Part $I$ of the proof); then we evaluate the part which will prove to be
asymptotically negligible as compared to the other (Part $I\!I$ of the proof)
and finally estimate the main part (Part $I\!I\!I$ of the proof).

{\it Part $I$}

\noindent The curve of integration will be the circle with radius $z_o(\xi)$, 
centered in 
the origin: $\gamma\,:\,z(\xi,t)=z_o(\xi)e^{it},\ t\in [-\pi,\pi]$.
We note that since $\xi\in [a,1-c]$ then  $z_o(\xi)\in [1-(1-c)^2,1-a^2]=:[\ovl a,
1-\ovl c]$.

We introduce the function $\ovl\Psi_\xi(t):=\Psi_\xi(z_o(\xi)e^{it})=
\Psi_\xi(z(\xi,t))$ and perform the change of variable $z=z(\xi,t)$ in the
integral~\eqref{y1}:
$$
%\eqalign{
 p^{(2n)}(2k,0) 
%& ={1\over 2\pi i}\int_{-\pi}^\pi{G(z_o(\xi)e^{it})\over
%      z_o(\xi)e^{it}}\exp\{n\ovl\Psi_\xi(t)\}iz_o(\xi)e^{it}\d t\cr &
       ={1\over 2\pi}\int_{-\pi}^\pi G(z(\xi,t))\exp\{n\ovl\Psi_\xi(t)\}\d t.
%}
  \autoeqno{t-int1}
$$
Now we want to write the Taylor expansion of $\ovl\Psi_\xi(t)$ with Lagrange
remainder, centered in $t=0$. This is possible since the third order derivative
of $\ovl\Psi_\xi(t)$ exists and is continuous in $t$, for all $\xi\in [a,1-c]$.
%,in fact
%$$
%\matrix{
%  \ovl\Psi_\xi(0) =\Psi_\xi(z_o(\xi))=\phi(\xi);\hfill
%   & \ovl\Psi_\xi^\prime(0)=0;\hfill\cr
%  \ovl\Psi_\xi^{\prime\prime}(0) =-{1-\xi^2\over 2\xi^2};\hfill
%    & \ovl\Psi_\xi^{\prime\prime\prime}(0) =-i{(1-\xi^2)(3-\xi^2)\over 4\xi^4};
%	    \hfill\cr 
%  \ovl\Psi_\xi^{\prime\prime\prime}(t) =
%       -i\xi{(1-\xi^2)e^{it}((1-\xi^2)e^{it}+2)\over 4(1-(1-\xi^2)e^{it})^{5/2}}.
%		\hfill\cr 
%  }$$
Hence we can write $\ovl\Psi_\xi(t)=\phi(\xi)-{1-\xi^2\over 4\xi^2}t^2+R(\xi,t)$,
where the remainder is $R(\xi,t)={\ovl\Psi^{\prime\prime\prime}_\xi
(\ovl t)\over 3!}t^3$, and $\ovl t$ is a point lying in the segment between
$0$ and $t$.

In particular we note that $\ovl\Psi_\xi^{\prime\prime\prime}(0)$ is bounded
and far from $0$ for $\xi\in [a,1-c]$, that is there exists $\eps>0$ such that
$|\ovl\Psi_\xi^{\prime\prime\prime}(0)|>\eps$ for all $\xi\in[a,1-c]$.
Then $R(\xi,t)=O(t^3)$ for $t\ra 0$, 
uniformly with respect to $\xi\in [a,1-c]$, that is, there exists $C>0$ such that
$|R(\xi,t)|\le Ct^3$ for all $\xi\in [a,1-c]$ and $t$ sufficiently small. 
Moreover, this implies that $R(\xi,t)=o(t^2)$
uniformly with respect to $\xi\in [a,1-c]$, that is
for all $\eps>0$ there exists $\delta=
\delta(\eps)>0$ such that $\left |{R(\xi,t)\over t^2}\right |<\eps$ for all
$(\xi,t)\in [a,1-c]\times [-\delta,\delta]$.

Hence we can choose $\alpha>0$ such that $|R(\xi,t)|\le 
-{\ovl\Psi_\xi^{\prime\prime}(0)\over 4}\,t^2$ (remember that 
$\ovl\Psi_\xi^{\prime\prime}(0)$ exists, is finite and strictly negative
for all $\xi\in [a,1-c]$).

Now we split the integral in \eqref{t-int1} into two parts: the first with
$t$ ranging from $-\alpha$ to $\alpha$ and the second will be the rest:
$$
\eqalign{
 {\rm\bf (A)} & :={1\over 2\pi}\int_{-\alpha}^\alpha G(z(\xi,t))\exp
	\{n\ovl\Psi_\xi(t)\}\d t,\cr
 {\rm\bf (B)} & :={1\over 2\pi}\int_{\alpha<|t|\le\pi}G(z(\xi,t))\exp
	\{n\ovl\Psi_\xi(t)\}\d t.
}$$

{\it Part $I\!I$}

\noindent We estimate {\bf (B)}. Let us note that, by the definitions of
$\phi$ and $z_o(\xi)$, 
$$\exp\{n\phi(\xi)\}=\exp\{n\Psi(z_o(\xi))\}={F(z_0(\xi))^{2k}\over 
z_o(\xi)^n}$$
from which we get $z_o(\xi)^n={F(z_0(\xi))^{2k}\over\exp\{n\phi(\xi)\}}$.
We also note that $\exp\{n\ovl\Psi_\xi(t)\}={F(z(\xi,t))^{2k}\over
z(\xi,t)^n}={F(z(\xi,t))^{2k}\over z_o(\xi)^ne^{int}}$.
Hence we rewrite {\bf (B)}:
$$
%\eqalign{
{\rm\bf (B)} 
%& ={1\over 2\pi}\int_{\alpha<|t|\le\pi}G(z(\xi,t))
%	{F(z(\xi,t))^{2k}\over z_o(\xi)^ne^{int}}\d t=\cr
%	& ={1\over 2\pi}\int_{\alpha<|t|\le\pi}G(z_o(\xi)e^{it})
%	\left({F(z_o(\xi)e^{it})\over F(z_o(\xi))}\right)^{2k}
%	e^{-int}e^{n\phi(\xi)}\d t=\cr
%	& 
	={e^{n\phi(\xi)}\over 2\pi}\int_{\alpha<|t|\le\pi}G(z_o(\xi)e^{it})
	\left({F(z_o(\xi)e^{it})\over F(z_o(\xi))}\right)^{2n\xi}
	e^{-int}\d t.
%}
$$
We want to give an upper estimate for the modulus of the integrand:
we claim that $|F(z)|\le F(|z|)$ and equality holds if and only if $z=|z|$
(apply Proposition~\lemmaref{seriepot}).
%Indeed, we note that in equation \eqref{F2-ol} if $i=1$,
%then $b_0={1\over 2}$, $b_1={3\over
%4}$
%hence applying Proposition~\lemmaref{seriepot} to $\sum_{n=0}^\infty b_{n}
%z^{n}$ we obtain
%$$
%  |F(z)|=|z^{1/2}|\left|\sum_{n=0}^\infty b_{n}z^{n}\right|=|z|^{1/2}
%  \left|\sum_{n=0}^\infty b_{n}z^{n}\right|
%  \le|z|^{1/2}\sum_{n=0}^\infty b_{n}|z|^{n}= F(|z|),
%$$
%where the inequality turns to an equality if and only if $z=|z|$.
Now we observe that $F(z)$ is a continuous mapping and $F(|z|)=0$ if and only
if $z=0$, then ${|F(z)|\over F(|z|)}$ is continuous in the compact set 
$K:=\{z\in\complessi\,:\,|z|\in [\ovl a,1-\ovl c],\alpha\le\arg(z)\le\pi\}$
where it attains a maximum $\lambda<1$.
Then since $\xi<1$,
$$
  |{\rm\bf (B)}|\le{e^{n\phi(\xi)}\over 2\pi}\int_{\alpha<|t|\le\pi}
	|G(z_o(\xi)e^{it})|\lambda^{2n}\d t.
$$
%We claim that 
Moreover, by Proposition~\lemmaref{seriepot}, $|G(z_o(\xi)e^{it})|\le 
G(1-\ovl c)$, from which we get the
upper estimate (uniform with respect to $\xi\in[a,1-c]$):
$$
  |{\rm\bf (B)}|\le e^{n\phi(\xi)}G(1-\ovl c)\lambda^{2n}.
\autoeqno{By1}
$$
%In fact,
%it suffices to apply Proposition~\lemmaref{seriepot} to
%$G(z)$, using equation \eqref{G-ol}
%Moreover, if $|z_1|\le|z_2|$, then
%$G(|z_1|)\le G(|z_2|)$ which implies that $|G(z_o(\xi)e^{it})|\le 
%G(\max_{\xi\in[a,1-c]}(z_o(\xi)))=G(1-\ovl c)$.

{\it Part $I\!I\!I$}

\noindent We estimate {\bf (A)}:
$$
{\rm\bf (A)}={e^{n\phi(\xi)}\over 2\pi}\int_{-\alpha}^\alpha G(z_o(\xi)e^{it})
	\exp\left\{-n{1-\xi^2\over 4\xi^2}t^2+nR(\xi,t)\right\}\d t.
$$
We perform a change of variable in order to stress the main term of the 
exponential: $\theta:=\sqrt{n}b(\xi)t$, where $b(\xi)=\sqrt{1-\xi^2\over 2\xi^2}
=\sqrt{-\ovl\Psi^{\prime\prime}_\xi(0)}$. Hence 
$\d t={n^{-1/2}\d\theta\over b(\xi)}$ and {\bf (A)} becomes
$$
%\eqalign{
{\rm\bf (A)} 
%& ={e^{n\phi(\xi)}\over 2\pi}{n^{-1/2}\over b(\xi)}
% \int_{-\alpha\sqrt{n}b(\xi)}^{\alpha\sqrt{n}b(\xi)} 
%    \exp\left\{-{\theta^2\over 2}+nR(\xi,t_n)\right\}
%    G(z_o(\xi))\left({G(z_o(\xi)e^{it_n})\over G(z_o(\xi))}\right)
%    \d\theta\cr &
   ={e^{n\phi(\xi)}\over 2\pi b(\xi)}G(z_o(\xi))n^{-1/2}
\int_{-\alpha\sqrt{n}b(\xi)}^{\alpha\sqrt{n}b(\xi)}
    \exp\left\{-{\theta^2\over 2}+nR(\xi,t_n)\right\}
    \left({G(z_o(\xi)e^{it_n})\over G(z_o(\xi))}\right)
    \d\theta,
%}
$$
where we put $t_n:=\theta/(\sqrt nb(\xi))$.
We want to give a uniform upper bound for the modulus of the integrand: we
will then be able to apply 
Theorem~\lemmaref{dominunif} . As we noted before, 
$\left|{G(z_o(\xi)e^{it_n})\over G(z_o(\xi))}\right|\le 1$, moreover,
by our choice of $\alpha$, we have that
$|nR(\xi,t_n)|\le n\,{-\ovl\Psi^{\prime\prime}_\xi(0)\over 4}
\,t_n^2=
{\theta^2\over 4}$. Then the modulus of the integrand is bounded by
$\exp\{-\theta^2\!/4\}$ for all $n$, uniformly with respect to $\xi\in [a,1-c]$;
the integrand converges pointwise to $\exp\{-\theta^2\!/2\}$, and the 
interval of integration converges to $\reali$. Then applying 
Theorem~\lemmaref{dominunif}
 and noting that $\int_\reali
\exp\{-\theta^2/2\}=\sqrt{2\pi}$, we have that
$$
 {\rm\bf (A)} \asym{n} {e^{n\phi(\xi)}\over\sqrt\pi}n^{-1/2}
	{\xi\over\sqrt{ 1-\xi^2}}G(z_o(\xi))=
  {e^{n\phi(\xi)}\over\sqrt\pi}n^{-1/2}
      \sqrt{\xi\over {1-\xi^2}}\sqrt{2\over 1+\xi},
$$
uniformly with respect to $\xi\in[a,1-c]$.
To conclude the proof, we only need to show the asymptotical negligibility
of {\bf (B)}, that is, that $\left|{{\rm\bf(B)\over (A)}}\right|\ra 0$ when
$n$ tends to infinity, uniformly with respect to $\xi\in[a,1-c]$.
In fact, by equation~\eqref{By1}
$$
\eqalign{
\left|{{\rm\bf(B)\over (A)}}\right| 
& \le{e^{n\phi(\xi)}\lambda^{2n}G(1-\ovl c)
	\over{\rm\bf (A)}}\asym{n}
	{\sqrt\pi{\sqrt{1-\xi^2}\over\xi}\lambda^{2n}{G(1-\ovl c)\over
	G(z_o(\xi))}n^{1/2}}\le\cr &
	 \le\sqrt\pi n^{1/2} \lambda^{2n}
	{G(1-\ovl c)\over G(\ovl a)}{\sqrt{1-a^2}\over a}
}
$$
and the last term tends to $0$ when $n$ tends to infinity, uniformly with
respect to $\xi\in[a,1-c]$.  \QED
%%%%%%%%%%%%%%%%%%%%%%%%%%%%%%%%%%%%%%%%%%%%%%%%%%%%%%%%%%%%%%%%%%%%%%%%%%%%%%%
%%% FINE DELLA PARTE DI y E \xi\in [a,1-c]
%%%%%%%%%%%%%%%%%%%%%%%%%%%%%%%%%%%%%%%%%%%%%%%%%%%%%%%%%%%%%%%%%%%%%%%%%%%%%%

%%%%%%%%%%%%%%%%%%%%%%%%%%%%%%%%%%%%%%%%%%%%%%%%%%%%%%%%%%%%%%%%%%%%%%%%%
%%%%%%%%%%%%%%%%%%%%%%YPROX0.TEX%%%%%%%%%%%%%%%%%%%%%%%%%%%%%%%%%%%%%%%%%
\autosez{xi<a}{Estimates along the $y$-axis: the case $\xi\in [0,a]$}

If $\xi$ is allowed to tend to zero, the preceding estimate is no longer
true. Then we have to choose a different curve of integration.
We perform a change of variable
$u=f(z)$. The curve of integration will be the union of two pieces: one
will be the inverse image of a suitable curve in the $u$-plane, and the other
will be a part of a circle centered in the origin such that the union of the
two pieces forms a connected circuit about the origin.

We choose $u:=\sqrt{1-z}$, which, by our choice of the determination of the
square root, has argument in
%%%%%=|1-z|^{1/2}\exp (i\arg(1-z)/2)$ where $\arg(1-z)\in
$[-\pi/2,\pi/2)$, and then if $\re(u)<0$ or $u=ib$, $b>0$, then
$u=-\sqrt{1-z}$. Anyway, $z=1-u^2$.
 The desired curve in the $u$-plane is
simply a vertical segment whose parametrization is $u(\xi,t)=u(\xi)-it$ where
$u(\xi):=\sqrt{1-z_o(\xi)}=\xi$ and $t$ ranges from $-\alpha$ to $\alpha$ 
($\alpha$ will be chosen
in the sequel). Note that the segment is oriented downwards in order to
produce a correctly oriented curve in the $z$-plane (see the pictures below)
and that it lies in the half plane where $u=\sqrt{1-z}$.

\beginpicture
%%
%%      DISEGNO SEGMENTO IN U
%%      E CORRISPONDENTE CURVA IN Z
%%
\setcoordinatesystem units <4mm,4mm> point at 0 0
\put{$\,$} at 5 6.5
\put{$\,$} at 5 -7.5
\setlinear
%%%%%PIANO U
\plot 8 0  18 0 /
\plot 9.5 -6    9.5 6 /
\put{0} [rt] at 9 -.5
\plot 12 -4    12 4 /  
\put{$u(\xi)$} [rt] at 11.75 -.5
%freccetta
\plot 12 -3 11.7 -2.3 /
\plot 12 -3 12.3 -2.3 /

%%%%%PIANO Z
\plot  23 0  33 0 /
\plot  24.5 -6    24.5 6 /
\put{0} [rt] at 24 -.5
\put{$\zo$} [rt] at 27.1 -.5
\setquadratic
\plot   31.33739 -4.250004      30.88364 -4.000004      30.45739 -3.750004
	30.05864 -3.500004      29.68739 -3.250004      29.34364 -3.000004
	29.02739 -2.750004      28.73864 -2.500004      28.47739 -2.250004
	28.24364 -2.000004      28.03739 -1.750004      27.85864 -1.500004
	27.70739 -1.250004      27.58364 -1.000004      27.48739 -.7500039
	27.41864 -.5000039      27.37739 -.2500039      
	27.36364 0
	27.37739 .25            27.41864 .5             27.48739 .75
	27.58364 1              27.70739 1.25           27.85864 1.5
	28.03739 1.75           28.24364 2              28.47739 2.25
	28.73864 2.5            29.02739 2.75           29.34364 3
	29.68739 3.25           30.05864 3.5            30.45739 3.75
	30.88364 4              31.33739 4.25           /
%freccetta
\setlinear
\plot 29.34364 3 28.54364 2.8 /
\plot 29.34364 3 29.10364 2.35 /

%%%%%%% stanghetta per 1, piano u
\setlinear
\linethickness=.3pt
\plot 15 -.1  15 .1 /
\put{1}  at 15 -.5
%%%%%%% stanghetta per 1, piano z
\setlinear
\linethickness=.3pt
\plot 28.5 -.1  28.5 .1 /
\put{1}  at 28.5 -.5

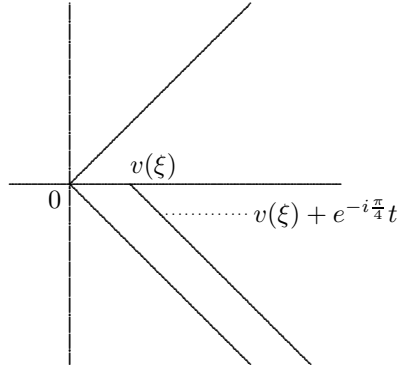
\put{\it Figure~2: the segment in the $u$-plane} [lt] at 8 -6.5
\put{\it Figure~3: the curve in the $z$-plane} [lt] at 23 -6.5
\endpicture

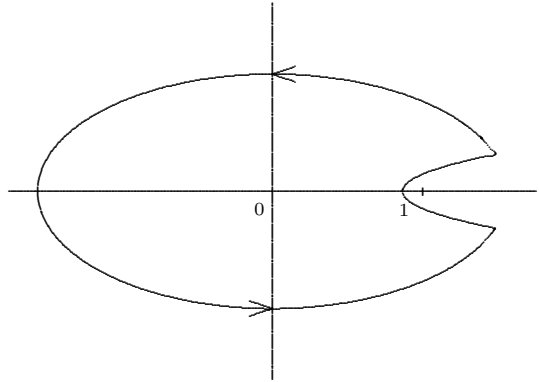
The curve of integration in the $z$-plane will be the union of $z(\xi,t):=
1-u(\xi,t)^2$ for $|t|\le\alpha$ and $\tilde z(\xi,s):=|z(\xi,\alpha)|
e^{is}$ for $\arg(z(\xi,\alpha))\le s\le 2\pi-\arg(z(\xi,\alpha))$
(note that $\arg(z(\xi,\alpha))=-\arg(z(\xi,-\alpha))$).
Hence
$$
p^{(2n)}(2k,0)={\rm\bf(A)}+{\rm\bf(B)}={1\over 2\pi i}
  \int_{\gamma_1}{G(z)F(z)^{2k}\over z^{n+1}}\d z
 +{1\over 2\pi i}\int_{\gamma_2}{G(z)F(z)^{2k}\over z^{n+1}}\d z
 \autoeqno{y-0}
$$
where $\gamma_1$ corresponds to $z(\xi,t)$, $|t|\le\alpha$,              
and $\gamma_2$ corresponds to $\tilde z(\xi,s)$,
$\arg(z(\xi,\alpha))\le s\le 2\pi-\arg(z(\xi,\alpha))$
(and {\bf(A)} and {\bf(B)} are the corresponding integrals).

Here is how the contour of integration appears in the $z$-plane (note
that the circumference is elliptic due to a different choice of measure
units on the horizontal and vertical axes).

%%
%%      FIGURA COL CONTORNO DI INT. IN Z COMPLETO DI CIRCONFERENZA
%%
\beginpicture
\setcoordinatesystem units <5mm,5mm> point at 0 0
\setlinear
\put{$\,$} at 0 5.5
\put{$\,$} at 0 -6  
\put{\it Figure~4: the curve of integration } [lt] at 2 0
\put{\it in the $z$-plane} [lt] at 5 -1
\plot  13 0  27 0 /
\plot  20 -5    20 5 /
\put{$\scriptstyle{0}$} [rt] at 19.8 -.3
\setquadratic
\plot   25.92423 -.98           25.62585 -.91875
	25.34673 -.8575         25.08685 -.7962499      24.84623 -.735
	24.62485 -.6737499      24.42273 -.6124999      24.23985 -.5512499
	24.07623 -.4899999      23.93185 -.4287499      23.80673 -.3674999
	23.70085 -.3062499      23.61423 -.2449999      23.54685 -.1837499
	23.49873 -.1224999      
	23.46023 0
	23.46985 .06125         23.49873 .1225          23.54685 .18375
	23.61423 .245           23.70085 .30625         23.80673 .3675
	23.93185 .42875         24.07623 .49            24.23985 .55125
	24.42273 .6125001       24.62485 .67375         24.84623 .7350001
	25.08685 .7962501       25.34673 .8575          25.62585 .91875
	25.92423 .98      
	25.83191 1.109838       25.58117 1.395403       
%25.69896 1.270822
	25.58117 1.395403       25.41675 1.548946       25.08041 1.811594
	24.69332 2.056142       24.25933 2.280144       23.78278 2.481365
	23.26843 2.657792       22.72143 2.807663       22.14724 2.929482
	21.55159 3.02203        20.94044 3.084382       20.31989 3.115917
	19.69615 3.116318       19.07544 3.085582       18.46397 3.024016
	17.86785 2.932235       17.29303 2.811156       16.74526 2.661989
	16.23001 2.486225       15.75243 2.285619       15.31729 2.062176
	14.92893 1.818128       14.59125 1.555914       14.3076 1.278154
	14.08084 .9876233       13.91321 .6872244       13.80641 .3799591
	13.76148 0.0688973      13.77889 -.2428528      13.85846 -.5521764
	13.9994 -.8559829       14.20028 -1.151237      14.45912 -1.434988
	14.77332 -1.704401      15.13975 -1.956784      15.55473 -2.189616
	16.01413 -2.400571      16.51336 -2.587539      17.04742 -2.748653
	17.61099 -2.882304      18.19843 -2.987156      18.80386 -3.062161
	19.42125 -3.106571      20.04442 -3.11994       20.66714 -3.102136
	21.2832 -3.053337       21.88644 -2.97403       22.47083 -2.865007
	23.03054 -2.727358      23.55996 -2.562459      24.05381 -2.371956
	24.50716 -2.157753      24.91547 -1.921991      25.27467 -1.667025
	25.58117 -1.395403      25.83191 -1.109838      25.92423 -.98     /
\setlinear
\plot 20  3.117 20.6 2.917 / 
\plot 20  3.117 20.6 3.317 /
\plot 20  -3.117 19.4 -2.917 / 
\plot 20  -3.117 19.4 -3.317 /

%%%%%%% stanghetta per 1
\setlinear
\linethickness=.3pt
\plot 24 -.1  24 .1 /
\put{$\scriptstyle{1}$}  at 23.5 -.5

\endpicture

We observe that the integral still makes sense, since $G(z)$ has a unique
singularity in $z=1$ and can be extended to an holomorphic function 
defined in an open set containing the integration domain. Similarly, also
$F(z)^{2k}$ can be extended to an holomorphic function defined in an open
set containing the integration domain (this appears obvious once we look at 
the explicit expressions for these functions instead of their power series
representation).

First we choose $a$ depending on $\alpha$ such that for all $\xi\in[0,a]$
we have $|z(\xi,\alpha)|\ge 1+\eps_o$ for some fixed $\eps_o>0$. Thanks
to this choice, the circular part of the curve of integration  will be
uniformly far from the singularity $z=1$ (that is its distance from 1
is greater than $\eps>0$, not depending on $\xi$).
This choice is possible since the mapping $\xi\mapsto z(\xi,\alpha)$ is
continuous and if $\xi\ra 0$ then $z(\xi,\alpha)\ra 1+\alpha^2$.
Then for all $\eps>0$ there exists a right neighbourhood $\cal U$ of $0$
such that whenever $\xi\in\cal U$ we have $|z(\xi,\alpha)|\ge 1+\alpha^2
-\eps$ and the last quantity is greater than $1+\eps_o$ for $\eps$
sufficiently small.

We observe that in this paragraph we will operate further choices of $a$,
namely $a$ will be chosen sufficiently small in order to satisfy all the
conditions we will find out to be necessary. In the sequel we will not stress
that when a new condition is introduced, if necessary $a$ is chosen
smaller than that of the preceding choice.

Now we estimate the integral {\bf(B)} that will appear to be
asymptotically negligible if compared to the integral {\bf(A)}.

\lemma{int-2-y}{Let $\gamma_2$ be the curve with parametrization $z(t):=
|z(\xi,\alpha)|\exp(it)$ where $|t|\ge\arg(z(\xi,\alpha))$.
Then there exists a sufficiently small $a$ such that
$$
  \left|{1\over 2\pi i}\int_{\gamma_2}{G(z)F(z)^{2k}\over z^{n+1}}\d z\right|
  \le C\,e^{n\phi(\xi)}\,\lambda^n
$$
for some $C>0$, $\lambda<1$ and for all $\xi\in [0,a]$.}

\proof
Let us write the integrand as a function of $t$:
$$
\eqalign{
& {1\over 2\pi i}\int_{\gamma_2}{G(z)F(z)^{2k}\over z^{n+1}}\d z=\cr
 %& ={1\over 2\pi}\int_{|t|\ge\arg(z(\xi,\alpha))}
 %  {G(|z(\xi,\alpha)|e^{it})\,(F(|z(\xi,\alpha)|e^{it})^2)^k\over
 %  |z(\xi,\alpha)|^ne^{itn}}\d t\cr
 & ={e^{n\phi(\xi)}\over 2\pi}
   \int_{|t|\ge\arg(z(\xi,\alpha))}
      {(F(|z(\xi,\alpha)|e^{it}))^{2n\xi}\over F(z_o(\xi))^{2n\xi}}\cdot
   {z_o(\xi)^n\over |z(\xi,\alpha)|^n}\cdot
   {G(|z(\xi,\alpha)|e^{it})\over e^{itn}}\d t}
   \autoeqno{trasc-y}
$$
where
% in the last equality 
we used $F(z_o(\xi))^{2k}=e^{n\phi(\xi)}\cdot
z_o(\xi)^n$.

We want to give an upper bound for $|G|$ and $|F^2|$, uniformly with respect
to $\xi\in [0,a]$ and $|t|\ge\arg(z(\xi,\alpha))$ (this upper bound
will depend on $a$).

This is possible since $\bigcup_{\xi\in [0,a]}\{(\xi,t)\,:\,|t|\ge
\arg(z(\xi,\alpha))\}$ is a compact subset of $\reali^2$ and the mapping
$(\xi,t)\mapsto |z(\xi,\alpha)|e^{it}$ is continuous, hence
$K:=\{z=|z(\xi,\alpha)|e^{it}\,:\,\xi\in [0,a],|t|\ge\arg
(z(\xi,\alpha))\}$ is a compact subset of $\complessi$.

%%
%%      FIGURA CON IL COMPATTO DI PAGINA 10
%%
\beginpicture
\setcoordinatesystem units <4mm,4mm> point at 0 0
\setlinear
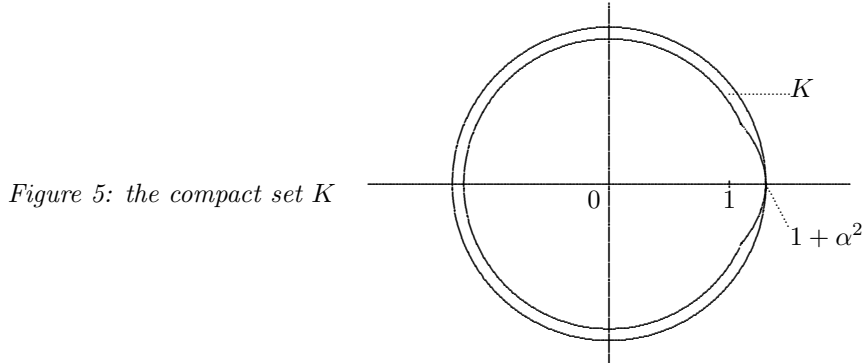
\put{\it Figure~5: the compact set $K$} [lt] at -20 0
\put{$\,$} at -22 6.5
\put{$\,$} at -22 -6
\plot  -8 0  8 0 /
\plot  0 -6    0 6 /
\put{0}  at -.5 -.5
\put{$1+\alpha^2$} [l] at 6 -1.7
\put{$K$} [l] at 6 3.2

\setquadratic
%curva da xi=0 a xi=.5  
\plot  5.21  0  5.2  .22        5.17  .44       5.12  .66
5.05  .8800001  4.96  1.1       4.85  1.32      4.72  1.54
4.57  1.76      4.4  1.98       4.4  1.98
% arco di circonferenza
4.180348  2.409376      3.918928  2.814677      3.618351  3.191856
3.28162  3.537142       2.9121  3.847086        2.513484  4.118592
2.089754  4.348946      1.645144  4.535847      1.184096  4.677427
.7112164  4.772271      .231231  4.819433       
-.2510647  4.81844
-.7308519  4.769304     -1.203337  4.672513     -1.663798  4.529037
-2.107635  4.340308     -2.530413  4.108212     -2.927908  3.835069
-3.296149  3.523607     -3.631455  3.176938     -3.930477  2.798526
-4.190228  2.392153     -4.408111  1.961878     -4.581949  1.512
-4.710007  1.047015     -4.791003  .571569      -4.82413  .0904117
-4.809055  -.3916488
-4.745931  -.8697962    -4.635386  -1.339253    -4.478527  -1.795328
-4.276919  -2.233465    -4.032578  -2.649286    -3.747945  -3.038636
-3.425864  -3.397626    -3.069552  -3.722667    -2.682571  -4.010513
-2.268787  -4.258287    -1.832333  -4.463514    -1.377572  -4.624143
-.9090461  -4.738569    -.4314376  -4.805649    
.0504817  -4.824713
.5318967  -4.795569     1.007997  -4.718511     1.474026  -4.594306
1.925327  -4.424197     2.35739  -4.209883      2.7659  -3.953505
3.146773  -3.657625     3.496206  -3.325199     3.810705  -2.959549
4.087129  -2.564328     4.322716  -2.143486     4.4  -1.98
% curva con xi da .5 a 0 e y negative
4.4  -1.98      4.57  -1.76     4.72  -1.54     4.85  -1.32
4.96  -1.1      5.05  -.8799998         5.12  -.6599998
5.17  -.4399998         5.2  -.2199998  5.21  0
% circonferenza di raggio 1+alpha^2 (alpha=.55)
5.183972  .5201321      5.106147  1.035067      4.977303  1.53966
4.798728  2.02887       4.572205  2.497807      4.299999  2.941787
3.984828  3.356374      3.629842  3.737426      3.238588  4.081133
2.814975  4.384064      2.363235  4.643191      1.887883  4.855924
1.393668  5.020138      .8855277  5.134193      .3685396  5.196949
-.1521309  5.207778
-.6712813  5.166574     -1.183725  5.073746     -1.68434  4.930223
-2.168126  4.737439     -2.630249  4.497321     -3.066091  4.212266
-3.471298  3.885124     -3.841821  3.519164     -4.173957  3.118041
-4.464389  2.685764     -4.710215  2.226651     -4.908978  1.745291
-5.058692  1.246492     -5.15786  .7352389      -5.205494  .2166395
-5.201116  -.3041246    
-5.144771  -.8218499    -5.03702  -1.331364     -4.878942  -1.827575
-4.672114  -2.305526    -4.418605  -2.76044     -4.120946  -3.187773
-3.782113  -3.583256    -3.40549  -3.942935     -2.994839  -4.263219
-2.554266  -4.540905    -2.088172  -4.773221    -1.601213  -4.957844
-1.098256  -5.09293     -.584325  -5.177129     -.0645558  -5.2096
.4558583  -5.190019     .9717177  -5.11858      1.477868  -4.995999
1.969252  -4.823499     2.44096  -4.602805      2.888279  -4.336121
3.306739  -4.026112     3.692159  -3.675876     4.040689  -3.288911
4.348845  -2.869085     4.613549  -2.420592     4.832156  -1.947913
5.002482  -1.455771     5.122825  -.9490842     5.191983  -.4329141
5.21  0 /
%%%%%%% stanghetta per 1
\setlinear
\linethickness=.3pt
\plot 4 -.1  4 .1 /
\put{1}  at 4 -.5

\setdots <1.5pt>
\plot 5.21 0 6 -1.5 /
\plot  4 3  6 3 /
\endpicture

Moreover $d(K,1)\ge \eps_o$ and $G$ and $F^2$ are holomorphic functions
in an open domain containing $K$, therefore there exists $C\ge1$ (depending on
$a$) such that
$$
\max_{z\in K}(|F^2(z)|,|G(z)|)\le C.\autoeqno{C(a)}
$$
Hence the modulus of the integrand in equation~\eqref{trasc-y} is bounded
by
$$
  \left({C\over F(z_o(a))}\right)^{na}\cdot
  \left({1\over 1+\eps_o}\right)^{n}\cdot C.
$$
We observe that we may write $C=C(a)$ and if we suppose to choose $C(a)$ to
be the smallest possible $C\ge 1$ satisfying \eqref{C(a)}, then $C(a)$ turns
out to be a continuous increasing function (since as $a$ increases, $K=K(a)$
becomes a larger subset of $\complessi$).
We want to show that, if $a$ is small enough, then
$$
 \lambda:= \left({C(a)\over F(z_0(a))}\right)^{a}\cdot\left({1\over 1+\eps_o}
 \right)<1.\autoeqno{lambda}
 $$

We note that $a\mapsto\left({C(a)\over F(z_0(a))}\right)^{a}$ is a continuous
increasing function, not smaller than $1$, and $\lim_{a\ra 0}z_o(a)=1$,
$F$ is continuous and $F(1)=1$, hence
$$
  \lim_{a\ra 0}\left({C(a)\over F(z_0(a))}\right)^{a}=1.
$$
As ${1\over 1+\eps_o}<1$, $a$ can be chosen small enough in order to
satisfy \eqref{lambda}.
The thesis follows using this estimate. \QED

In order to estimate the integral {\bf(A)}, we rewrite it stressing
the exponential part of the integrand:
$$
{\rm\bf(A)}
%={1\over 2\pi i}\int_{\gamma_1}{G(z)F(z)^{2k}\over z^{n+1}}\d z
={1\over \pi}\int^{\alpha}_{-\alpha}
\exp\{n\ovl\Psi_\xi(t)\}{G(z(\xi,t))\over z(\xi,t)}(u(\xi)-it)\d t
$$
where   $z(\xi,t)=1-(u(\xi)-it)^2$ and $\ovl\Psi_\xi(t):=\Psi_\xi(z(\xi,t))$.

\lemma{psit}{The function $\ovl\Psi_\xi(t)$ has a Taylor expansion centered
in $0$:
$$
  \ovl\Psi_\xi(t)=\phi(\xi)-{1\over 2}\left({2\over 1-u(\xi)^2}\right)t^2
  +R(\xi,t)
$$
where $|R(\xi,t)|\le C|t|^3$ for all $\xi\in [0,a]$, $|t|\le\alpha$, and for 
some $C>0$. ($C$ depends on $a$ and $\alpha$).}

\proof
We calculate the first derivative of $\ovl\Psi_\xi$: its Taylor series
will lead us to the Taylor expansion of the primitive function:
$$
   \ovl\Psi_\xi^\prime(t) 
	 =-i\left\{ {1\over 1+i{t\over 1-u(\xi)} }-
	 {1\over 1-i{t\over 1+u(\xi)}}
	\right\}=
    \sum_{n\ge 1}\left[{(-i)^{n+1}\over(1-u(\xi))^n}+
    {(i)^{n+1}\over(1+u(\xi))^n}\right]\,t^n,
$$
where the series converges to the function itself,
%This function has a Taylor series expansion (converging to the function
%itself), 
provided that $\left|{t\over 1\pm u(\xi)}\right|<1$ 
(but this is true for $a$ and $\alpha$
sufficiently small).
%$$
%  \eqalign{
%   \ovl\Psi_\xi^\prime(t) & =-i\sum_{n\ge 0}{(-it)^n\over(1-u(\xi))^n}+
%   i\sum_{n\ge 0}{(it)^n\over(1+u(\xi))^n}=\cr
%   &=\sum_{n\ge 0}\left[{(-i)^{n+1}\over(1-u(\xi))^n}+
%    {(i)^{n+1}\over(1+u(\xi))^n}\right]\,t^n=
%    \sum_{n\ge 1}\left[{(-i)^{n+1}\over(1-u(\xi))^n}+
%    {(i)^{n+1}\over(1+u(\xi))^n}\right]\,t^n
%}$$
%from which we derive the Taylor series of $\ovl\Psi_\xi(t)$:
Then
$$
   \ovl\Psi_\xi(t) =
    \sum_{n\ge 1}\left[{(-i)^{n+1}\over(1-u(\xi))^n}+
    {(i)^{n+1}\over(1+u(\xi))^n}\right]\,{t^{n+1}\over n+1}+\phi(\xi)
$$
where we added the constant $\phi(\xi)=\ovl\Psi_\xi(0)$,
and
%Then
$$
   \ovl\Psi_\xi(t) 
     =\phi(\xi)-{1\over 1-u(\xi)^2}t^2+{4\over 3}i{u(\xi)\over(1-u(\xi)^2
    )^2}t^3+ t^4\,S(\xi,t)
$$
where $S(\xi,t)=
    \sum_{n\ge 4}{1\over n}\left[{(-i)^{n}\over(1-u(\xi))^{n-1}}+
    {(i)^{n}\over(1+u(\xi))^{n-1}}\right]\,t^{n-4}$.
We want to show that $|S(\xi,t)|\le C$  where $C>0$ does not
depend on $\xi\in [0,a]$ nor on $t\in[-\alpha,\alpha]$. This will conclude the
proof, since
$$
R(\xi,t):=\left[{4\over 3}i{u(\xi)\over (1-u(\xi)^2)^2}+tS(\xi,t)\right]\,t^3.
$$
%and $\left|{4\over 3}i{u(\xi)\over(1-u(\xi)^2)^2}\right|\le 
%{4\over 3}{a\over(1-a^2)^2}$. Then let
%us estimate $|S(\xi,t)|$:
But
$$
\eqalign{
	|S(\xi,t)| & \le
   % \sum_{n= 4}^\infty{1\over n}\left[{1\over(1-u(\xi))^{n-1}}+
   % {1\over(1+u(\xi))^{n-1}}\right]\,|t|^{n-4}\le\cr
    \sum_{n= 4}^\infty{1\over n}\left[{1\over(1-a)^{n-1}}+
    1\right]\,|t|^{n-4}
%  = \sum_{n= 4}^\infty{1\over n}\left[{1+(1-a)^{n-1}\over(1-a)^{n-1}}
%    \right]\,|t|^{n-4} \le   \cr
%    & \le
%    \sum_{n= 4}^\infty{1\over n}\left({2-a\over 1-a}\right)^{n-1}\,|t|^{n-4}
%    \le  \sum_{n= 4}^\infty{1\over n}
%    \left({2-a\over 1-a}\right)^{n-1}\,\alpha^{n-4}\le\cr
%    & =
\le\sum_{n= 0}^\infty \left(
    {2-a\over 1-a}\right)^{n+3}\,{\alpha^{n}\over n+4}.
    }$$
%where for the second inequality we used $1-u(\xi)\ge 1-a$ and $1+u(\xi)\ge 1$
%and for the third inequality we applied $x^n+y^n\le (x+y)^n$ for all
%$x,y\ge 0$. 
Now the last term is equal to a constant $C=C(a,\alpha)$
not depending on $\xi$, provided that the power series converges in
$\alpha$. 
%Then we evaluate its radius of convergence $\rho$:
%$$
%  1/\rho=\limsup_{n\ra +\infty}{1\over (n+4)^{1\over n}}
%  \,\left({2-a\over 1-a}
%  \right)^{1+{3\over n}} ={2-a\over 1-a}\le 4
%$$
Moreover the radius $\rho$ of this power series is ${2-a\over 1-a}\le 4$
%where the last inequality holds 
if $a\le {1\over 2}$ and therefore
$\rho\ge {1\over 4}$. Choosing $\alpha\le {1\over 4}$ we get to the
conclusion. \QED

The asymptotic estimate turns out to be different whether $\xi$ is allowed
to tend very fast to $0$ or not, that is we have to distinguish two subcases.

\theorem{y-xilento}{If $\xi\in[n^{-1/4}, a]$ for some $a\in(0,1)$, then
$$
  p^{(2n)}(2k,0)\asym{n}\sqrt{2\over\pi}{e^{n\phi(\xi)}\sqrt\xi\over
  \sqrt{(1+\xi)(1-\xi^2)} }\,n^{-1/2}
$$
uniformly with respect to $\xi\in [n^{-1/4},a]$.}

\proof
We rewrite the integral {\bf(A)} using the Taylor expansion of
$\ovl\Psi_\xi(t)$:
$$
{\rm\bf(A)}={e^{n\phi(\xi)}\over \pi}\int^{\alpha}_{-\alpha}
{G(z(\xi,t))\over z(\xi,t)}
\exp\left\{n\left[-{1\over 2}\left({2\over 1-u(\xi)^2}\right)t^2+R(\xi,t)
\right]\right\}(u(\xi)-it)\d t. 
\autoeqno{y-inta}
$$
Now we perform a change of variable in order to stress the main term of the
exponential: $\theta:=\sqrt{n}b(\xi)t$, 
where $b(\xi)=\sqrt{2\over 1-u(\xi)^2}$.
Moreover we denote by $t_n=\theta/(\sqrt nb(\xi))$. 
 Then we choose $\alpha$ (possibly smaller
than the preceding choice) such that
$$
  |R(\xi,t)|\le{1\over 4}b(\xi)^2t^2\qquad\forall |t|\le\alpha
  \autoeqno{sceltaalfa}
$$
uniformly with respect to $\xi\in[0,a]$.
This is possible thanks to a similar argument as in the case $\xi>a$ (see
Theorem~\lemmaref{y-xinonzero} Part I), noting that $b(\xi)\ge b(0)=\sqrt 2$.
After the change of variable
$$
\eqalign{
{\rm\bf(A)} 
%& ={e^{n\phi(\xi)}\over\pi\sqrt nb(\xi)}
%  \int_{-\alpha\sqrt nb(\xi)}^{\alpha\sqrt nb(\xi)}
%  {G\left(z\left(\xi,t_n\right)\right)\over
%	  z\left(\xi,t_n\right)}\cdot\cr
%     &\qquad\qquad\qquad\qquad \cdot\exp\left\{-{1\over 2}\theta^2+n
%     R\left(\xi,t_n\right)\right\}\,
%     \left(u(\xi)-i
%     t_n\right)\d\theta=\cr 
& ={e^{n\phi(\xi)}\over\pi\sqrt nb(\xi)}{G(z_o(\xi))\over z_o(\xi)}
  \int_{-\alpha\sqrt nb(\xi)}^{\alpha\sqrt nb(\xi)}
  {G\left(z\left(\xi,t_n\right)\right)\over
  G(z_o(\xi))}{z_o(\xi)\over
	  z\left(\xi,t_n\right)}\cdot\cr
   & \qquad\qquad\qquad\qquad\qquad \cdot \exp\left\{-{1\over 2}\theta^2+n
     R\left(\xi,t_n\right)\right\}\,
     \left(u(\xi)-i
     t_n\right)\d\theta.
}$$
We want to give a uniform (with respect to $\xi\in [n^{-1/4},a]$) upper bound
for the modulus of the integrand
in order to apply Theorem~\lemmaref{dominunif}.
The exponential part is easily bounded:
$$
\left|\exp\left\{-{1\over 2}\theta^2+n
     R\left(\xi,t_n\right)\right\}\right|\le\exp\left\{
     -{1\over 2}\theta^2+n
     {b(\xi)^2\over 4}t_n^2\right\}
     =\exp\{-\theta^2\!/4\},
     \autoeqno{y2-1}
$$
while $u(\xi)-it_n\,{\buildrel n\over\sim}\,u(\xi)=\xi$
uniformly with respect to $\xi\in [n^{-1/4},a]$, for every fixed $\theta$,
as $n$ tends to infinity.
%In fact
%$$
%u(\xi)-it_n=u(\xi)\left(1-i{\theta\over u(\xi)
%  \sqrt nb(\xi)}\right)
%$$
%and
%$$
%\left|{\theta\over u(\xi)\sqrt nb(\xi)}\right|\le {|\theta|\over
%  \sqrt 2n^{1/4}}                \autoeqno{y2-2}
%$$
%where we used $u(\xi)=\xi\ge n^{-1/4}$, $b(\xi)\ge\sqrt 2$, and the last
%quantity tends uniformly to $0$ (with respect to $\xi\in [n^{-1/4},a]$)
%for every fixed $\theta$, as $n$ tends to infinity.

Let us note that $z_o(\xi)=|z_o(\xi)|\le|z(\xi,t)|$ for all
$|t|\le\alpha$, $\xi\in [n^{-1/4},a]$, hence
$$
  \left|{z_0(\xi)\over z(\xi,t_n)}\right|\le 1   \autoeqno{y2-3}
$$
uniformly with respect to $\xi\in [n^{-1/4},a]$.
Now we only have to evaluate $\left|{G(z(\xi,t_n)\over G(z_o(\xi))}\right|$.
But, using $u(\xi)=\xi$, we have 
$$
\left|{G(z(\xi,t_n))\over G(z_o(\xi))}\right|=
{\sqrt{\xi(1+\xi)}\over |\sqrt{(\xi-it)(1+\xi-it)}|}\le 1
\autoeqno{y2-4}
$$
for all $\xi\in[0,a]$, $|t|\le\alpha$.
Then
$$
\eqalign{
{\rm\bf(A)} & =
 {e^{n\phi(\xi)}u(\xi)\over\pi\sqrt nb(\xi)}{G(z_o(\xi))\over z_o(\xi)}
  \int_{-\alpha\sqrt nb(\xi)}^{\alpha\sqrt nb(\xi)}
  {G(z(\xi,t_n))\over
  G(z_o(\xi))}{z_o(\xi)\over
	  z(\xi,t_n)}\cdot\cr
   & \qquad\qquad\qquad\qquad\qquad \cdot \exp\left\{-{1\over 2}\theta^2+n
     R(\xi,t_n)\right\}\,
     \left(1-i
     {\theta\over u(\xi)\sqrt nb(\xi)}\right)\d\theta
}
$$
where the modulus of the integrand is uniformly bounded 
with respect to $\xi\in [n^{-1/4},a]$ by $\exp
\{-\theta^2\!/4\}(1+c|\theta|)$ (use \eqref{y2-1}, \eqref{y2-3}, 
\eqref{y2-4}) 
and the integrand converges pointwise to
$e^{-{1\over 2}\theta^2}$.
% (apply \eqref{y2-2} and note that $t_n$ tends 
%to zero). 
Then we apply 
Theorem~\lemmaref{dominunif} and we obtain
$$
{\rm\bf(A)}\,{\buildrel n\over\sim}\,
 {e^{n\phi(\xi)}u(\xi)\over\pi\sqrt nb(\xi)}{G(z_o(\xi))\over z_o(\xi)}
 \int_\reali e^{-{1\over 2}\theta^2}\d\theta=
\sqrt{2\over\pi}{e^{n\phi(\xi)}\sqrt\xi\over
  \sqrt{(1+\xi)(1-\xi^2)} }\,n^{-1/2}
$$
uniformly with respect to $\xi\in [n^{-1/4},a]$.
To prove the thesis, we only have to show asymptotic negligibility of
{\bf (B)} with respect to {\bf (A)}. This follows from 
Lemma~\lemmaref{int-2-y}:
$$
  \left|{{\rm\bf(B)}\over{\rm\bf(A)}}\right|\le
  { e^{n\phi(\xi)}\lambda^nC\over 
 \sqrt{2\over\pi}{e^{n\phi(\xi)}\sqrt\xi\over
  \sqrt{(1+\xi)(1-\xi^2)}} }\,n^{-1/2}\le C\lambda^n n^{1/2+1/8}
$$
where the last term tends to zero as $n$ tends to infinity, uniformly with 
respect to $\xi\in [n^{-1/4},a]$. \QED

When $\xi\in[0,n^{-1/4}]$ it is no longer true that $u(\xi)-it_n=\xi-it_n\,
\sim\,\xi$ and the technique will be slightly different. A useful tool
will appear to be the computation of the real and imaginary parts of
$\sqrt{\xi-it_n}$.

\lemma{xi-it}{Let $\sqrt{\xi-it}=a(\xi,t)+ib(\xi,t)$. Then
$$
  a(\xi,t)=\sqrt{{\sqrt{\xi^2+t^2}+\xi\over 2}},\qquad\qquad
  b(\xi,t)={\rm sign}(-t)\sqrt{{\sqrt{\xi^2+t^2}-\xi\over 2}}
$$
and both these terms are $O(\sqrt\xi)+O(\sqrt{|t|})$.}

%\proof
%We first establish the sign of $a$ and $b$: by our choice of the
%determination of the square root, $\arg(\sqrt{\xi-it})$ lies in $[-\pi/2,0)$ 
%if $t>0$ and in $[0,\pi/2)$ otherwise. Hence $a\ge 0$ for all $\xi$ and
%$t$, and ${\rm sign}(b)=-{\rm sign}(t)$.
%
%Elementary computations give us the explicit expressions for $a$ and $b$.
%To estimate these quantities, we study them first for $\xi\le |t|$, 
%then for $\xi>|t|$.
%In the first case $\xi^2+t^2\le 2t^2$ and it is not difficult to show
%that both $a$ and $|b|$ are less or equal to $\sqrt{\sqrt 2+1\over 2}
%\sqrt{|t|}$. In the other case, 
% $\xi^2+t^2\le 2\xi^2$ and both $a$ and $|b|$ are less or equal to 
% $\sqrt{\sqrt 2+1\over 2}\sqrt\xi$. In any case, $a$ and $|b|$ are
% less or equal to $C(\sqrt{|t|}+\sqrt\xi)$. \QED
%%%%%%%%%%%%%%%%%%%%%%%%%%%%%%%%%%%%%%%%%%%%%%%

\theorem{y-xiveloce}{If $\xi\in[0,n^{-1/4}]$ then
$$
  p^{(2n)}(2k,0)\asym{n}{e^{n\phi(\xi)}\over
  \sqrt 2\pi }I\left(\sqrt n\xi\right)\,n^{-3/4},
$$
where $I(t):={\displaystyle{\int_\reali}}
 e^{-{\theta^2\over 2}}\sqrt{\sqrt{t^2+{\theta^2
\over 2}}+t}\,\d\theta$, and the estimate is
uniform with respect to $\xi\in [0,n^{-1/4}]$.
Moreover, if $\xi\in[0,n^{-1/2-\eps}]$ for some $\eps>0$, then
$$
  p^{(2n)}(2k,0)\asym{n}{\sqrt2\,e^{n\phi(\xi)}\over
  \Gamma\left({1\over 4}\right)}\,n^{-3/4},
$$
uniformly with respect to $\xi\in [0,n^{-1/2-\eps}]$.}

\proof
The integral we have to estimate is still {\bf(A)} of equation \eqref{y-inta},
$\alpha$ is chosen as in Theorem~\lemmaref{y-xilento} and we perform the 
same  change of variable $\theta:=\sqrt nb(\xi)t$.
Then we proceed 
differently: in Part $I$ we show a decomposition of $G(z)$ into its singular
and regular parts; Part $I\! I$ is devoted to the estimate of the real part
of the integrand in {\bf (A)}: we stress only the terms which are not
$o(\xi)$ or $o(t_n^2)$ (uniformly with respect to $\xi\in[0,n^{-1/4}]$).
In Part $I\! I\! I$ we write ${\rm\bf (A)}={\rm\bf (A_1)}+{\rm\bf (A_2)}+
{\rm\bf (A_3)}$.
In Part $I\!V$ we estimate ${\rm\bf (A_1)}$ and describe some properties of $I(t)$.
Finally, in Part $V$ we show negligibility of {\bf (B)}, ${\rm\bf (A_2)}$ and
${\rm\bf (A_3)}$.

{\it Part $I$}

\noindent We write a decomposition for $G(z)$:
$$
  \eqalign{
   G(z) &= 
	{\sqrt 2\over (1-z)^{1/4}}
   (1+\sqrt{1-z})^{-1/2}=
   {\sqrt 2\over (1-z)^{1/4}}\cdot\sum_{n=0}^\infty\pmatrix{-1/2\cr n}
   (1-z)^{n/2}=\cr
  & = (1-z)^{-1/4}H(z)+(1-z)^{1/4}K(z)
  }     \autoeqno{decG}
$$
where $H(z):=\sqrt 2\displaystyle{\sum_{n=0}^\infty\pmatrix{-1/2\cr 2n}(1-z)^{n}}$ and
$K(z):=\sqrt 2\displaystyle{\sum_{n=0}^\infty\pmatrix{-1/2\cr 2n+1}(1-z)^{n}}$ are two
holomorphic functions defined in a disc centered in $z=1$ and with 
greater or equal to $2/3$ radius.
Note that the expansion of $(1+\sqrt{1-z})^{-1/2}$ as a power series is 
admissible provided that $|1-z|<1$, and this condition surely holds in the integration
domain if $\alpha$ is sufficiently small.

We decompose $H$ and $K$ into their real and imaginary parts: $H=H_0+iH_1$ and
$K=K_0+iK_1$. The following properties hold: $H(1)=H_0(1)=\sqrt 2$,
$K(1)=K_0(1)=-1/\sqrt 2$, $H_1(z)=O(|1-z|)$, $K_1(z)=O(|1-z|)$ and 
$O(|1-z|)=O(\xi^2)+O(t_n^2)$. 
%Let us prove it for $H_1(z)$:
%$$
%  {|H_1(z)|\over |1-z|}\le\sqrt 2\sum_{n=1}^\infty\pmatrix{-1/2\cr 2n}
%  |1-z|^{n-1}={3\sqrt 2\over 8} +|1-z|\sqrt 2\sum_{n=1}^\infty\pmatrix{-1/2\cr 2(n+1)}
%  |1-z|^{n-1}
%$$
%which proves the estimate for $H_1(z)$, similarly one can prove the analogous for
%$K_1(z)$. 
%Now we rewrite $O(|1-z|)$ as a 
%function of $\xi$ and $t_n$, since $1-z=(u(\xi)-it_n)^2=(\xi-it_n)^2$:
%$$
%  |1-z|=\xi^2+t_n^2\ \Ra\ O(|1-z|)=O(\xi^2+t_n^2)=O(\xi^2)+O(t_n^2).
%$$

{\it Part $I\! I$}

\noindent We can write {\bf (A)} as follows:
$$
  \eqalign{
   {\rm\bf(A)} %& ={e^{n\phi(\xi)}\over\pi\sqrt nb(\xi)}
%   \int_{-\alpha\sqrt nb(\xi)}^{\alpha\sqrt nb(\xi)}
%   \exp\left\{-{1\over 2}\theta^2+nR(\xi,t_n)\right\}{\xi-it_n\over z(\xi,t_n)}
%   \cdot\cr
%   & \qquad\qquad\qquad\qquad\cdot\left({1\over\sqrt{\xi-it_n}}H(z(\xi,t_n))+
%   \sqrt{\xi-it_n}K(z(\xi,t_n))
%   \right)\d\theta\cr
   & ={e^{n\phi(\xi)}\over\pi\sqrt nb(\xi)}
   \int_{-\alpha\sqrt nb(\xi)}^{\alpha\sqrt nb(\xi)}
   \exp\left\{-{1\over 2}\theta^2+nR(\xi,t_n)\right\}{\sqrt{\xi-it_n}
   \over z(\xi,t_n)}\cdot\cr
   & \qquad\qquad\qquad\qquad\cdot\left(H(z(\xi,t_n))+
    (\xi-it_n)K(z(\xi,t_n))\right)\d\theta.
  }
$$
Since the transition probabilities are non negative quantities, 
we are interested
only in the real part of the last integral (in fact the imaginary part
must be $0$, but taking into account only the real part of the integrand
avoids useless computation).

Our aim is to apply Theorem~\lemmaref{dominscon1}, then we want to estimate
the function (depending on $\xi$, $n$ and $\theta$) to which the real part
of our integrand is asymptotic (uniformly with respect to $\xi\in[0,
n^{-1/4}]$). Hence we rewrite the integrand using $\sqrt{\xi-it_n}=
a(\xi,t_n)+ib(\xi,t_n)$ (see Lemma~\lemmaref{xi-it}):
$$
\eqalign{
   & \ident_{[-\alpha\sqrt nb(\xi),\alpha\sqrt nb(\xi)]}(\theta)\cdot
   \exp\left\{-{1\over 2}\theta^2+nR_0(\xi,t_n)\right\}
   \exp\left\{inR_1(\xi,t_n)\right\}\cdot\cr
   & \qquad\qquad\cdot{1\over z_n}
   (a(\xi,t_n)+ib(\xi,t_n))\left(H(z_n)+
    (a(\xi,t_n)+ib(\xi,t_n))^2K(z_n)\right),
}
$$
where $R_0(\xi,t_n)$ and $R_1(\xi,t_n)$ are respectively the real and 
imaginary part of $R(\xi,t_n)$ and $z_n:=z(\xi,t_n)$.

First we want to estimate  the main term of
$$
  \re\left\{\exp\left\{inR_1(\xi,t_n)\right\}
   (a(\xi,t_n)+ib(\xi,t_n))\left(H(z_n)+
    (a(\xi,t_n)+ib(\xi,t_n))^2K(z_n)\right)
  \right\}. \autoeqno{prod1}
$$
It is useful to note that 
$$
nR_1(\xi,t_n)=t_n(O(\xi)+O(t_n))
\autoeqno{y-nr1}
$$
for every fixed $\theta$, uniformly with respect to $\xi\in[0,n^{-1/4}]$.
In fact, 
$$
R_1(\xi,t)=\left[{4\over 3}{u(\xi)\over (1-u(\xi))^2}+
t\im S(\xi,t)\right]t^3=[O(\xi)+O(t)]t^3
\autoeqno{y-r1}
$$
since 
$$
S(\xi,t)={1\over 4}\left[(1-u(\xi))^{-3}+(1+u(\xi)^{-3}\right]+
    t\sum_{n\ge 0}{1\over n+5}\left[{(-i)^{n+1}\over(1-u(\xi))^{n+4}}+
    {(i)^{n+1}\over(1+u(\xi))^{n+4}}\right]\,t^{n}
$$
and hence $\im S(\xi,t)=O(t)$ uniformly with respect to $\xi\in[0,n^{-1/4}]$.

Then from \eqref{y-r1} we get
$$
nR_1(\xi,t_n)=nt_n^2\cdot t_n(O(\xi)+O(t_n))={\theta^2\over b(\xi)^2}
t_n(O(\xi)+O(t_n))
$$ 
which leads to \eqref{y-nr1} since $b(\xi)\ge\sqrt 2$ for all $\xi\in[0,a]$.

Now we expand the product in \eqref{prod1}, where we distinguish the
real and imaginary parts of $H$ and $K$ and write $a_n$ and $b_n$
for $a(\xi,t_n)$ and $b(\xi,t_n)$ respectively.
We claim that this product can be written as follows:
$$
%  \eqalign{
%  & \re\{(\cos(n\r1)+i\sin(n\r1))(a_n+ib_n)\cdot\cr
%  &\cdot(H_0(z_n)+iH_1(z_n)+(\xi-it_n)(K_0(z_n)+iK_1(z_n)))\}=\cr
%  & =\cos(n\r1)\{a_nH_0(z_n)+a_n\xi K_0(z_n)+a_nt_nK_1(z_n)+\cr
%  &\quad -b_nH_1(z_n)-b_n\xi K_1(z_n)+b_nt_nK_0(z_n)\}+\cr
%  & \quad -\sin(n\r1)\{a_nH_1(z_n)+a_n\xi K_1(z_n)-a_nt_nK_0(z_n)+\cr
%  & \quad +b_nH_0(z_n)+b_n\xi K_0(z_n)+b_nt_nK_1(z_n)\}=\cr
%  & =
\cos(n\r1)a_nH_0(z_n)+\cos(n\r1)b_nt_nK_0(z_n)+o(\xi)+o(t_n^2),
%}
\autoeqno{y-est-re}$$
%We have to prove the last equality, 
where $o(\xi)$ and $o(t_n^2)$ are uniform with respect to $\xi\in[0,n^{-1/4}]$.

The proof of \eqref{y-est-re} is tedious but straightforward; particular
care should be put only in estimating the following terms:
$$\eqalign{
	\xi\cdot  O(\sqrt{|t_n|}) & =o(\xi),\cr
	O(\sqrt\xi)\cdot O(t_n^2) & =o(t_n^2),\cr
	O(\xi)\cdot O(|t_n|^{3/2}) & =o(\xi),}\autoeqno{opiccoli2}
$$
and these estimates are uniform with respect to $\xi\in[0,n^{-1/4}]$,
since %Now we remark that 
we consider $n$ tending to infinity, which implies
that both $\xi$ and $t_n$ tend to $0$ (for every fixed $\theta$).
%, uniformly
%with respect to $\xi\in[0,n^{-1/4}]$), hence
%$$\eqalign{
%	\xi\cdot  O(\sqrt{|t_n|}) & =o(\xi)\cr
%	O(\sqrt\xi)\cdot O(t_n^2) & =o(t_n^2)\cr
%	O(\xi)\cdot O(|t_n|^{3/2}) & =o(\xi).}\autoeqno{opiccoli2}
%$$

Finally, we note that $z_n^{-1}=1+O(\xi^2)+O(t_n^2)$ for every fixed $\theta$, 
uniformly with respect to $\xi\in[0,n^{-1/4}]$:
$$
  \left|{1\over z_n}-1\right|={|\xi-it_n|^2\over |1-(\xi-it_n)^2|}
	\le C\,(\xi^2+t_n^2).
	\autoeqno{y-1/zn}
$$
Then by \eqref{y-est-re} and \eqref{y-1/zn} we can write {\bf (A)} as follows:
$$
  \eqalign{
   {\rm\bf(A)} & =
   {e^{n\phi(\xi)}\over\pi\sqrt nb(\xi)}
   \int_{-\alpha\sqrt nb(\xi)}^{\alpha\sqrt nb(\xi)}
   \exp\left\{-{1\over 2}\theta^2+nR_0(\xi,t_n)\right\}
     (1+O(\xi^2)+O(t_n^2))\cdot\cr
   &\qquad\cdot
   \left\{\cos(n\r1)a_nH_0(z_n)+\cos(n\r1)b_nt_nK_0(z_n)+o(\xi)+o(t_n^2)
   \right\} 
    \d\theta,
  }
$$
where all the ``$O$'' and ``$o$'' are uniform with respect to $\xi\in[0,
n^{-1/4}]$, for every fixed $\theta$. 

{\it Part $I\! I\! I$}

We divide {\bf (A)} in three parts, and
we write $a_n$ and $b_n$ as functions of $\theta$, $\xi$ and $n$:
$$
  \eqalign{
  a_n
%=a\left(\xi,{\theta\over\sqrt nb(\xi)}\right)& =
%  {1\over\sqrt 2}\sqrt{\sqrt{\xi^2+{\theta^2(1-\xi^2)\over 2n}}+\xi}=\cr
	& ={(1-\xi^2)^{1/4}\over\sqrt 2\,n^{1/4}}
	\sqrt{\sqrt{{\theta^2\over 2}+{n\xi^2\over (1-\xi^2)}}+
	{\sqrt n\xi\over\sqrt{1-\xi^2}}}\cr
  b_n
%=b\left(\xi,{\theta\over\sqrt nb(\xi)}\right)& =
%	{{\rm sign}(-\theta)\over\sqrt 2}
%	\sqrt{\sqrt{\xi^2+{\theta^2(1-\xi^2)\over 2n}}-\xi}=\cr
	& =-{\theta(1-\xi^2)^{1/4}\over 2n^{1/4}
	\sqrt{\sqrt{{\theta^2\over 2}+{n\xi^2\over (1-\xi^2)}}+
	{\sqrt n\xi\over\sqrt{1-\xi^2}}}  },
}$$
obtaining
$$
  \eqalign{
   {\rm\bf(A_1)} & =
   {e^{n\phi(\xi)}(1-\xi^2)^{3/4}\over\pi\sqrt 2\, n^{3/4}}
   \int_{-\alpha\sqrt nb(\xi)}^{\alpha\sqrt nb(\xi)}
   \exp\left\{-{1\over 2}\theta^2+nR_0(\xi,t_n)\right\}
     (1+O(\xi^2)+O(t_n^2))\cdot\cr
   &\qquad\qquad\cdot     \cos(n\r1)
       {1\over\sqrt 2} \sqrt{\sqrt{{\theta^2\over 2}+{n\xi^2\over (1-\xi^2)}}+
	{\sqrt n\xi\over\sqrt{1-\xi^2}}}
      H_0(z_n)\d\theta,\cr
   {\rm\bf(A_2)} & =
   {e^{n\phi(\xi)}(1-\xi^2)^{3/4}\over\pi\sqrt 2\, n^{3/4}}
   \int_{-\alpha\sqrt nb(\xi)}^{\alpha\sqrt nb(\xi)}
   \exp\left\{-{1\over 2}\theta^2+nR_0(\xi,t_n)\right\}
     (1+O(\xi^2)+O(t_n^2))\cdot\cr
   &\qquad\qquad\cdot
   \cos(n\r1)
	{\theta\over 2
	\sqrt{ \sqrt{{\theta^2\over 2}+{n\xi^2\over (1-\xi^2)}}+
	{\sqrt n\xi\over\sqrt{1-\xi^2}} }  }
   {\theta\sqrt{1-\xi^2}\over\sqrt{2n}}(-K_0(z_n))  
    \d\theta\cr
   {\rm\bf(A_3)} & =
   {e^{n\phi(\xi)}(1-\xi^2)^{1/2}\over\pi\sqrt 2\, n^{1/2}}
   \int_{-\alpha\sqrt nb(\xi)}^{\alpha\sqrt nb(\xi)}
   \exp\left\{-{1\over 2}\theta^2+nR_0(\xi,t_n)\right\}
     (1+O(\xi^2)+O(t_n^2))\cdot\cr
   &\qquad\cdot
   (o(\xi)+o(t_n^2))\d\theta.
		       }   \autoeqno{a1ea2}
$$

{\it Part $IV$}

We notice that the integral in ${\rm\bf(A_1)}$ seems to be asymptotic
to
$$
I\left(\sqrt n\xi\right)
=\int_\reali e^{-{\theta^2\over 2}}\sqrt{\sqrt{
n\xi^2+{\theta^2
\over 2}}+\sqrt n\xi}\d\theta,
$$
then it will be useful to state some properties of $I(t)$:
$I(t)$ exists and is finite for every $t\ge 0$,
it is continuous and increasing in $[0,+\infty)$,
differentiable in $(0,+\infty)$, $I(t)\ge 2\sqrt{t\pi}$ and if $t$ tends to infinity,
$I(t)\,\sim\,2\sqrt{t\pi}$. 

Now we want to apply Theorem~\lemmaref{dominscon1} to the integrand of
${\rm\bf(A_1)}$ in \eqref{a1ea2}, but a new problem raises: the quantity
$\sqrt n\xi$ may tend to $0$ or to $+\infty$ as well.
Hence we introduce a new function, namely $Q(n,\xi):=\max\left\{1,
\sqrt n\xi\right\}\ge 1$ ($Q(n,\xi)$ tends to infinity
if and only if $\sqrt n\xi\ra +\infty$).
Once again we rewrite ${\rm\bf(A_1)}$:
$$
  \eqalign{
   {e^{n\phi(\xi)}(1-\xi^2)^{3/4}\sqrt{Q(n,\xi)}\over\pi\sqrt 2 n^{3/4}}
 &  \int_{-\alpha\sqrt nb(\xi)}^{\alpha\sqrt nb(\xi)}
    \exp \left\{-{1\over 2}\theta^2+nR_0(\xi,t_n)\right\}\cdot\cr
   & \cdot (1+O(\xi^2)+O(t_n^2))cos(n\r1)
   {1\over\sqrt 2}H_0(z_n)\cdot\cr
   & \cdot\
    \sqrt{\sqrt{{\theta^2\over 2Q^2(n,\xi)}+{n\xi^2\over (1-\xi^2)Q^2(n,\xi)}}
    + {\sqrt n\xi\over\sqrt{1-\xi^2}Q(n,\xi)}}
      \d\theta.
}$$
We must first evaluate the (uniform) asymptotic value of the integrand:
recalling that $|nR_0(\xi,t_n)|\ra 0$ for every fixed $\theta$, uniformly
with respect to $\xi\in[0,n^{-1/4}]$, it is easy to show that
$$
  \eqalign{
  f_n(\theta,\xi) & :=\ident_{[-\alpha\sqrt nb(\xi),
  \alpha\sqrt nb(\xi)]}(\theta)
   \exp \left\{-{1\over 2}\theta^2+nR_0(\xi,t_n)\right\}\cdot\cr
   & \cdot (1+O(\xi^2)+O(t_n^2))cos(n\r1)
   {1\over\sqrt 2}H_0(z_n)\cdot\cr
   & \cdot\
    \sqrt{\sqrt{{\theta^2\over 2Q^2(n,\xi)}+{n\xi^2\over (1-\xi^2)Q^2(n,\xi)}}
    + {\sqrt n\xi\over\sqrt{1-\xi^2}Q(n,\xi)}}\,\sim\cr
    & \sim\,h_n(\theta,\xi):=\exp(-\theta^2\!/2)\cdot\cr
   & \cdot\
    \sqrt{\sqrt{{\theta^2\over 2Q^2(n,\xi)}+{n\xi^2\over (1-\xi^2)Q^2(n,\xi)}}
    + {\sqrt n\xi\over\sqrt{1-\xi^2}Q(n,\xi)}}.
}$$
As for a uniform bound for $f_n$ and $h_n$, they are both uniformly bounded
by
$$
C\,\exp(-\theta^2\!/4)\sqrt{\sqrt{1+{\theta^2\!/2}}+1}.
$$
Finally, $|\int_\reali h_n(\theta,\xi)\d\theta|>c$ for some $c>0$, for all
$n$ and $\xi\in[0,n^{-1/4}]$, since $h_n(\theta,\xi)\ge e^{-\theta^2/2}\cdot
\min\{c,2^{-1/4}\sqrt{|\theta|}\}$.

Hence, by Theorem~\lemmaref{dominscon1},
$$
%\eqalign{
{\rm\bf(A_1)}\,{\buildrel n\over\sim}\,
%   {e^{n\phi(\xi)}\over\pi\sqrt 2\, n^{3/4}} &
%   \sqrt{Q(n,\xi)}  \int_\reali
%   \exp (-\theta^2\!/2)\cdot\cr
%   & \cdot\
%    \sqrt{\sqrt{{\theta^2\over 2Q^2(n,\xi)}+{n\xi^2\over (1-\xi^2)Q^2(n,\xi)}}
%    + {\sqrt n\xi\over\sqrt{1-\xi^2}Q(n,\xi)}}
%      \d\theta=\cr &
	 = {e^{n\phi(\xi)}\over\pi\sqrt 2 n^{3/4}}
	   I\left(\sqrt n\xi\right).
%}
$$

{\it Part $V$}

We are now able to show that {\bf (B)} is asymptotically negligible if
compared to ${\rm\bf (A_1)}$, using Lemma~\lemmaref{int-2-y}:
$$
  \left|{{\rm\bf (B)}}\over {{\rm\bf (A_1)}}\right|\le{C e^{n\phi(\xi)}
	\lambda^n\over e^{n\phi(\xi)}n^{-3/4}I(0)}\le C\lambda^nn^{3/4}
$$
and the last term tends to $0$ uniformly with respect to $\xi\in[0,n^{-1/4}]$.

We estimate ${\rm\bf (A_2)}$ and ${\rm\bf (A_3)}$ to show their negligibility
too:
$$
\eqalign{
   |{\rm\bf(A_2)}| & =
   {e^{n\phi(\xi)}(1-\xi^2)^{5/4}\over 4\pi n^{5/4}}
   \int_{-\alpha\sqrt nb(\xi)}^{\alpha\sqrt nb(\xi)}
   \exp\left\{-{1\over 2}\theta^2+nR_0(\xi,t_n)\right\}
     (1+O(\xi^2)+O(t_n^2))\cdot\cr
   &\qquad\qquad\cdot
   \cos(n\r1)
	{ \theta^2\over 
	\sqrt{ \sqrt{{\theta^2\over 2}+{n\xi^2\over (1-\xi^2)}}+
	{\sqrt n\xi\over\sqrt{1-\xi^2}} }  }
	(-K_0(z_n))      \d\theta\cr
    & \le C{e^{n\phi(\xi)}\over n^{5/4}}
   \int_\reali
   \exp(-\theta^2\!/4)
	{ \theta^2\over |\theta|^{1/2} }    \d\theta\cr
}$$
since $\sqrt{\sqrt{t^2+\theta^2\!/2}+t}\ge 2^{-1/4} |\theta|^{1/2}$ for all
$t\ge 0$.
Then
$$
  \left|{{\rm\bf (A_2)}}\over {{\rm\bf (A_1)}}\right|
  \le C { e^{n\phi(\xi)}n^{-5/4}\over e^{n\phi(\xi)}n^{-3/4}}\le Cn^{-1/2}
$$
which tends uniformly to $0$ with respect to
$\xi\in[0,n^{-1/4}]$.

As for ${\rm\bf (A_3)}$, we observe that every $o(\xi)$
besides those in \eqref{opiccoli2}
is also equal to $\xi^{1/8}o(\xi)$, while for the $o(\xi)$ in \eqref{opiccoli2}
we have that 
$$\eqalign{
& {\xi\cdot O(\sqrt{|t_n|})\over\xi}\le C\sqrt{|t_n|}\le 
	C{\theta^{1/2}\over n^{1/4}},\cr
& {O(\xi)\cdot O(|t_n|^{3/2})\over\xi}\le C|t_n|^{3/2}\le 
	C{\theta^{3/2}\over n^{3/4}}.}
$$
Moreover, $|O(\xi^2)|$ and $|O(t_n^2)|$ are smaller than some constant $C>0$ (recall
that $|t_n|\le\alpha$ and $\xi\le n^{-1/4}$), and $|o(t_n^2)|\le C\,{\theta\over n}$,
hence we can majorize $\rm\bf |(A_3)|$:
$$
%\eqalign{  
|{\rm\bf (A_3)}| 
%&\le C\,e^{n\phi(\xi)}\xi n^{-1/2}
%	\,\int_\reali e^{-\theta^2\!/4}\left(\xi^{1/8}+{\theta^{1/2}\over n^{1/4}}
%	+{\theta^{3/2}\over n^{3/4}} \right)\d\theta+\cr
%	& +C\,e^{n\phi(\xi)} n^{-1/2}
%	\,\int_\reali e^{-\theta^2\!/4}{\theta\over n}\d\theta\le\cr &
	 \le C\,e^{n\phi(\xi)}n^{-3/4}
	\left(\xi^{1/8}+n^{-1/4} \right)
%}
$$
whence
$$
  \left|{{\rm\bf (A_3)}\over{\rm\bf (A_1)}}\right|\le C\,
	{e^{n\phi(\xi)}n^{-3/4}\left(\xi^{1/8}+n^{-1/4} \right)
	\over e^{n\phi(\xi)}n^{-3/4}
	I(0) }\le C\,\left(\xi^{1/8}+n^{-1/4}\right)
$$
and the  last term tends uniformly to $0$ with respect to
$\xi\in[0,n^{-1/4}]$.

The statement for $\xi\in[0,n^{-1/2-\eps}]$ is proved in the same way,
once we note that
$$
I(0)=\sqrt2\Gamma\left({3\over 4}\right),\hbox{ and }
\Gamma\left({3\over 4}\right)={\sqrt2\,\pi\over\Gamma\left({1\over 4}
\right)}
$$
(use $\Gamma\left({1\over 2}\right)=\sqrt\pi$ and duplication formula
$\Gamma(2z)=(2\pi)^{-1/2}2^{2z-1/2}\Gamma(z)\,\Gamma\left(z+{1\over 2}
\right)$).
\QED

%% file: x.tex
%%%%%%%%%%%%%%%%%%%%%%%COMBX.TEX%%%%%%%%%%%%%%%%%%%%%%%%%%%%%%%%%%%%%%%%%
%%%%%%%%%%%%%%%%%%%%%%%%%%%%%%%%%%%%%%%%%%%%%%%%%%%%%%%%%%%%%%%%%%%%%%%%%%%%%
%% INIZIANO I CONTI LUNGO x
%%%%%%%%%%%%%%%%%%%%%%%%%%%%%%%%%%%%%%%%%%%%%%%%%%%%%%%%%%%%%%%%%%%%%%%%%%%%%
\autosez{x}{Estimates along the $x$-axis: the case $\xi\in[a,1-c]$}

In this section we want to give an asymptotic estimate 
for the transition probabilities
$p^{(n)}((k,0),(0,0))$, which we will rename $p^{(n)}(k,0)$.
We first note that if $n$ and $k$ do not have the same parity, then
$p^{(n)}(k,0)=0$, but we cannot repeat the trick we used on the $y$-axis
to derive $p^{(2n+1)}(2k+1,0)$ from $p^{(2n)}(2k,0)$.
We estimate only the second type of transition probabilities (the
first ones can be derived in a similar way --- see Paragraph~\sref{fr}).

The basic idea underlying our proofs here is essentially the same as for
the case of the $y$-axis, nevertheless much more technical difficulties will
arise and the techniques we will employ appear more involved.
We want to point out that  here we will use the same symbols (such as
$F$, $\Psi_\xi$, $\phi$) for functions and values which are not the same
but play the same role as the analogous for the $y$-axis.

By Cauchy's integral formula,
$$
  p^{(2n)}(2k,0)={1\over 2\pi i}\int_{\gamma^\prime}{\widetilde G(z)\cdot
	\widetilde F_2(z)^{2k}\over z^{2n+1}}\d z
$$
where $\gamma^\prime$ is a positively
oriented, simple closed curve in $\complessi$ which has $0$ in its interior
and $1$ in its exterior.

The substitution $u=z^2$ appears to be useful and when $\gamma^\prime
\,:\,z=z(t)$
describes one circuit about $0$ then $\gamma\,:\,u=u(t)$ describes 
two times the corresponding circuit, then
$$
%\eqalign{
  p^{(2n)}(2k,0) 
%& ={2\over 2\pi i}
%  \int_{\gamma}{G(u)F(u)^{2k}\over u^n\sqrt u}\cdot
%	{\d u\over 2\sqrt u}\cr
%	& 
={1\over 2\pi i}\int_\gamma{G(z)F(z)^{2k}\over z^{n+1}}\d z
%}
$$
where $F(z)=F_2(z)$.

We first stress the ``exponential part'' of the integrand, rewriting
$$
%\eqalign{
 {F(z)^{2k}\over z^n} %& =\exp\{2k\log(F(z))-n\log(z)\}=\cr
%	& =\exp\{2k\log(1+\sqrt{1-z}-\sqrt 2\sqrt{1-z+\sqrt{1-z}})
%	-(k+n)\log(z)\}=\cr &
	 =\exp\{n[2\xi\log(1+\sqrt{1-z}-\sqrt 2\sqrt{1-z+\sqrt{1-z}})
	-(\xi+1)\log(z)]\},
%}
$$
where $\xi:={k\over n}$.
We now introduce the auxiliary function 
$$\Psi_\xi(z):=2\xi\log(1+\sqrt{1-z}-\sqrt 2\sqrt{1-z+\sqrt{1-z}})-
(\xi+1)\log(z),$$ 
and rewrite the expression for the transition probabilities:
$$
  p^{(2n)}(2k,0) ={1\over 2\pi i}
  \int_\gamma{G(z)\over z}\exp\{n\Psi_\xi(z)\}\d z.
\autoeqno{x1}
$$
\lemma{phixi1}{Let $\phi(\xi):=\min\{\Psi_\xi(z)\,:\,0\le z\le1\}$ and let 
$z_o(\xi)$ be the (unique) point where this minimum is attained. Then $z_o(\xi)=
1-u_o(\xi)^2$, where
$u_o(\xi)=\xi^{2/3}[54+6\sqrt{81-6\xi^2}]^{1/3}/6+
\xi^{4/3}[54+6\sqrt{81-6\xi^2}]^{-1/3}$.
%\lemma{derPsix}{
The second and the third order derivatives of $\Psi_\xi(z)$
take the following values in $z_o(\xi)$:
$$
\eqalign{
  \Psi^{\prime\prime}_\xi(z_o(\xi)) & ={1+2z_o(\xi)-\sqrt{1-z_o(\xi)}
	\over 4z_o(\xi)^2(1-z_o(\xi))}=\cr
	& ={2u_o(\xi)+3\over 4u_o(\xi)^2(1-u_o(\xi))(1+u_o(\xi))^2};\cr
  \Psi^{\prime\prime\prime}_\xi(z_o(\xi)) & =
	{-16\uo^2-24\uo+19\over 16\uo^4(1-\uo)^2(1+\uo)};
}$$
(the expression for the third order derivative appears much more
complex as a function of $\zo$).}
The following theorem is the analogous of Theorem~\lemmaref{y-xinonzero}
for the $x$-axis (and the proof is just the same).

\theorem{x-xinonzero}{Let $a,c$ be positive numbers such that $a<1-c$ and let
$\xi\in[a,1-c]$. Then
$$
  p^{(2n)}(2k,0)\,{\buildrel n\over\sim}\,\sqrt{2\over\pi}
	{e^{n\phi(\xi)}\sqrt{1-\zo}G(z_o(\xi))
	\over\sqrt{1+2\zo-\sqrt{1-\zo}}}\,
	n^{-1/2}
$$
uniformly with respect to $\xi\in [a,1-c]$.}

%%%%%%%%%%%%%%%%%%%%%%%%%%%%%%%%%%%%%%%%%%%%%%%%%%%%%%%%%%%%%%%%%%%%%%%%%%%%%%%
%%% FINE DELLA PARTE DI x E \xi\in [a,1-c]
%%%%%%%%%%%%%%%%%%%%%%%%%%%%%%%%%%%%%%%%%%%%%%%%%%%%%%%%%%%%%%%%%%%%%%%%%%%%%%
%%%%%%%%%%%%%%%%%%%%%%%%%%%%XPROX0.TEX%%%%%%%%%%%%%%%%%%%%%%%%%%%%%%%%%%%%%%%%%%%%%
\autosez{x-xi<a}{Estimates along the $x$-axis: the case $\xi\in [0,a]$}

In order to choose the curve of integration, we perform a change of variable
 which is the inverse mapping of $z=f(v)$. The curve of integration will be
the union of two pieces: one
will be the image of a suitable curve in the $v$-plane, and the other
will be a part of a circle centered in the origin such that the union of the
two pieces forms a connected circuit about the origin, that is, as in the case 
of the $y$-axis,
$$
p^{(2n)}(2k,0)={\rm\bf(A)}+{\rm\bf(B)}={1\over 2\pi i}
  \int_{\gamma_1}{G(z)F(z)^{2k}\over z^{n+1}}\d z
 +{1\over 2\pi i}\int_{\gamma_2}{G(z)F(z)^{2k}\over z^{n+1}}\d z,
 \autoeqno{xab}
$$
where $\gamma_1$ corresponds to a particular curve $z(\xi,t)$, $|t|\le\alpha$,
and $\gamma_2$ corresponds to $|z(\xi,\alpha)|e^{is}$, $s\in[\arg(z(\xi,
\alpha)),2\pi-\arg(z(\xi,\alpha))]$.

We choose $z=1-v^4$, hence $v:= \root 4\of{1-z}$ is the inverse mapping
if, by our choice of the determination of the
root, $v$ has argument in $[-\pi/4,\pi/4)$.

 The desired curve in the $v$-plane is a
segment whose parametrization is 
$$
  v(\xi,t)=\cases{
  v(\xi)+e^{i\beta}t  & if $t\in[0,\alpha]$ \cr
  v(\xi)-e^{-i\beta}t & if $t\in[-\alpha,0)$,\cr 
}$$  
where $v(\xi):=\sqrt\uo$, and $\alpha$ and $\beta$
will be chosen in the sequel.
The corresponding curve in the $z$-plane is $z(\xi,t):=1-v(\xi,t)^4$.

In order to obtain the Taylor series of the function $t\mapsto
\ovl\Psi_\xi(t):=\Psi_\xi(z(\xi,t))$ we introduce the auxiliary
function 
$$
  \widetilde\Psi_\xi(v):=\Psi_\xi(1-v^4)=
  2\xi\log(1+\sqrt{v^4}-\sqrt 2\sqrt{v^4+\sqrt{v^4}})
  -(\xi+1)\log(1-v^4).
$$
Particular attention must be paid in computing the roots ($\sqrt{v^4}$
is not necessarily equal to $v^2$ nor $\sqrt{a\cdot b}$ is always equal to
$\sqrt a\cdot\sqrt b$).
For our purposes we only need the explicit expressions of 
$\widetilde\Psi_\xi(v)$ for $\arg(v)\in [-\pi/4,0]$ (where $\sqrt{v^4}=v^2$,
$\sqrt{v^4+v^2}=v\sqrt{1+v^2}$) and for $\arg(v)\in [\pi/4,\pi/2)$ 
(where $\sqrt{v^4}=-v^2$,
$\sqrt{v^4-v^2}=-iv\sqrt{1-v^2}$).

Hence
$$\eqalign{
     \widetilde\Psi_{\xi,1} (v)& :=2\xi\log(1+v^2-\sqrt 2\,v\sqrt{1+v^2})-
		(\xi+1)\log(1-v^4), \cr
     \widetilde\Psi_{\xi,2} (v)& :=2\xi\log(1-v^2+\sqrt 2\,iv\sqrt{1-v^2})-
	       (\xi+1)\log(1-v^4) ,
}$$
where $\widetilde\Psi_{\xi,1}(v)=\widetilde\Psi_{\xi}(v)$ if 
$\arg(v)\in [-\pi/4,0]$ and $\widetilde\Psi_{\xi,2}(v)=
\widetilde\Psi_{\xi}(v)$ if $\arg(v)\in [\pi/4,\pi/2)$.

\lemma{psiv}{The functions $\widetilde\Psi_{\xi,1}(v)$ and
$\widetilde\Psi_{\xi,2}(v)$ have a Taylor series expansion centered in $0$,
with convergence radius equal to $1$:
$$
  \widetilde\Psi_{\xi,1}(v) =\sum_{n=1}^\infty d_n(\xi)\, v^n,\qquad\qquad
  \widetilde\Psi_{\xi,2}(v) =\sum_{n=1}^\infty d_n^\prime(\xi)\,v^n
$$
where $d_{4n}(\xi)=d^\prime_{4n}(\xi)=1/n$ for every $n\ge 1$;
$d_{2n+1}(\xi)=-2\sqrt 2\,\xi
{b_n\over 2n+1}$, $d^\prime_{2n+1}(\xi)=2\sqrt 2\,i\xi\,{(-1)^nb_n\over 2n+1}$
for every $n\ge0$, $b_n:=\sum_{i=0}^{[n/2]}\pmatrix{1/2\cr n-2i}$ and
all other coefficients are equal to zero.
}

\lemma{psitx}{The function $\ovl\Psi_\xi(t)=\widetilde\Psi_\xi(v(\xi,t))$ has
a Taylor series expansion centered in $0$, with positive radius of 
convergence $\rho$ not depending on $\xi$:
$$
  \ovl\Psi_\xi(t)=\sum_{m=0}^\infty e^{im\beta}g_m(\xi)\,t^m,
  \autoeqno{psitx}
$$
where 
$$g_m(\xi)=\sum_{h=m}^\infty \pmatrix{h\cr m}d_h(\xi)v(\xi)^{h-m}$$
if $\arg(v(\xi,t))\in[-\pi/4,0]$ and
$$g_m(\xi)=\sum_{h=m}^\infty \pmatrix{h\cr m}d^\prime_h(\xi)v(\xi)^{h-m}$$
if $\arg(v(\xi,t))\in[\pi/4,\pi/2]$. 
(Note that the series defining $g_m(\xi)$
both converge for every $\xi\in[0,a]$, provided that $a$ is sufficiently
small).}

To avoid misunderstandings, we define $\ovl\Psi_{\xi,j}(t)$, $j=1,2$, 
as follows:
$$
 \ovl\Psi_\xi(t)=\cases{
  \ovl\Psi_{\xi,1}(t)  & if $\arg(v(\xi,t))\in[-\pi/4,0]$ \cr
  \ovl\Psi_{\xi,2}(t)  & if $\arg(v(\xi,t))\in[\pi/4,\pi/2]$.
}$$

\lemma{psitx-sv}{The following equalities hold:
$$
\eqalign{
\ovl\Psi_{\xi,1}(t) & =\phi(\xi)+e^{2i\beta}a_2(\xi)t^2+e^{3i\beta}a_3(\xi)t^3
	+e^{4i\beta}a_4(\xi)t^4+R(\xi,t)\cr
\ovl\Psi_{\xi,2}(t) & =\phi(\xi)+e^{2i\beta}a^\prime_2(\xi)t^2
	+e^{3i\beta}a_3^\prime(\xi)t^3
	+e^{4i\beta}a^\prime_4(\xi)t^4+R(\xi,t)
}$$
where the remainder term  (which is not necessarily equal in the two 
expressions) is \break $R(\xi,t)=O(t^5)$,
uniformly with respect to $\xi\in[0,a]$, $|t|\le
\alpha$. Moreover, if $\xi\ra 0$
$$
  \matrix{
	a_2(\xi)\,\sim\,3\cdot 2^{2/3}\xi^{2/3};\hfill
	&  a_2^\prime(\xi)\,\sim\,3\cdot 2^{2/3}\xi^{2/3}-i2^{1/3}\xi^{4/3};
	\hfill\cr
	a_3(\xi)\,\sim\,2^{11/6}\xi^{1/3};\hfill
	  & a_3^\prime(\xi)\,\sim\,2^{11/6}\xi^{1/3}-i\sqrt 2/3\,\xi;
	  \hfill\cr
	a_4(\xi)\,\sim\,1;\hfill
	  & a_4^\prime(\xi)\,\sim\,1+i\cdot7\cdot 2^{-3/2}\xi^{4/3}.\hfill\cr
}$$}

%Now we rewrite equation~\eqref{xab}:
%$$
%  \eqalign{
% p^{(2n)}(2k,0) & ={\rm\bf(A)}+{\rm\bf(B)}=\cr
%  & ={1\over \pi} \im\,
%	\int_0^{\alpha}{G(\zt)\over\zt}\exp\{n\ovl\Psi_\xi(t)\}
%	[-4e^{i\beta}(v(\xi)+e^{i\beta}t)^3]\d t+\cr
%	& +{1\over 2\pi i}\int_{|s|\ge\arg(z(\xi,\alpha))}
%	{G(|z(\xi,\alpha)|e^{is})\,
%	(F(|z(\xi,\alpha)|e^{is})^2)^k\over
 %  |z(\xi,\alpha)|^ne^{isn}}\d s.
%}$$
Now we have to choose the curve of integration (that is, $\beta$) in order
to obtain that:
\item{(a)} $\ovl\Psi_\xi(t):=\widetilde\Psi_\xi(v(\xi,t))$ has a Taylor expansion
	in a neighbourhood of $t=0$;
\item{(b)} the curve $z(\xi,t)=1-v(\xi,t)^4$ has some ``good properties''
	(for instance there exists $\alpha>0$ such that $|z(\xi,\alpha)|\ge 
	1+\eps_0$ for some $\eps_0>0$ and for every $\xi$ in the considered 
	range);
\item{(c)}we can bound $\exp\{n\ovl\Psi_\xi(t)\}$ with an integrable function.

We observe that, not considering the integral {\bf(B)} in 
equation~\eqref{xab}, we have 
to integrate on the curve $z_1(\xi,t):=1-(v(\xi)+
e^{i\beta}t)^4$ for $t\in [0,\alpha]$ and on $z_2(\xi,t):=1-(v(\xi)+
e^{i\beta}t)^4$ for $t\in [-\alpha,0)$, but $\ovl{z_1(\xi,t)}=
z_2(\xi,-t)$, that is the second curve is simply the conjugate of the
first one. Then we denote by $\zt$ the curve
$z_1(\xi,t)$ and the sum of the two corresponding
integrals will be 
$$
%  \eqalign{
  {\rm\bf(A)} % & ={1\over 2\pi i} \bigg\{
%	\int_0^{\alpha}{G(\zt)\over\zt}\exp\{n\ovl\Psi_\xi(t)\}
%	 [-4e^{i\beta}(v(\xi)+e^{i\beta}t)^3]\d t+\cr
%	& -\int_0^{\alpha}{G(\ovl\zt)\over\ovl\zt}
%	\exp\{n\ovl\Psi_\xi(-t)\}
%	[-4e^{-i\beta}(v(\xi)+e^{-i\beta}t)^3]\d t\bigg\}=\cr
  ={1\over \pi}  \im\,
	\int_0^{\alpha}{G(\zt)\over\zt}\exp\{n\ovl\Psi_\xi(t)\}
	[-4e^{i\beta}(v(\xi)+e^{i\beta}t)^3]\d t, %}
$$
since $G(\ovl z)=\ovl{G(z)}$ and $\exp\{n\ovl\Psi_\xi(-t)\}=
\ovl{\exp\{n\ovl\Psi_\xi(t)\}}$ (remember that $\exp\{n\ovl\Psi_\xi(t)\}=
F(\zt)^{2k}/\zt^n$, and $F(\ovl z)^2=\ovl{F(z)^2}$).

This allows us to consider only the case $t\ge0$ (that is, the part of the
curve that lies in the first quadrant).
%%%%%%%%%%%%%%%%%%%%%%%%%%%%%%%%%%%%%
%% FORSE E` IL CASO DI SOTTOLINEARE QUI COME NEL CASO y
%% CHE LA PSI(t) E`SVILUPPABILE IN SERIE DI TAYLOR IN QUANTO
%% COMPOSIZIONE DI FUNZIONI ANALITICHE (UNA SOMMA DI ANALITICHE 
%% COMPOSTA CON UN POLINOMIO
%%%%%%%%%%%%%%%%%%%%%%%%%%%%%%%%%%%%

As we would like to give an upper bound for the exponential part in
{\rm\bf(A)}, we restrict the range of $\beta$. In fact we require
that the real parts of the expansion in Lemma~\lemmaref{psitx-sv}
of at least $\ovl\Psi_{\xi,1}(t)$ have non positive coefficients, that is,
that $\cos(2\beta)$, $\cos(3\beta)$ and $\cos(4\beta)$ are all non
positive.
It is easy to show that this corresponds to $\beta\in\Delta:=
[\pi/4,3\pi/8]\cup
[-3\pi/8,-\pi/4]$, that is, $e^{i\beta}$ lies in two conjugate sectors.

{\bf Remark}
%\lemma{curv}{
Among the curves 
of the family ${\cal F}:=\{\zt\,:\,t\ge 0,\ \beta\in\Delta\}$, the ones
with $\beta=\pm\pi/4$ have exactly one intersection with the real axis:
 $z(\xi,0)=1-v(\xi)^4$; while the others have exactly two
intersections with the same axis: one for $t=0$ and the other for
$t={v(\xi)\over \sin\beta-\cos\beta}$ if $\beta\in(\pi/4,3\pi/8]$ or
for $t=-{v(\xi)\over \sin\beta+\cos\beta}$ if $\beta\in[-3\pi/8,-\pi/4]$.
%}
\medskip

We present the plot of some of the curves we will use.

%%%%CURVA CON BETA=\PI/4
\beginpicture
\setcoordinatesystem units <5mm,5mm> point at 0 0
\setlinear
\put{$\,$} at -2 4.5
\put{$\,$} at -2 -6.5
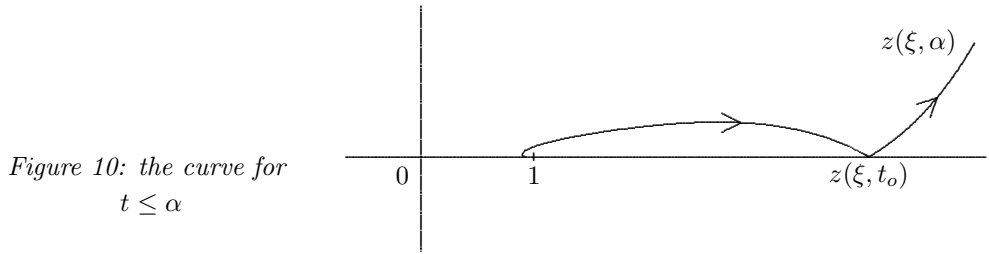
\put{\it Figure~6: the curve for} [lt] at 0 0
\put{$\beta={\pi\over 4}$} [lt] at 3 -1
\plot  9 0  26 0 /
\plot  11 -6    11 4 /
\put{0}  at 10.5 -.5
\put{1}  at 13 -.5
%freccetta
\plot 18.85603 -3.861168  18.40603 -3.461168 /
\plot 18.85603 -3.861168  18.30603 -3.861168 /
\setquadratic
\plot   12.7822 0               12.729 -.0239718        12.67841 -.0562826
	12.63332 -.0978462      12.59692 -.1495768      12.57269 -.2123882
	12.56443 -.2871944      12.57621 -.3749093      12.61243 -.4764468
	12.67777 -.592721       12.77722 -.7246459      12.91606 -.8731352
	13.09989 -1.039103      13.3346 -1.223463       13.62636 -1.42713
	13.98166 -1.651017      14.40731 -1.896039      14.91037 -2.163109
	15.49824 -2.45314       16.1786 -2.767049       16.95944 -3.105747
	17.84906 -3.470149      18.85603 -3.861168      19.98924 -4.27972
	21.25788 -4.726718      22.67145 -5.203076      24.23971 -5.709707 /
%%%%%%% stanghetta per 1
\setlinear
\linethickness=.3pt
\plot 13 -.1  13 .1 /
\endpicture

The plot shows that the curve for $\beta={\pi\over 4}$ turns clockwise (with
respect to the origin),
then we will prefer its conjugate, that is the curve for 
$\beta=-{\pi\over 4}$.

%%%%%CURVA CON BETA=-\PI/4
\beginpicture
\setcoordinatesystem units <5mm,5mm> point at 0 0
\setlinear
\put{$\,$} at -2 6.5
\put{\it Figure~7: the curve for} [lt] at 0 0
\put{$\beta=-{\pi\over 4}$} [lt] at 3 -1
\plot  9 0  26 0 /
\plot  11 -3.5    11 6 /
\put{0} at 10.5 -.5
\put{1}  at 13 -.5
%freccetta
\plot 18.85603 3.861168  18.40603 3.461168 /
\plot 18.85603 3.861168  18.30603 3.861168 /

\setquadratic
\plot   12.7822 0               12.729 .0239718         12.67841 .0562826
	12.63332 .0978462       12.59692 .1495768       12.57269 .2123882
	12.56443 .2871944       12.57621 .3749093       12.61243 .4764468
	12.67777 .592721        12.77722 .7246459       12.91606 .8731352
	13.09989 1.039103       13.3346 1.223463        13.62636 1.42713
	13.98166 1.651017       14.40731 1.896039       14.91037 2.163109
	15.49824 2.45314        16.1786 2.767049        16.95944 3.105747
	17.84906 3.470149       18.85603 3.861168       19.98924 4.27972
	21.25788 4.726718       22.67145 5.203076       24.23971 5.709707 /

%%%%%%% stanghetta per 1
\setlinear
\linethickness=.3pt
\plot 13 -.1  13 .1 /

\put{$\,$} at 0 -4
\endpicture

Here is the curve for $\beta={\pi\over 3}$, where we used
a logarithmic scale on the horizontal axis in order to show the
``pathological'' behaviour of the curves of the family with 
$\beta\neq \pm\pi/4$.

\beginpicture
\setcoordinatesystem units <5mm,5mm> point at 0 0
\setlinear
%%%%%CURVA CON BETA=\PI/3
\put{$\,$} at -2 6.5
\put{$\,$} at -2 -5.5
\put{\it Figure~8: the curve for} [lt] at 0 0
\put{$\beta={\pi\over 3}$} [lt] at 3 -1
\plot  9 0  26 0 /
\plot  11 -5    11 6 /
\put{0}  at 10.5 .5
%\put{$\times$} at 14 0
\put{1}  at 14 .5
\setquadratic
\plot   13.75768 0              13.72046 -.0222696      13.69853 -.04963793
	13.69639 -.08198807     13.71849 -.1191252      13.76898 -.1607764
	13.85151 -.2065911      13.96894 -.2561405      14.12311 -.308918
	14.31473 -.3643391      14.54332 -.4217412      14.80728 -.4803839
	15.10404 -.5394488      15.43028 -.5980395      15.78222 -.6551818
	16.15586 -.7098234      16.5472 -.7608343       16.95241 -.8070062
	17.36796 -.8470531      17.79067 -.8796111      18.2178 -.9032381
	18.64697 -.9164144      19.07621 -.9175422      19.50387 -.9049456
	19.92863 -.876871       20.34944 -.8314868      20.76545 -.7668833
	21.17604 -.6810731      21.58073 -.5719907      21.97918 -.4374927
	22.37115 -.2753577      22.75651 -.0832864      23.13516 .1410984
	23.5071 .4002518        23.87236 .696707        24.23099 1.033075
	24.58308 1.412045       24.92874 1.836384       25.2681 2.308936
	25.60129 2.832625       25.92846 3.410451 /
%%%%%%% freccette
\setlinear
% frecce su una retta a 45 gradi: se (x0,y0)=punta
%       e (x1,y1)=coda bassa, (x2,y2)=coda alta
%       allora  x1=x0-(.25)*1/sqrt(2)  ;y1=y0-(.75)*1/sqrt(2)
%               x2=x0-(.75)*1/sqrt(2)  ;y2=y0-(.25)*1/sqrt(2)
\plot   24.54148 1.041186    24.71826 1.571517 /
\plot   24.18793 1.394740    24.71826 1.571517 /

%%%%%%% stanghetta per 1
\setlinear
\linethickness=.3pt
\plot 14 -.1  14 .1 /
%%%%%%%%% intersezioni con i reali
\setdots <2pt>
\plot 13.75768 -.2 13.75768 -3 /
\put{$1-v(\xi)^4$} at 13.75768 -3.5

\plot 22.90456 -.2  22.90456 -3 /
\put{$1+4v(\xi)^4{9\over (\sqrt 3-1)^4}$} at 22.90456 -3.5

\endpicture

Now we have to choose a proper change of variable to perform
in the integral {\rm\bf (A)}. The choice is between three substitutions,
each stressing a different piece of the exponential part of the integrand:
\item{(a)} $na_2(\xi)t^2=\theta^2$
\item{(b)} $na_3(\xi)t^3=\theta^3$
\item{(c)} $na_4(\xi)t^4=\theta^4$

\lemma{subs}{If we want to have an upper bound for $n\ovl\Psi_\xi(t(\theta))$,
then substitution (a) will do in the case $\xi\ge n^{-3/4}$, substitution
(c) in the case $\xi\le n^{-3/4}$.}

\proof
With substitution (a) and $t=t_n=t_n(\theta):={\theta\over\sqrt{na_2(\xi)}}$,
$\exp\{n\ovl\Psi_\xi(t)\}$ becomes either
$$
  e^{n\phi(\xi)}\exp\left\{e^{2i\beta}\theta^2+e^{3i\beta}{a_3(\xi)\over
  \sqrt n\,(a_2(\xi))^{3/2}}\,\theta^3+ e^{4i\beta}{a_4(\xi)\over 
  n(a_2(\xi))^{2}}\,\theta^4+nR(\xi,t_n)    \right\}
$$
or
$$
  e^{n\phi(\xi)}\exp\left\{e^{2i\beta}{a_2^\prime(\xi)\over a_2(\xi)}\,
  \theta^2+e^{3i\beta}{a_3^\prime(\xi)\over
  \sqrt n\,(a_2(\xi))^{3/2}}\,\theta^3+ e^{4i\beta}{a_4^\prime(\xi)\over 
  n(a_2(\xi))^{2}}\,\theta^4+nR(\xi,t_n)    \right\}
$$
(depending on where the point $v(\xi,t)$ lies in the complex plane).
We choose $a$ sufficiently small that for every $\xi\in[0,a]$,
$a_i\in[c\cdot h_i(\xi),C\cdot h_i(\xi)]$ $i=2,3,4$, where $h_i(\xi)$ is the asymptotic 
value of $a_i(\xi)$ as $\xi$ tends to zero, $0<c<C$, and the same holds
for the modulus of the real and imaginary parts of $a^\prime_i(\xi)$.

Then
$$
  \matrix{
   {a_3(\xi)\over\sqrt n\,(a_2(\xi))^{3/2}}\le {C\over n^{1/2}\xi^{2/3}};
	\hfill
   & {a_4(\xi)\over n(a_2(\xi))^{2}}\le{C\over n\xi^{4/3}};\hfill\cr
   \left|{\re a_3^\prime(\xi)\over\sqrt n\,(a_2(\xi))^{3/2}}\right|
	\le {C\over n^{1/2}\xi^{2/3}};\hfill
   & \left|{\re a_4^\prime(\xi)\over n(a_2(\xi))^{2}}\right|
	\le{C\over n\xi^{4/3}};\hfill\cr
}$$
which are all surely bounded if $\xi\ge n^{-3/4}$, while
$$
   \matrix{
   \left|{\im a_3^\prime(\xi)\over\sqrt n\,(a_2(\xi))^{3/2}}\right|
	\le {C\over n^{1/2}};\hfill
   & \left|{\im a_4^\prime(\xi)\over n(a_2(\xi))^{2}}\right|
	\le{C\over n}\hfill\cr
}$$
both tend to zero uniformly with respect to $\xi$.

With substitution (b), $t=t_n:={\theta\over\root 3\of{na_3(\xi)}}$, 
and the coefficient of $\theta^2$ is bounded if $\xi\le n^{-3/4}$,
while the coefficient of $\theta^4$ is bounded if $\xi\ge n^{-3/4}$
(the coefficient of $\theta^3$ is $e^{3i\beta}$).
This makes substitution (b) not a suitable one.

Finally, with substitution (c), 
$t=t_n:={\theta\over\root 4\of{na_4(\xi)}}$, and
the coefficients of $\theta^2$ and $\theta^3$ are bounded for 
$\xi\le n^{-3/4}$, while the coefficient of $\theta^4$ is
$e^{4i\beta}$.
\QED

\autosez{xi<n3/4}{Estimates along the $x$-axis: the case 
$\xi\in[0, n^{-3/4}]$}

We choose $\beta=-{\pi\over 4}$: the curve of integration will be similar 
to that in Figure 8, even if the parabolic piece is here substituted by the
union of the arcs of the curves in Figures 12 and 13 corresponding to
$t\in[0,\alpha]$.

We choose $a$ (eventually smaller than the preceding choices)
depending on $\alpha$ such that for all $\xi\le a$
we have $|z(\xi,\alpha)|\ge 1+\eps_o$ for some fixed $\eps_o>0$
(where $z(\xi,t)=1-(v(\xi)+e^{-\pi/4}t)^4$). 
Thanks
to this choice, the circular part of the curve of integration  will be
far from the singularity $z=1$.
This choice is possible since the mapping $\xi\mapsto z(\xi,\alpha)$ is
continuous and if $\xi\ra 0$ then $z(\xi,\alpha)\ra 1+\alpha^4$.
Then for all $\eps>0$ there exists a right neighbourhood $\cal U$ of $0$
such that whenever $\xi\in\cal U$ we have $|z(\xi,\alpha)|\ge 1+\alpha^2
-\eps$ and the last quantity is greater than $1+\eps_o$ for $\eps$
sufficiently small (for instance $\eps_o=\alpha^4/2$).

Negligibility of {\bf (B)} is proved in the same way as in 
Lemma~\lemmaref{int-2-y}.
\lemma{int-2-x1}{
There exists a sufficiently small $a$ such that
$$
  \left|{1\over 2\pi i}\int_{|s|\ge\arg(z(\xi,\alpha))}
	{G(|z(\xi,\alpha)|e^{is})\,
	(F(|z(\xi,\alpha)|e^{is})^2)^k\over
   |z(\xi,\alpha)|^ne^{isn}}\d s
  \right|
  \le C\,e^{n\phi(\xi)}\,\lambda^n
$$
for some $C>0$, $\lambda<1$ and for all $\xi\in [0,a]$.}

\theorem{ximin}
{If $\xi\in[0, n^{-3/4}]$ then
$$
  p^{(2n)}(2k,0)\asym{n}{4\,e^{n\phi(\xi)}\over
  \pi }I(n^{1/4}\xi^{1/3})\,n^{-3/4},
$$
where 
$$
\eqalign{
I(t)  := \int_0^{+\infty} &
 \exp(-2^{4/3}t\theta^3-\theta^4)\cdot[\cos(
 3\cdot 2^{2/3}t^2\theta^2-2^{4/3}t\theta^3)(2^{-1/3}t^2+\theta^2+2^{4/3}
 t\theta)+\cr
 & +\sin(3\cdot 2^{2/3}t^2\theta^2-2^{4/3}t\theta^3)(\theta^2-2^{-1/3}t^2)]
   \d\theta, 
}$$   
   and the estimate is
uniform with respect to $\xi\in [0,n^{-3/4}]$.
Moreover, if $\xi\in [0,n^{-3/4-\eps}]$ for some $\eps>0$, then  
$$
  p^{(2n)}(2k,0)\asym{n}{\sqrt 2\,e^{n\phi(\xi)}\over
  \Gamma\left({1\over 4}\right)}\,n^{-3/4},
$$
uniformly with respect to $\xi\in [0,n^{-3/4-\eps}]$.
}

\proof
We  perform the change of variable (c):
$$
 \eqalign{
 {\rm\bf (A)} %& ={1\over\pi}\im\int_0^{\alpha}{G(\zt)\over\zt}
% \exp\{n\ovl\Psi_\xi(t)\}[-4e^{-i\pi/4}(v(\xi)+e^{-i\pi/4}t)^3]\d t=\cr
 &= {-4\, e^{n\phi(\xi)}\over\pi\root 4\of{na_4(\xi)}} \im \bigg\{{e^{-i\pi/4}}
 \int_0^{n^{1/4}a_4(\xi)^{1/4}\alpha}
 {G(z(\xi,t_n))\over z(\xi,t_n)}\cdot\cr
 &\qquad\cdot\exp\left\{
 -i\,{n^{1/2}a_2(\xi)\over a_4(\xi)^{1/2}}\,\theta^2
 -\left(  {1\over\sqrt 2}+i{1\over\sqrt 2}  \right) 
 {n^{1/4}a_3(\xi)\over a_4(\xi)^{3/4}}\,\theta^3-
 \theta^4+nR(\xi,t_n)
 \right\}\cdot\cr
 &\qquad\qquad \cdot(v(\xi)+e^{-i\pi/4}t_n)^3\d\theta\bigg\}=\cr
 & {-2\sqrt 2\,e^{n\phi(\xi)}\over\pi\root 4\of{na_4(\xi)}}
  \bigg(\im-\re\bigg) 
 \int_0^{n^{1/4}a_4(\xi)^{1/4}\alpha}
 {G(z(\xi,t_n))\over z(\xi,t_n)}\cdot\cr
 &\qquad\cdot\exp\left\{
 -i\,{n^{1/2}a_2(\xi)\over a_4(\xi)^{1/2}}\,\theta^2
 -\left(  {1\over\sqrt 2}+i{1\over\sqrt 2}  \right) 
 {n^{1/4}a_3(\xi)\over a_4(\xi)^{3/4}}\,\theta^3-
 \theta^4+nR(\xi,t_n)
 \right\}\cdot\cr
 &\qquad\qquad \cdot(v(\xi)+e^{-i\pi/4}t_n)^3\d\theta.
 }\autoeqno{A-x1}
$$
We note that $\arg(v(\xi)+e^{-i\pi/4}t)\in[-\pi/4,0]$ for all $t\ge0$, hence
we will only need expression $\ovl\Psi_{\xi,1}(t)$ for $\ovl\Psi_\xi(t)$.

%%%%%%%%%%%%%%%%%% figura con i quadranti e la semiretta in v
%%%%%SEGMENTI IN V PER BETA=-PI/4
\beginpicture
\setcoordinatesystem units <4mm,4mm> point at 0 0
\setlinear
\put{$\,$} at -4 6.5
\put{$\,$} at -4 -6.5
\put{\it Figure~9: the segment in the} [lt] at -2 0
\put{\it $v$-plane} [lt] at 1 -1

\plot  13 0  24 0 /
\plot  15 -6    15 6 /
\put{0}  at 14.5 -.5
\plot  15 0  21 6 /
\plot  15 0  21 -6 /
\plot 17 0 23 -6 /
\put{$v(\xi)$} [l] at 17 .6

\setdots <2.5pt>
\plot 18 -1 21 -1 /
\put{$v(\xi)+e^{-i{\pi\over 4}}t$} [l] at 21.1 -1
\endpicture

Now we choose $\alpha$ (smaller than $1\over\sqrt 6$) 
such that $|R(\xi,t)|\le{a_4(\xi)t^4\over 2}$
for all $\xi\in[0,a]$ and $t\in[0,\alpha]$.
This choice is possible since  $|R(\xi,t)|\le C|t|^5$ and $a_4(\xi)>c$ for
all $\xi\le a$. Then $|nR(\xi,t_n)|\le{\theta^4\over 2}$ and
we are able to guarantee integrability, since the modulus of the integrand
in equation~\eqref{A-x1} is bounded by
$$
  C\theta\exp\left\{-{n^{1/4}a_3(\xi)\over\sqrt 2\,a_4(\xi)^{3/4}}\,
  \theta^3-{\theta^4\over 2}\right\}\le 
  C\theta\exp\left\{-C^\prime\theta^3-{\theta^4\over 2}\right\}
$$
if $\xi\in[0,n^{-3/4}]$ (we used that $|v(\xi)+e^{-i\pi/4}t_n|\le
C\theta$).

We proceed as in Theorem~\lemmaref{y-xiveloce}: we recall the decomposition
of $G(z)$
$$
\eqalign{
  G(z) & =(1-z)^{-1/4}H(z)+(1-z)^{1/4}K(z)=\cr
	&=(v(\xi)+e^{-i\pi/4}t)^{-1}H(\zt)+(v(\xi)+e^{-i\pi/4}t)K(\zt)
}$$
where we used $z=\zt=1-(v(\xi)+e^{-i\pi/4}t)^{4}$.
Hence the integral in \eqref{A-x1} can be written as:
$$
 \eqalign{
% & \int_0^{n^{1/4}a_4(\xi)^{1/4}\alpha}
% {G(z(\xi,t_n))\over z(\xi,t_n)}\cdot\cr
% & \qquad\cdot\exp\left\{
% -i\,{n^{1/2}a_2(\xi)\over a_4(\xi)^{1/2}}\, \theta^2
% -\left(  {1\over\sqrt 2}+i{1\over\sqrt 2}  \right) 
% {n^{1/4}a_3(\xi)\over a_4(\xi)^{3/4}}\, \theta^3-
% \theta^4+nR(\xi,t_n)
% \right\}\cdot\cr
% &\qquad\qquad \cdot(v(\xi)+e^{-i\pi/4}t_n)^3\d\theta=\cr
 & \int_0^{n^{1/4}a_4(\xi)^{1/4}\alpha}
  \exp\left\{
 -{n^{1/4}a_3(\xi)\over \sqrt 2\,a_4(\xi)^{3/4}}\,\theta^3-
 \theta^4+nR_0(\xi,t_n)
 \right\}\cdot\cr
  & \qquad\qquad\cdot\exp\left\{
 -i\,{n^{1/2}a_2(\xi)\over a_4(\xi)^{1/2}}\,\theta^2
 -i\,{n^{1/4}a_3(\xi)\over \sqrt 2\,a_4(\xi)^{3/4}}\, \theta^3
 +inR_1(\xi,t_n)
 \right\}\cdot\cr
 &\qquad\qquad\cdot{1\over 1-(v(\xi)+e^{-i\pi/4}t_n)^4}\cdot\cr
 & \qquad\qquad \cdot    \left\{ (v(\xi)+e^{-i\pi/4}t_n)^2H(z_n)+
 (v(\xi)+e^{-i\pi/4}t_n)^4K(z_n)\right\}\d\theta
 }\autoeqno{int-asint?}
$$
where $z_n:=z(\xi,t_n)$ and $R_0$ and $R_1$ are respectively the real and
imaginary part of $R$.

We note that both $H_1(z)$ and $K_1(z)$ are $O(|1-z|)$, that is
$O(|v(\xi)+e^{-i\pi/4}t|^{4})=o(v(\xi)^2)+o(t^2)$.
%, since
%$|1-z|=[v(\xi)^2+\sqrt 2\,tv(\xi)+t^2]^2$.
Moreover, it is easy to show that 
$$
{1\over 1-(v(\xi)+e^{-i\pi/4}t_n)^4}=1+O(v(\xi)^4)+O(t_n^4)=
1+o(v(\xi)^2)+o(t_n^2),
$$
and $o(v(\xi)^2)$ and $o(t^2)$ are uniform with respect to 
$\xi\in[0,n^{-3/4}]$.

Then we want to pick the main terms in the following product:
$$
  \eqalign{
  \exp & \left\{
 -i\,{n^{1/2}a_2(\xi)\over a_4(\xi)^{1/2}}\, \theta^2
 -\,i{n^{1/4}a_3(\xi)\over \sqrt 2\,a_4(\xi)^{3/4}}\, \theta^3
 +inR_1(\xi,t_n)
 \right\}\cdot\cr
 &   \cdot(v(\xi)+e^{-i\pi/4}t_n)^2\cdot \left\{H(z_n)+
 (v(\xi)+e^{-i\pi/4}t_n)^2K(z_n)\right\}.
}\autoeqno{prodtrasc}$$
For simplicity's sake we call $iA_n=iA(n,\xi,\theta)$ the function in 
the exponential.
%We note that, since $|v(\xi)+e^{-i\pi/4}t_n|^2 =O(v(\xi)^2)+O(t_n^2)$,
%we can write the product in equation~\eqref{prodtrasc} as
%$$
%  g(\theta,\xi,n)\cdot\{H(z_n)+[O(v(\xi)^2)+O(t_n^2)]\cdot K(z_n)\}
%$$
%where $g(\theta,\xi,n)$ is a suitable function that can easily be deduced
%from \eqref{prodtrasc}. It is very useful to note that the terms involving
%$K(z_n)$ are asymptotically negligible if compared to the ones involving
%$H(z_n)$: for every fixed $\theta$, $z_n$ tends to $1$ as $n$ tends to 
%infinity, $H(z_n)$ and $K(z_n)$ are finite and greater than $C>0$, hence
%$|K(z_n)|/|H(z_n)|\le C$ for some $C>0$.

%Then we develop 
Explicit computation shows that the product in \eqref{prodtrasc} is:
$$
\eqalign{
 & \bigg\{ \cos(A_n)\cdot(v(\xi)^2+\sqrt 2\,v(\xi)t_n)\cdot H_0(z_n)+\cr
	&+\sin(A_n)\cdot(\sqrt 2\,v(\xi)t_n+t_n^2)\cdot H_0(z_n)+o(v(\xi)^2)+
	o(t_n^2)\bigg\}+\cr
  &+i\bigg\{-\cos(A_n)\cdot(\sqrt 2\,v(\xi)t_n+t_n^2)\cdot H_0(z_n)+\cr
	&+\sin(A_n)\cdot(v(\xi)^2\sqrt 2+v(\xi)t_n)\cdot H_0(z_n)+o(v(\xi)^2)+
	o(t_n^2)\bigg\}.
}$$
%Then we rewrite {\bf (A)}:
%$$
%\eqalign{
 %  {-2\sqrt 2\,e^{n\phi(\xi)}\over\pi\root 4\of{na_4(\xi)}}&\im \Bigg\{(1-i)
%   \int_0^{n^{1/4}a_4(\xi)^{1/4}\alpha}\exp\left\{
%   -{n^{1/4}a_3(\xi)\over\sqrt 2\,a_4(\xi)^{3/4}}\,
%   \theta^3-\theta^4+nR_0(\xi,t_n)\right\}\cdot\cr
%   & \cdot(1+o(v(\xi)^2)+o(t_n^2))\cdot\bigg\{
%   \cos(A_n)\cdot(v(\xi)^2+\sqrt 2\,v(\xi)t_n)\cdot H_0(z_n)+\cr
%	&+\sin(A_n)\cdot(\sqrt 2\,v(\xi)t_n+t_n^2)\cdot H_0(z_n)+o(v(\xi)^2)+
%	o(t_n^2)+\cr
 % &+i\big[-\cos(A_n)\cdot(\sqrt 2\,v(\xi)t_n+t_n^2)\cdot H_0(z_n)+\cr
%	&+\sin(A_n)\cdot(v(\xi)^2\sqrt 2+ v(\xi)t_n)\cdot H_0(z_n)+o(v(\xi)^2)+
%	o(t_n^2)\big]\bigg\}\d\theta\Bigg\}.
%}$$
%But the imaginary part of the integral multiplied by $(1-i)$ is simply the 
%integral of the imaginary part of the integrand minus the real one.
The integral {\bf (A)} becomes:
%$$
%\eqalign{
 %  {-2\sqrt 2\,e^{n\phi(\xi)}\over\pi\root 4\of{na_4(\xi)}}&
%   \int_0^{n^{1/4}a_4(\xi)^{1/4}\alpha}\exp\left\{
 %  -{n^{1/4}a_3(\xi)\over\sqrt 2\,a_4(\xi)^{3/4}}\,
%   \theta^3-\theta^4+nR_0(\xi,t_n)\right\}\cdot\cr
 %  & \cdot(1+o(v(\xi)^2)+o(t_n^2))\cdot\bigg\{
 %  -\cos(A_n)\cdot(v(\xi)^2+2\sqrt 2\,v(\xi)t_n+t_n^2)\cdot H_0(z_n)+\cr
%	&+\sin(A_n)\cdot(v(\xi)^2-t_n^2)\cdot H_0(z_n)+o(v(\xi)^2)+
%	o(t_n^2)\bigg\}\d\theta=\cr
 %  ={2\sqrt 2\,e^{n\phi(\xi)}\over\pi\root 4\of{na_4(\xi)}}&
%   \int_0^{n^{1/4}a_4(\xi)^{1/4}\alpha}\exp\left\{
%   -{n^{1/4}a_3(\xi)\over\sqrt 2\,a_4(\xi)^{3/4}}
 %  \theta^3-\theta^4+nR_0(\xi,t_n)\right\}\cdot\cr
%   & \cdot(1+o(v(\xi)^2)+o(t_n^2))\cdot\bigg\{
%   \cos(A_n)\cdot(v(\xi)^2+2\sqrt 2\,v(\xi)t_n+t_n^2)\cdot H_0(z_n)+\cr
%	&+\sin(A_n)\cdot(t_n^2-v(\xi)^2)\cdot H_0(z_n)+o(v(\xi)^2)+
%	o(t_n^2)\bigg\}\d\theta,
%}$$
$$\eqalign{
   {2\sqrt 2\,e^{n\phi(\xi)}\over\pi\root 4\of{a_4(\xi)}\,n^{3/4}}&
   \int_0^{n^{1/4}a_4(\xi)^{1/4}\alpha}\exp\left\{
   -{n^{1/4}a_3(\xi)\over\sqrt 2\,a_4(\xi)^{3/4}}\,
   \theta^3-\theta^4+nR_0(\xi,t_n)\right\}\cdot\cr
   &\cdot(1+o(v(\xi)^2)+o(t_n^2))\cdot\cr
   &\cdot\bigg\{
   \cos(A_n)\cdot(\sqrt n\,v(\xi)^2+2\sqrt 2\sqrt n\,v(\xi)t_n+
   \sqrt n\,t_n^2)\cdot H_0(z_n)+\cr
	& +\sin(A_n)\cdot(\sqrt  n\,t_n^2-\sqrt  n\,v(\xi)^2)\cdot H_0(z_n)+
	\sqrt  n\,o(v(\xi)^2)+\sqrt  n\, o(t_n^2)\bigg\}\d\theta,
	}
\autoeqno{x-int2}$$
where $o(v(\xi)^2)$ and $o(t_n^2)$ are uniform with respect to $\xi\in[0,
n^{-3/4}]$ (observe that we divided it by $n^{1/2}$
 and multiplied the integrand
by the same quantity).

Our aim is to apply Theorem~\lemmaref{dominscon1}, then if 
$f_n(\theta,\xi)$ is the integrand in equation~\eqref{x-int2}, we have 
to evaluate its uniform
asymptotic $h_n(\theta,\xi)$, show that there exist $g(\theta)$ and
$g_1(\theta)$ in ${\rm L}^1(\reali^+)$ such that
$$
  |h_n(\theta,\xi)|\le g(\theta)\hbox{ and }\quad\qquad
  |f_n(\theta,\xi)|\le g_1(\theta),
$$
for every $\theta$ and $\xi\le n^{-3/4}$ and that there exists $c>0$ such that
$$
	\left|\int_{\reali^+}h_n(\theta,\xi)\d\theta\right|\ge c
\autoeqno{intnonzero}
$$
for all $n$ sufficiently large and $\xi\le n^{-3/4}$.
Then we will get that 
$$
  \int_{\reali^+}f_n(\theta,\xi)\d\theta\ \sim\ 
  \int_{\reali^+}h_n(\theta,\xi)\d\theta
$$
uniformly with respect to $\xi\in[0,n^{-3/4}]$.

Now recalling the definitions of $t_n$, $A_n$ and $z_n$, the asymptotic
values of $a_i(\xi)$ and $v(\xi)$, and 
that $H(z_n)$ tends to $\sqrt 2$ as $n$
tends to infinity,
$$
\eqalign{
   f_n(\theta,\xi) 
%=& \ident_{[0,n^{1/4}a_4(\xi)^{1/4}\alpha]}(\theta)
 %  \exp\left\{
%   -{n^{1/4}a_3(\xi)\over\sqrt 2\,a_4(\xi)^{3/4}}
 %  \theta^3-\theta^4+nR_0(\xi,t_n)\right\}\cdot\cr
%   & \cdot(1+o(v(\xi)^2)+o(t_n^2))\cdot\bigg\{
%   \cos(A_n)\cdot(\sqrt  n\,v(\xi)^2+2\sqrt 2\sqrt  n\,v(\xi)t_n+
%	\sqrt n\,t_n^2)\cdot  \cr
%	&\cdot H_0(z_n)
%	+\sin(A_n)\cdot(\sqrt n\,t_n^2-\sqrt n\,v(\xi)^2)\cdot H_0(z_n)+
%	\sqrt n\, o(v(\xi)^2)+
%	\sqrt n\, o(t_n^2)\bigg\}\cr
   \sim h_n(\theta,\xi)= &
	   \sqrt 2\exp\left\{
	-{2^{4/3}n^{1/4}\xi^{1/3}}\,
	\theta^3-\theta^4\right\}\cdot\cr
	&\cdot\bigg\{
	\cos\left(3\cdot 2^{2/3}n^{1/2}\xi^{2/3}\theta^2-
	{2^{4/3}n^{1/4}\xi^{1/3}}\, \theta^3\right)\cdot\cr
	&\cdot\left(2^{-1/3} n^{1/2}\xi^{2/3}+2^{4/3}n^{1/4}\xi^{1/3}
	{\theta}
	+\theta^2\right)+\cr
	&+\sin\left(3\cdot 2^{2/3}n^{1/2}\xi^{2/3}\theta^2-
	{2^{4/3}n^{1/4}\xi^{1/3}}\, \theta^3\right)
	\cdot({\theta^2}-2^{-1/3}n^{1/2}\xi^{2/3})\bigg\}.
}$$
From these formulas it is easy to see that we can use
$$
\eqalign{
& g_1(\theta) = C\exp\left\{-C\theta^3-{\theta^4\over 2}\right\}
	(1+\theta+\theta^2);\cr
& g(\theta) = C\exp\{-C\theta^3-{\theta^4}\}   
	(1+\theta+\theta^2).
}$$
Now we put $\delta:=n^{1/4}\xi^{1/3}$:
clearly $\delta\in[0,1]$ and assumption~\eqref{intnonzero}
is equivalent to $I(\delta)>0$ for every $\delta\in[0,1]$.
But by numerical computation performed with {\it Derive for Windows}
this appears obvious.

Hence by Theorem~\lemmaref{dominscon1}
$$
{\rm\bf (A)}\asym{n}
  {4e^{n\phi(\xi)}\over
  \pi }n^{-3/4}I(n^{1/4}\xi^{1/3}),
$$
and this is the uniform asymptotic estimate for the $p^{(2n)}(2k,0)$
since by Lemma~\lemmaref{int-2-x1}
$$
\left|{ {\rm\bf(B)}\over{\rm\bf(A)} }\right|\le C{e^{n\phi(\xi)}\lambda^n
	\over e^{n\phi(\xi)}n^{3/4}\min_{\delta\in[0,1]}I(\delta)}\le
	C\lambda^n n^{-3/4}
$$
which tends to zero uniformly with respect to $\xi\in[0,n^{-3/4}]$.
Finally we observe that the same argument conducted this far shows that,
if $\xi\in[0,n^{-3/4-\eps}]$ for some $\eps>0$, the uniform asymptotic
estimate is
$$
p^{(2n)}(2k,0)\asym{n}
  {4e^{n\phi(\xi)}\over
  \pi }n^{-3/4}I(0),
$$
and the thesis follows from $I(0)=\int_0^{+\infty}\theta^2e^{-\theta^4}
\d\theta={\pi\sqrt 2\over 4 \Gamma\left({1\over 4}\right)}$. 
We observe that this last result can be obtained under the weaker 
assumption that
$n^{1/4}\xi^{1/3}\ra 0$ uniformly with respect to $\xi$.
\QED
%%%%%%%%%%%%%%%%%%%%%%%%%XBRUTTO.TEX%%%%%%%%%%%%%%%%%%%%%%%%%%%%%%%%%%%%%
\autosez{xi>n3/4}{Estimates along the $x$-axis: the case 
$\xi\in [n^{-{3\over 4}},a]$}

We observe that since in this case we have to use substitution (a), the
curve of integration we used so far does not fit. In fact in the 
exponential, $\theta^2$ would have a pure imaginary coefficient.
Hence we are forced to seek another solution, that will be choosing
a curve of integration made of more pieces: the curve
$z_1(\xi,t)=1-(v(\xi)+e^{-i{\pi\over 3}}t)^4$ for $0\le t\le t_o$,
$z_2(\xi,t)=1-(v(\xi)+e^{i{\pi\over 3}}t)^4$ for $t_o\le t\le \alpha$,
the conjugate of these two curves and a circular arc to make the whole
connected.
Here is how these curves appear.

%%%%%%%%%%%%%%% CURVA DI INTEGRAZIONE PER t<alpha
%%%%%%%%%%%%%%% CASO xi>n^{-3/4}
\beginpicture
\setcoordinatesystem units <5mm,5mm> point at 0 0
\setlinear
\put{$\,$} at -2 4.5 
\put{$\,$} at -2 -3
\put{\it Figure~10: the curve for} [lt] at 0 0
\put{$t\le\alpha$} [lt] at 3 -1
\plot  9 0  26 0 /
\plot  11 -2.5    11 4 /
\put{$0$}  at 10.5 -.5
\put{$z(\xi,t_o)$} at 22.90456 -.5
\put{$z(\xi,\alpha)$} [r] at 25.22923 3.050731 
\setquadratic
\plot   13.75768 0              13.72046 .0222696       13.69853 .04963793
	13.69639 .08198807      13.71849 .1191252       13.76898 .1607764
	13.85151 .2065911       13.96894 .2561405       14.12311 .308918
	14.31473 .3643391       14.54332 .4217412       14.80728 .4803839
	15.10404 .5394488       15.43028 .5980395       15.78222 .6551818
	16.15586 .7098234       16.5472 .7608343        16.95241 .8070062
	17.36796 .8470531       17.79067 .8796111       18.2178 .9032381
	18.64697 .9164144       19.07621 .9175422       19.50387 .9049456
	19.92863 .876871        20.34944 .8314868       20.76545 .7668833
	21.17604 .6810731       21.58073 .5719907       21.97918 .4374927
	22.37115 .2753577       22.75651 .0832864       22.90456 0
	23.2806 .2376096        23.64994 .5109637       24.01261 .8226258
	24.36868 1.175237       24.71826 1.571517       25.06145 2.014263
	25.39839 2.50635        25.72923 3.050731 /
%%%%%%% freccette
\setlinear
% frecce orizzontali : xshift=.5 yshift=.25
\plot  19.00387 .6549456       19.50387 .9049456 /
\plot  19.00387 1.1549456      19.50387 .9049456 /

% frecce su una retta a 45 gradi: se (x0,y0)=punta
%       e (x1,y1)=coda bassa, (x2,y2)=coda alta
%       allora  x1=x0-(.25)*1/sqrt(2)  ;y1=y0-(.75)*1/sqrt(2)
%               x2=x0-(.75)*1/sqrt(2)  ;y2=y0-(.25)*1/sqrt(2)
\plot   24.54148 1.041186    24.71826 1.571517 /
\plot   24.18793 1.394740    24.71826 1.571517 /

%%%%%%% stanghetta per 1
\setlinear
\linethickness=.3pt
\plot 14 -.1  14 .1 /
\put{1}  at 14 -.5
\endpicture

And here is a sketch for the curve of integration: we substituted
the circumference with an ellipsis (not centered in the origin), only for
aesthetical reasons.

\beginpicture
\setcoordinatesystem units <4mm,4mm> point at 0 0
\setlinear
\put{$\,$} at 0 4.5
\plot  2 0  26 0 /
\plot  11 -4    11 4 /
\put{$0$}  at 10.5 -.5
\setquadratic
\plot   25.39839 -2.50635       25.06145 -2.014263      24.71826 -1.571517
	24.36868 -1.175237      24.01261 -.8226258      23.64994 -.5109637
	23.2806 -.2376096       22.90456 0              22.75651 -.0832864
	22.37115 -.2753577      21.97918 -.4374927      21.58073 -.5719907
	21.17604 -.6810731      20.76545 -.7668833      20.34944 -.8314868
	19.92863 -.876871       19.50387 -.9049456      19.07621 -.9175422      
	18.64697 -.9164144      18.2178 -.9032381       17.79067 -.8796111
	17.36796 -.8470531      16.95241 -.8070062      16.5472 -.7608343
	16.15586 -.7098234      15.78222 -.6551818      15.43028 -.5980395
	15.10404 -.5394488      14.80728 -.4803839      14.54332 -.4217412
	14.31473 -.3643391      14.12311 -.308918       13.96894 -.2561405      
	13.85151 -.2065911      13.76898 -.1607764      13.71849 -.1191252      
	13.69639 -.08198807     13.69853 -.04963793     13.72046 -.0222696      
	13.75768 0              
	13.72046 .0222696       13.69853 .04963793
	13.69639 .08198807      13.71849 .1191252       13.76898 .1607764
	13.85151 .2065911       13.96894 .2561405       14.12311 .308918
	14.31473 .3643391       14.54332 .4217412       14.80728 .4803839
	15.10404 .5394488       15.43028 .5980395       15.78222 .6551818
	16.15586 .7098234       16.5472 .7608343        16.95241 .8070062
	17.36796 .8470531       17.79067 .8796111       18.2178 .9032381
	18.64697 .9164144       19.07621 .9175422       19.50387 .9049456
	19.92863 .876871        20.34944 .8314868       20.76545 .7668833
	21.17604 .6810731       21.58073 .5719907       21.97918 .4374927
	22.37115 .2753577       22.75651 .0832864       22.90456 0
	23.2806 .2376096        23.64994 .5109637       24.01261 .8226258
	24.36868 1.175237       24.71826 1.571517       25.06145 2.014263
	25.39839 2.50635 /
%%%%%%% freccette
\setlinear
% frecce orizzontali : xshift=.5 yshift=.25
\plot  19.00387 .6549456       19.50387 .9049456 /
\plot  19.00387 1.1549456      19.50387 .9049456 /

\plot  19.57621 -1.1675422       19.07621 -.9175422 /
\plot  19.57621 -.6675422       19.07621 -.9175422 /

% frecce su una retta a 45 gradi: se (x0,y0)=punta
%       e (x1,y1)=coda bassa, (x2,y2)=coda alta
%       allora  x1=x0-(.25)*1/sqrt(2)  ;y1=y0-(.75)*1/sqrt(2)
%               x2=x0-(.75)*1/sqrt(2)  ;y2=y0-(.25)*1/sqrt(2)
\plot   24.54148 1.041186    24.71826 1.571517 /
\plot   24.18793 1.394740    24.71826 1.571517 /

\plot   24.89901 -1.352013  24.36868 -1.175237 /
\plot   24.54545 -1.705567  24.36868 -1.175237 /
%%%%%%% stanghetta per 1
\setlinear
\linethickness=.3pt
\plot 14 -.1  14 .1 /
\put{1}  at 14 -.5

\ellipticalarc axes ratio 4:1 266 degrees from 25.39839 2.50635 
	center at 16 0 

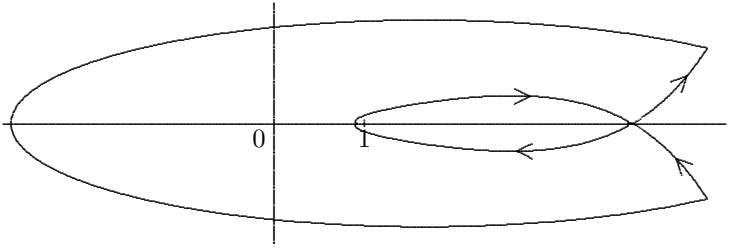
\put{\it Figure~11: The curve of integration in the case 
	$\xi\ge n^{-{3\over 4}}$} 
	[lt] at 4 -5.5
\put{$\,$} at 0 -7.5
\endpicture

%%%%
We let $z(\xi,t)$ be the expression for the curve of integration
for $-\alpha\le t\le\alpha$, that is
$$
z(\xi,t)=\cases{
    1-(v(\xi)+e^{-i{\pi\over 3}}t)^4 & if $0\le t\le t_o$,\cr
    1-(v(\xi)+e^{i{\pi\over 3}}t)^4  & if $t_o< t\le \alpha$,\cr
    1-(v(\xi)-e^{i{\pi\over 3}}t)^4 & if $-t_o\le t\le 0$,\cr
    1-(v(\xi)-e^{-i{\pi\over 3}}t)^4  & if $-\alpha\le t< -t_o$.
}$$
The transition probabilities are given now by equation \eqref{xab} where
$\gamma_1$ is this new $z(\xi,t)$ and $\gamma_2$ is the corresponding
circular arc.
As we already remarked, integral {\bf (A)} is equal to
$$
  {1\over\pi}\im\int_0^{\alpha}
  {G(\zt)\over \zt}\exp\{n\ovl\Psi_\xi(t)\}
  [-4\,e^{i\beta(t)}(v(\xi)+e^{i\beta(t)}t)^3]\d t,
\autoeqno{int-circ+baf}$$
where $\beta(t)=-\pi/3$ if $t\in [0,t_0]$, $\beta(t)=\pi/3$ 
if $t\in (t_0,\alpha]$ and $t_0=2/(\sqrt 3 -1)v(\xi)$.

Now the explicit expression for $\ovl\Psi_\xi(t)$ depends on where in the 
complex plane the curves $v(\xi,t)=v(\xi)+ e^{\pm i{\pi\over 3}}t$ lie.

%%%%%SEGMENTI IN V PER BETA=PI/3
\beginpicture
\setcoordinatesystem units <4mm,4mm> point at 0 0
\setlinear
\put{$\,$} at -2 7.5
\put{$\,$} at 0 -8
\put{\it Figure~12: the segments in the} [lt] at 0 0
\put{\it $v$-plane} [lt] at 3 -1

\plot  13 0  24 0 /
\plot  15 -7    15 7 /
\put{0}  at 14.5 -.5
\plot  15 0  22 7 /
\plot  15 0  22 -7 /

\plot  17.36603 2.36603  20.04145 7 /
\plot 16 0  17.36603 -2.36603 /
\put{$v(\xi)$} [l] at 17 .6

\setdashes <3pt>
\plot  16 0   17.36603 2.36603   /
\plot   17.36603 -2.36603  20.04145 -7 /

\setdots <2.5pt>
\plot 16 0 17 .6 /
\plot  16.78 -1.183015  20 -1.183015 /
\put{$v(\xi)+e^{-i{\pi\over 3}}t$} [l] at 20.1 -1.183015
\plot 17.46603 -2.36603  20 -2.36603 /
\put{$v(\xi)+e^{-i{\pi\over 3}}t_o$} [l] at  20.1 -2.36603
\plot  18.12 3.5        20 3.5 /
\put{$v(\xi)+e^{i{\pi\over 3}}t$} [l] at 20.1 3.5
\plot 17.46603 2.36603  20 2.36603 /
\put{$v(\xi)+e^{i{\pi\over 3}}t_o$} [l] at  20.1 2.36603

\endpicture
From the picture above we can see that if $t\in[0,t_o]$ then $v(\xi)+
e^{i{\pi\over 3}}t$ has argument in the interval $[-\pi/4,0]$, while
if $t\ge t_o$ then $v(\xi)+
e^{-i{\pi\over 3}}t$ has argument in the interval $[\pi/4,\pi/2)$. 
Hence in equation~\eqref{int-circ+baf} we can substitute $\ovl\Psi_\xi(t)$
with the expansion of $\ovl\Psi_{\xi,1}(t)$ if $t\in[0,t_o]$ and with
the expansion of $\ovl\Psi_{\xi,2}(t)$ otherwise (see 
Lemma~\lemmaref{psitx-sv}).
Then our goal is the estimate of the imaginary part of
$$
  \eqalign{
  -4\int_0^{\alpha} & {G(z_1(\xi,t))\over z_1(\xi,t)}e^{-i{\pi\over 3}}\,
	(v(\xi)+e^{-i{\pi\over 3}}t)^3
	\ident_{[0,t_o]}(t)\cdot\cr
	&\cdot\exp\left\{n\phi(\xi)+ne^{-{2\over 3}\pi i}a_2(\xi)t^2-na_3(\xi)t^3+
	ne^{-{4\over 3}\pi i}a_4(\xi)t^4 +nR(\xi,t)\right\}+\cr
	& + {G(z_2(\xi,t))\over z_2(\xi,t)}e^{i{\pi\over 3}}\,
	(v(\xi)+e^{i{\pi\over 3}}t)^3
	\ident_{[t_o,\alpha]}(t)\cdot\cr
	&\cdot\exp\left\{n\phi(\xi)+ne^{{2\over 3}\pi i}a_2^\prime(\xi)t^2-
	na_3^\prime(\xi)t^3+
	ne^{{4\over 3}\pi i}a_4^\prime(\xi)t^4 +nR^\prime(\xi,t)\right\}
	\d t,
}$$
where we wrote $R^\prime(\xi,t)$ for the second remainder term to
distinguish it from the first one.

We perform the change of variable $\theta=\sqrt{na_2(\xi)}\,t$ and
write $t_n$ for ${\theta\over\sqrt{na_2(\xi)}}$. 

We choose $\alpha$ such that $|R^\prime(\xi,t)|\le
{\re(a_4^\prime(\xi))\over 4}\, t^4$ and $|R(\xi,t)|\le
{a_4(\xi)\over 4}\, t^4$ for all $\xi\le a$ and $|t|\le\alpha$.
Moreover we require that $\alpha>\max_\xi t_o={2\over\sqrt 3-1}v(a)$,
which is surely true if $a$ is sufficiently small.

We distinguish two subcases: $\xi\in[n^{-3/4},n^{-3/4+\eps}]$ and $\xi\in
[n^{-3/4+\eps},a]$.

\theorem{xi>n3/4-1}{If $\xi\in[n^{-3/4+\eps},a]$, for some $\eps>0$ and
for sufficiently small $a$, then
$$
p^{(2n)}(2k,0)\asym{n}{4\,e^{n\phi(\xi)}v(\xi)^3
\over\pi\sqrt{a_2(\xi)}}{G(\zo)\over\zo}\,n^{-1/2}\int_0^{+\infty}e^{-{\theta^2
\over 2}}\cos\left({\pi\over 6}-{\sqrt3\over 2}\theta^2\right)\d\theta
$$ 
uniformly with respect to $\xi\in[n^{-3/4+\eps},a]$.}

\proof
We apply Theorem~\lemmaref{dominunif}, observing that the
modulus of the integrand is 
majorized by
$$
C(1+\theta)^3 \exp\left\{{\theta^2\over 2}-C\theta^3-
      C\theta^4\right\},
$$
for every $\theta\ge0$, $\xi\in[n^{-3/4+\eps},a]$ and $n$.
Moreover, the integrand converges pointwise, uniformly with respect to
$\xi\in[n^{-3/4+\eps},a]$, to
$$
	e^{-i{\pi\over 3}}\exp\left\{e^{-{2\over 3}\pi i}\theta^2\right\}.
$$
Negligibility of {\bf (B)} is shown as usual.
\QED

\theorem{xi>n3/4-2}{If $\xi\in[n^{-3/4},n^{-3/4+\eps}]$, for some 
$\eps>0$, then
$$
p^{(2n)}(2k,0)\, {\buildrel n\over\sim}\,
 {2^{1/6} e^{n\phi(\xi)}\over \pi\sqrt3 }
	{G(\zo)\over\zo}\,I(n^{-1/2}\xi^{-2/3})\,\xi^{2/3}n^{-1/2},
$$ 
uniformly with respect to $\xi\in[n^{-3/4},n^{-3/4+\eps}]$,
where $I(t)$ is an integral function, defined as follows:
$$\eqalign{
  I(t)& :=\int_0^{+\infty} 
	\exp\left\{-{\theta^2\over 2}-{2^{11/6}t\over
	6\sqrt 3 }\,\theta^3
      -{t^2\over 9\cdot 2^{7/3}}\,
	\theta^4        \right\}\cdot\cr
	&\cdot\bigg\{
	\sqrt3\,
	\cos\left(-{\sqrt3\,\theta^2\over 2}+{\sqrt3\,t^2\over 
	9\cdot 2^{7/3}}\, \theta^4\right)
	\cdot (1+\sqrt3\,2^{-1/6}t\theta
	-(3\sqrt6)^{-1}t^{3}\theta^3)                        +\cr
	 & -\sin\left(-{\sqrt3\,\theta^2\over 2}+{\sqrt3\,t^2\over 
	9\cdot 2^{7/3}}\, \theta^4\right)
	 \cdot (1-\sqrt3\,2^{-1/6}t\theta
	-2^{2/3}t^{2}\theta^2
	-(3\sqrt6)^{-1}t^{3}\theta^3)\bigg\}\cdot\cr  
	&\cdot\bigg\{ \ident_{[0,{2^{7/6}\over\sqrt 3-1}
	{\sqrt3\over t}]}(\theta)       -
	\ident_{[{2^{7/6}\over\sqrt 3-1}
	{\sqrt3\over t},+\infty)}(\theta)
	\bigg\}.
}$$
}

\proof
The substantial difference with the previous case is that 
we cannot find a pointwise limit for the integrand, but only
a pointwise asymptotic estimate. Hence we apply 
Theorem~\lemmaref{dominscon1}, exactly as in the proof of 
Theorem~\lemmaref{ximin}.
\QED

%% file: fr.tex
\autosez{fr}{Final remarks}

First, we want to spend a few words to explain how one can obtain an
estimate for the transition probabilities with odd time and space parameters.

From the results on the $y$-axis we obtain also
$p^{(2n+1)}((0,2k+1),(0,0))$, since 
$$
  p^{(2n+1)}((0,2k+1),(0,0)) ={1\over 2}\{ p^{(2n)}(0,2k+2),(0,0))+
  p^{(2n)}(0,2k),(0,0))\},
$$
as we noted in Paragraph~\sref{y}.

The estimate
of $p^{(2n+1)}((2k+1,0),(0,0))$ requires further calculation, but
it is not much different from what we did this far.
We briefly point out what has to be done.

Let $G(z)$, $F(z)$ and $\Psi_\xi(z)$ be as in Paragraph~\sref{x},
then
$$
\eqalign{
 p^{(2n+1)}((2k+1,0)& ,(0,0)) ={1\over 2\pi i}\oint
 {G(z)\over z}{F(z)\over\sqrt z}\cdot\exp\{n\Psi_\xi(z)\}\d z=\cr
         & ={1\over 2\pi i}\oint
 {G(z)\over z}{1+\sqrt{1-z}-\sqrt2\sqrt{1-z+\sqrt{1-z}}\over z}
 \cdot\exp\{n\Psi_\xi(z)\}\d z.
 }$$
Thus, if $\xi\in[a,1-c]$ the proof will differ from that of 
Theorem~\lemmaref{x-xinonzero} because one should multiply and divide
the integrand not by $G(z_o(\xi))\over z_o(\xi)$ but by 
$G(z_o(\xi))\cdot F(z_o(\xi))\over z_o(\xi)^2$. The rest of the proof is
completely analogous to that of Theorem~\lemmaref{x-xinonzero}.

The same observation can be made in the case $\xi\ge n^{-3/4}$.

In the case $\xi\in[0,n^{-3/4}]$ one will need a decomposition (similar to
the decomposition \eqref{decG} of $G(z)$) of $F(z)\cdot\sqrt z$:
$$
1+\sqrt{1-z}-\sqrt2\sqrt{1-z+\sqrt{1-z}}=1+\sqrt{1-z}-
        \root 4\of{1-z}J(z)+\root 4\of{(1-z)^3}L(z),
$$
where $J(z)$ and $L(z)$ are two holomorphic functions. Then
$$
  \eqalign{
  G(z)& \cdot(1+\sqrt{1-z}-\sqrt2\sqrt{1-z+\sqrt{1-z}})=\cr
  =& (1-z)^{-1/4}H(z)+(1-z)^{1/4}(H(z)+K(z))+(1-z)^{1/2}(J(z)\,K(z)+
  H(z)\,L(z))+\cr
  &+(1-z)^{3/4}K(z)+\hbox{\it holomorphic functions},
}$$
and one can proceed as in Theorem~\lemmaref{x-xinonzero}.
\medskip
Moreover, by symmetry, translation invariance and reversibility
from our estimates of $p^{(n)}((k,0),(0,0))$ and $p^{(n)}((0,k),(0,0))$
(for $k>0$) one easily derives uniform estimates for 
$p^{(n)}((k,0),(k_1,0))$, $p^{(n)}((k_1,k),(k_1,0))$, $p^{(n)}((k_1,0),(k_1,k))$
for $k,k_1\in\zzz$.

We note that it is not possible to extend these estimates to uniform
estimates for any starting and ending point, as we did in equation{loc-lim2}
for local estimates since the asymptotic we use (picked from \cite{Ber-Zuc}
Paragraph 6) does not hold uniformly.

%%%%%%%%%%%%%%%%%%%%%%%%%%%%%%%%%%%%%%%%%%%%%%%%%%%%%%%%%%%%%%%%%%%%%%%%%%
We now remark some particular features of the transition probabilities
of the simple random walk on $C_2$ that appear after our calculations.
We recall our estimates 
for $\xi$ tending very fast to zero, and we
stress the dependence on $\xi$.
If $k/n=\xi\in [0,n^{-3/4-\eps}]$ for some $\eps>0$, then
$$
  p^{(2n)}(x_k,o)\asym{n}
	\cases{{\sqrt2\,e^{-n\xi^2}\over
  \Gamma\left({1\over 4}\right)}\,n^{-3/4}, & if $x_k=(0,2k)$\cr
{\sqrt 2\,e^{-Cn\xi^{4/3}}\over
  \Gamma\left({1\over 4}\right)}\,n^{-3/4}, & if $x_k=(2k,0)$\cr
}
\autoeqno{xiveloce}$$
 uniformly with respect to $\xi\in [0,n^{-3/4-\eps}]$.

It is worth noting that, even if both these results agree with
the local limit estimate (substitute $\xi=0$), they show a different dependence on
$\xi$ along the two axes.
In particular this reflects in the impossibility of finding for $C_2$ an estimate
of the type Jones found for the $2$-dimensional Sierpi\'nski graph (\cite{Jo}),
and Barlow and Bass for the graphical Sierpi\'nski carpet (\cite{Ba-Ba2}).
Recall that associated with a graph with polynomial growth there are three 
positive constants: the {\it spectral dimension} $\delta_s$, the {\it
fractal dimension} $\delta_f$ and the {\it walk dimension} $\delta_w$.
In typical cases, one has the so-called ``Einstein relation'':
$\delta_s\cdot\delta_w=2\delta_f$ (see Telcs \cite{Te1}, \cite{Te2}).
For $C_2$, $\delta_s=3/2$ (see, for instance Cassi and Regina
\cite{Cassi-Regina}) and $\delta_f=2$. From the Einstein relation one would
therefore expect $\delta_w=8/3$. 
Here is Jones' estimate:
for some positive constants
$c_1,c_2,c_3,c_4$ for sufficiently large $n$ and
for all $x,y$ in the two-dimensional
Sierpi\'nski graph,
$$
c_1\,n^{-\delta_s/2}\exp\left(c_2\,n\xi^{\delta_w/(\delta_w-1)}
        \right)\le p^{(n)}(x,y)\le
c_3\,n^{-\delta_s/2}\exp\left(c_4\,n\xi^{\delta_w/(\delta_w-1)}
        \right),
$$
where $p^{(n)}(x,y)$ are the transition probabilities of the simple random
walk and \break $\xi=d(x,y)/n$.

Now, equation~\eqref{xiveloce} shows that if $x=(0,2k)$ and
$y=(0,0)$ one would get a ``Jones''-type estimate with $\delta_w=2$,
while if $x=(0,2k)$ and $y=(0,0)$ the estimate would do if $\delta_w=4$.
Observe that equation~\eqref{xiveloce} implies that $p^{(2n)}((2k,0),(0,0))/
p^{(2n)}((0,2k),(0,0))$ is not bounded below by a positive constant, as 
Jones' estimate would imply.

We finally observe that the method we used here to provide a uniform
asymptotic estimate of transition probabilities has already been
employed for homogeneous trees and free groups (see Lalley~\cite{Lal1},
\cite{Lal2} and \cite{Woess-giallo}). One of the aims of this paper was
also shedding new light onto this method. It seems that one could
use the described technique for more general graphs and random walks
(provided the knowledge of the Green function). As we noted, the major
difficulty appears to be in the choice of a proper contour of integration,
for which no recipe is known (we investigated only a restricted family
of curves in our case).

A natural extension of our results could be the analogues for $d$-combs.
However, technical difficulties increase notably already for $d=3$, since
$$
  G_3(z)={3\over\sqrt{3(1-z^2)+2\sqrt{1-z^2}+2\sqrt2\sqrt{1-z^2+\sqrt{1-z^2}}}}.
$$